%% file: LAA18ddinv.tex
\documentclass{siamltex}

\input preamble.tex

\begin{document}

\maketitle

\input abstract.tex
\input intro.tex
\input problem.tex

\input fwd_prec.tex

\input inv_prec.tex
\input conclusions.tex
\section*{Acknowledgments}
This material is based upon work supported by NSF awards CCF-1817048 and CCF-1725743; by NIH award 5R01NS042645-14;  by the U.S. Department of Energy, Office of Science, Office of Advanced Scientific Computing Research, Applied Mathematics program under Award Number DE-SC0019393; and by the U.S. Air Force Office of Scientific Research award FA9550-17-1-0190.  Any opinions, findings, and conclusions or recommendations expressed herein are those of the authors and do not necessarily reflect the views of the AFOSR,  DOE, NIH, and NSF. Computing time on the Texas Advanced Computing Centers Stampede system was provided by an allocation from TACC and the NSF. The authors would also like to thank Marcus Sarkis for several useful conversations.


\input appendix.tex

\bibliography{./Biblionew}

\end{document}

%% file: preamble.tex
\usepackage[
  hyperindex  = {true},
  colorlinks  = {true},
  linkcolor   = {blue},
  citecolor   = {blue}
]{hyperref}
\usepackage{multirow}
\usepackage{color, colortbl}
\usepackage{caption}
\usepackage{subcaption}
\usepackage{latexsym,graphicx, amssymb, amsfonts, setspace, geometry, amsmath}
\usepackage{pgfplots}
\pgfplotsset{compat=1.8}
\usepgfplotslibrary{statistics}
\usepackage{grffile}
\usepackage{pgfplotstable}
\usepackage{slashbox}
\usepackage{mathtools}
\usepackage[utf8]{inputenc}
\usepackage[english]{babel}
\usepackage{float}
\usepackage{esvect}
\usepackage{algorithm,algorithmic}
\usepackage{enumitem}
\usepackage{hhline}
\usepackage{varwidth}
\usepackage{tikz}
\usetikzlibrary{backgrounds,patterns,calc}
\usepackage{xparse}

\pgfplotsset{compat=newest}

\makeatletter
\newcommand{\doublewidetilde}[1]{{%
  \mathpalette\double@widetilde{#1}%
}}
\newcommand{\double@widetilde}[2]{%
  \sbox\z@{$\m@th#1\widetilde{#2}$}%
  \ht\z@=.9\ht\z@
  \widetilde{\box\z@}%
}
\makeatother

\newcommand{\xb}{\mathbf{x}}
\newcommand{\ub}{\mathbf{u}}

\newcommand{\yb}{\mathbf{y}}

\newcommand{\db}{\mathbf{d}}

\newcommand{\Ib}{\mathbf{I}}
\newcommand{\bb}{\mathbf{b}}
\newcommand{\Ab}{\mathbf{A}}

\newcommand{\Fb}{\mathbf{F}}
\newcommand{\Gb}{\mathbf{G}}
\newcommand{\Hb}{\mathbf{H}}
\newcommand{\Jb}{\mathbf{J}}
\newcommand{\Qb}{\mathbf{Q}}
\newcommand{\Rb}{\mathbf{R}}
\newcommand{\Sb}{\mathbf{S}}

\newcommand{\Ub}{\mathbf{U}}
\newcommand{\Vb}{\mathbf{V}}

\newcommand{\parB}{{\partial \mathcal{B}}}

\newcommand{\tran}{^{\top\kern-\scriptspace}}

\DeclareMathOperator*{\argmin}{arg\,min}

\definecolor{Gray}{gray}{0.9}
\newcolumntype{g}{>{\columncolor{Gray}}c}

\pgfmathdeclarefunction{gauss}{3}{%
  \pgfmathparse{1/(#3*sqrt(2*pi))*exp(-((#1-#2)^2)/(2*#3^2))}%
}

\title{A domain decomposition preconditioning for the integral equation formulation of the inverse scattering problem}

\author{Carlos Borges\thanks{Department of Mathematics, University of Central Florida, Orlando, FL}
\and
George Biros\thanks{Department of Mechanical Engineering and Institute for Computation Engineering and Sciences, University of Texas, Austin, TX}
}

%% file: abstract.tex
\begin{abstract}
We propose domain decomposition preconditioners for the solution of an integral equation formulation of forward and inverse acoustic scattering problems with point scatterers. We study both forward and inverse problems and propose preconditioning techniques to accelerate the iterative solvers. For the forward scattering problem, we extend the domain decomposition based preconditioning techniques presented for partial differential equations in {\em ``A restricted additive Schwarz preconditioner for general sparse linear systems", SIAM Journal on Scientific Computing, 21 (1999), pp. 792--797}, to integral equations.  We combine this domain decomposition preconditioner with a low-rank correction, which is easy to construct, forming a new preconditioner. For the inverse scattering problem, we use the forward problem preconditioner as a building block for constructing a preconditioner for the Gauss-Newton Hessian. We present numerical results that demonstrate the performance of both preconditioning strategies.
\end{abstract}

%% file: intro.tex
\section{Introduction}\label{s:intro}
 
\begin{figure}[H]
\begin{subfigure}[t]{0.48\textwidth}
\center
\begin{tikzpicture}[scale=0.50]

\draw[ultra thick,->] (-3-1,2.5) -- (-1-1,2.5);
\draw (-2.5-1,1.5) -- (-2.5-1,3.5);
\draw (-2-1,1.5) -- (-2-1,3.5);
\draw (-1.5-1,1.5) -- (-1.5-1,3.5);

\node at (0,-.1) { };
\node at (0,5.1) { };


\draw[] (2.5,2.5) circle [radius=4];

\draw[black,dotted,line width=3pt, line cap=round, dash pattern=on 0pt off 6\pgflinewidth] (2.5,2.5) circle [radius=4];


\begin{scope}
\pgfsetfillpattern{dots}{red}
  \filldraw  (1, 1)
  .. controls ++(165:-1) and ++(240: 1) .. ( 4, 1)
  .. controls ++(240:-1) and ++(165:-1) .. ( 4, 4)
  .. controls ++(165: 1) and ++(175:-2) .. (2, 4.5)
  .. controls ++(175: 2) and ++(165: 1) .. ( 1, 1);
\end{scope}

\node at (-2-1,1) {$u^{\emph{inc}}$};
\node at (8,1) {$u^{\emph{scat}}=?$};
\node at (2.5,2.5) {\bf $\text{supp(q)}$};
\node at (6,-0.5) {$\parB$};
\node at (7,6.5) {receivers};

\draw[ultra thick,->] (6+1,2) -- (7.75+1,2.75);
\draw (6.5+1,1.25) -- (6+1,3);
\draw (7+1,1.5) -- (6.5+1,3.25);
\draw (7.5+1,1.75) -- (7+1,3.5);

\draw[thick,->] (7,6.2) -- (6.1,4.6);

\end{tikzpicture}

\caption{Forward scattering problem}\label{fig:fwd_prob}
\end{subfigure}
\begin{subfigure}[t]{0.48\textwidth}
\center
\begin{tikzpicture}[scale=0.50]

\draw[ultra thick,->] (-3-1,2.5) -- (-1-1,2.5);
\draw (-2.5-1,1.5) -- (-2.5-1,3.5);
\draw (-2-1,1.5) -- (-2-1,3.5);
\draw (-1.5-1,1.5) -- (-1.5-1,3.5);

\node at (0,-.1) { };
\node at (0,5.1) { };


\draw[] (2.5,2.5) circle [radius=4];

\draw[black,dotted,line width=3pt, line cap=round, dash pattern=on 0pt off 6\pgflinewidth] (2.5,2.5) circle [radius=4];


\begin{scope}
\pgfsetfillpattern{dots}{red}
\filldraw (0,0) rectangle (5,5);
\end{scope}

\node at (-2-1,1) {$u^{\emph{inc}}$};
\node at (8,1) {$u^{\emph{meas}}$};
\node at (2.5,2.5) {\bf $\text{supp(q)=?}$};
\node at (6,-0.5) {$\parB$};
\node at (7,6.5) {receivers};

\draw[ultra thick,->] (6+1,2) -- (7.75+1,2.75);
\draw (6.5+1,1.25) -- (6+1,3);
\draw (7+1,1.5) -- (6.5+1,3.25);
\draw (7.5+1,1.75) -- (7+1,3.5);

\draw[thick,->] (7,6.2) -- (6.1,4.6);

\end{tikzpicture}

\caption{Inverse scattering problem}\label{fig:inv_prob}
\end{subfigure}
\caption{Scattering from point scatterers represented by $q(\xb)=\sum_{i=1}^N\delta(\xb-\xb_i)q_i$. In the {\em forward scattering problem}, $q(\xb)$ is known and one seeks to compute the scattered field given the incident field at receivers on the boundary $\partial\mathcal{B}$ of a disk, see Figure \ref{fig:fwd_prob}. In the {\em inverse scattering problem}, $q(\xb)$ is unknown and we seek to determine it from measurements of the scattered field at the receivers located on $\partial\mathcal{B}$, as depicted in Figure \ref{fig:inv_prob}.}\label{fig:problems_intro}
\end{figure}
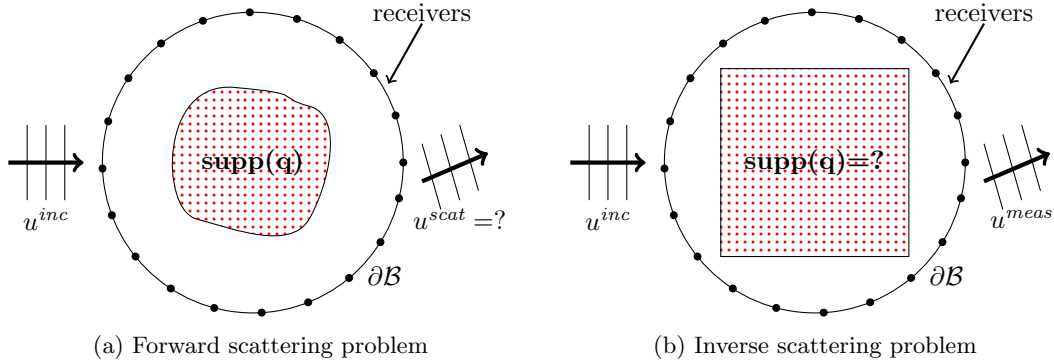

We consider the forward and inverse scattering problems in two dimensions. The problem setup is summarized in Figure \ref{fig:problems_intro}. We define $q(\xb)$ to represent a collection of point scatterers. We define the forward scattering operator $\Fb$ that maps $q$ into the scattered field by
\begin{equation}
\Fb (q;u^{\emph{inc}}) = u^{\emph{scat}},
  \label{eq:analytic_operator}
\end{equation}
where $u^{\emph{inc}}$ is the incident field and $u^{\emph{scat}}$ is the scattered field. The operator $\Fb$ in \eqref{eq:analytic_operator} is well-defined since the forward scattering problem is well-posed \cite{Colton}. 

We compute the scattered field using the Lippmann-Schwinger Equation \cite{Colton}. For a domain formed by $N$ point scatterers, the Lippman-Schwinger Equation becomes a finite dimensional problem \eqref{eq:discrete_forward}. To obtain the scattered field it amounts to solving an $N\times N$ dense linear system of equations using a preconditioned Krylov subspace iterative method such as GMRES \cite{Saad86} and combining it with the fast multipole method (FMM) \cite{Crutchfield,CHENG06}. In this case, the computational complexity necessary for the convergence of the GMRES method to a set tolerance $\epsilon_F$ is $\mathcal{O}(\kappa_F N)$ for low frequency problems and $\mathcal{O}(\kappa_F N\log(N))$ for high frequency problems, where $\kappa_F$ is the number of Krylov iterations.


In the inverse scattering problem, given observations of the scattered field we seek an approximation $\tilde{q}$ of the exact scatterers by solving
\begin{equation}
\tilde{q}=\argmin_q \frac{1}{2}\biggr\rvert\biggr\rvert\db-\Fb\left(q\right)\biggr\rvert_{\partial \mathcal{B}}\biggr\rvert\biggr\rvert^2,
\label{eq:invprob}
\end{equation}
where $\db$ is the vector with components being the measurements of the scattered field $u^{\emph{meas}}$ at receivers, $\Fb(q)\biggr\rvert_{\partial \mathcal{B}}$ is the vector with values of the forward operator evaluated at the receivers. This inverse problem is nonlinear and ill-posed. To solve it, we use a Gauss-Newton method combined with Tykhonov regularization. At each step of the Gauss-Newton iteration we solve the linear system with the Gauss Newton Hessian \eqref{eq:Hessian_prob_1}, where $\Hb=\Jb^\ast \Jb+\beta \Ib\in\mathbb{R}^{N\times N}$, $\Jb$ is the $N_dN_r\times N$ matrix representing the Fr\'{e}chet derivative of the forward operator with respect to $q$, $\Jb^\ast$ is the adjoint of $\Jb$, $N_d$ is the number of incident waves used to obtain measurements of the scattered field, $N_r$ is the number of receivers where the scattered field is measured, and $\beta \Ib$ is the regularization term. As with the forward scattering problem, constructing $\Hb$ is expensive for large $N$. So, to solve systems with $\Hb$ we use a matrix-free iterative method. In each application of the Fr\'{e}chet derivative, we need to solve $N_d$ forward scattering problems. Supposing that each of those problems takes approximately the same number of iterations $\kappa_F$ to converge, the total number of operations for the multiplication of $\Jb$ by a vector is $\mathcal{O}(\kappa_F N_d N_r N)$. The cost of inverting $\Hb$ is equal to the cost of applying $\Hb$ $\kappa_I$ times, which is proportional to $\kappa_I \kappa_F N_d N$. 


Our only concern in this paper is how to efficiently solve \eqref{eq:invprob} by proposing a novel preconditioning scheme. We do not consider issues with nonlinear solvers and more sophisticated regularization techniques. In a nutshell, our contributions in this paper are summarized as follows:
\begin{itemize}
\item For the forward problem, we apply and compare the domain decomposition based overlapping preconditioners presented in \cite{Marcus99} on the Lippmann-Schwinger Equation for point scatterers on domains up to 40 wavelengths. We propose and test the Additive Schwarz (AS), the Restricted Additive Schwarz (RAS), the Additive Harmonic Schwarz(AHS) and the Symmetric Restricted Additive Schwarz (SRAS) preconditioners for integral equations. We also propose a low rank correction to the above schemes. We term it the RC preconditioner (for the rank-corrected preconditioner).
\item For the inverse problem, we use the RC preconditioner proposed for the forward problem to construct an approximation of the Gauss-Newton Hessian, which is then used to precondition $\Hb$. We term this the HRC preconditioner. We use the preconditioners developed in the solution of a multifrequency full aperture inverse scattering problem. Our key contribution is the preconditioner for the Gauss-Newton Hessian operator. We compare this preconditioner with the low-rank preconditioners. However, in our case, the Hessian's rank can be high and our preconditioner scales better.
\end{itemize}

A list of the numerical experiments with their respective descriptions and results is given in Table \ref{table:experiments}. Experiments F.1-5 are for the forward problem and experiments I.1-5 for the inverse problem.

\begin{table} 
\center
\caption{List of numerical experiments with the respective Tables and Figures.}\label{table:experiments}
{
\begin{center}
\begin{tabular}{|c|l|c|c|}
\hline
Experiment & Description & Tables & Figures\\
\hline\hline
\multirow{2}{*}{F.1} & GMRES using DD preconditioning on          & \multirow{2}{*}{\ref{table:fwd_exp1}, \ref{table:fwd_exp1_sup}} & \multirow{2}{*}{\ref{fig:domain_example}, \ref{fig:q4_scat_field}, \ref{fig:eig_example_f1_1}, \ref{fig:eig_example_f1_2}} \\
      &  forward problem								       &									& \\\hline
\multirow{2}{*}{F.2} & Comparison of different partition geometries                          & \multirow{2}{*}{\ref{table:fwd_exp2}, \ref{table:fwd_exp2_sup}} & \multirow{2}{*}{\ref{fig:q4_scat_field}}  \\
      & on DD forward preconditioning								   &									     &   \\\hline 
\multirow{2}{*}{F.3} & Comparison of different scatterers on DD                               & \multirow{2}{*}{\ref{table:fwd_exp3}, \ref{table:fwd_exp3_sup}} & \multirow{2}{*}{\ref{fig:domain_example_16}, \ref{fig:q16_scat_field}} \\
      & forward preconditioning											&									     &   \\\hline 
\multirow{2}{*}{F.4} & Scalability of the DD preconditioning on    & \multirow{2}{*}{\ref{table:fwd_exp4}, \ref{table:fwd_exp4_sup}} & \multirow{2}{*}{X} \\
      &  forward problem											      &										&    \\\hline
F.5 & RC preconditioner on forward problem					      & \ref{table:fwd_exp5}, \ref{table:fwd_exp5_sup} & X \\\hline\hline
\multirow{2}{*}{I.1} & Influence of the overlapping on inverse                    & \multirow{2}{*}{\ref{tab:error_exp_1_2}}   & \multirow{2}{*}{\ref{fig:errorH_inv_ol}, \ref{fig:error_inv_ol_ft}} \\
     &preconditioning										      &									&	\\\hline
\multirow{2}{*}{I.2} & Influence of the number of subdomains on              & \multirow{2}{*}{\ref{tab:error_exp_1_2}}  & \multirow{2}{*}{\ref{fig:errorH_inv_sd}, \ref{fig:error_inv_sd_ft}} \\
			     & inverse preconditioning					      & 							      & \\\hline 
I.3 & Scalability of the inverse preconditioning 		           				& \ref{tab:inv_scal_it}, \ref{tab:error_tf}, \ref{tab:inv_scal_error} & X \\\hline
\multirow{2}{*}{I.4} & Comparison with low-rank inverse                                                     &\multirow{2}{*}{X}  & \multirow{2}{*}{\ref{fig:Hcomparison_I4}, \ref{fig:Heig_example_I4_k5}, \ref{fig:Heig_example_I4_k20}} \\
      & preconditioner											      &										&    \\\hline
I.5 & Nonlinear inverse problem 									& \ref{tab:inv_full} & \ref{fig:domain_example_I4}, \ref{fig:inv_full_sol}, \ref{fig:inv_full_bench}\\\hline
\end{tabular}
\end{center}
}
\end{table}

{\bf Notation:} We present the most common symbols used in this paper in Table \ref{table:symbol}. 

\begin{table}
\caption{List of main symbols used in this article.}\label{table:symbol}
{\small
\begin{center}
\begin{tabular}{ll}
\hline
Symbol & Description \\
\hline\hline
$N$                           & Number of point scatterers \\
$q$                            & Collection of $N$ point scatterers representing the medium ($q\in \mathbb{R}^N$) \\
$\Qb$		         & $N\times N$ diagonal matrix with diagonal $q$ \\
$N_r$                        & Number of receivers \\
$N_\lambda$            & Number of singular values of the forward operator used for the RC preconditioner\\
$N_s$                       & Number of subdomains used in the partition of the domain\\
$\theta$                     & Incident direction of plane wave $u^{\emph{inc}}$ \\
$k$                           & Wavenumber (or frequency) of the incident plane wave \\
$N_d$                       & Number of incident waves \\
$\Fb$                        & Forward scattering operator mapping $q$ to the scattered field (for given illumination) \\
$\Gb$                        & $N\times N$ complex matrix with entries defined in Equation \eqref{eq:green_source_source}\\
$\Gb_r$                     & $N_r\times N$ complex matrix with entries defined in Equation \eqref{eq:green_rec_source}\\
$\Ib$                          & Identity matrix conforming in size to the other matrices in the equation\\
$\Ab$                        & $N\times N$ complex matrix $\Ib+k^2\Gb\Qb$ of the forward scattering problem \\
$\Jb$                        & $N_dN_r\times N$ matrix of the Fr\'{e}chet derivative of the forward operator $\Fb$\\
$\Hb$                        & $N\times N$ Gauss-Newton Hessian such that $\Hb=\Jb^\ast \Jb+\beta \Ib$ \\
$\beta$                     & Regularization parameter for the Gauss-Newton Hessian\\
$\parB$                     & Circle enclosing support of the scatterer where $N_r$ receivers are located \\
$u^{\emph{inc}}$      & Incident plane wave \\
$u^{\emph{scat}}$    & Scattered field off of $q$ \\
$u^{\emph{meas}}$  & Scattered field off of $q$ measured at the receivers\\
$\db$                        & Data measurements of the scattered field on the receivers ($\db\in \mathbb{R}^{N_dN_r}$)\\
$\Omega$                & Support of $q$ \\
$\Omega_i^\delta$   & Subdomain with overlap $\delta$ \\
$\Rb_i^\delta$          & Restriction operator for subdomain $i$ with overlap $\delta$ \\
$\tilde{\Ab}^{-1}$      & Preconditioner for the forward problem \\
$\tilde{\Hb}^{-1}$      & Preconditioner for the inverse problem \\
$\kappa_F$              & Number of iterations for unpreconditioned GMRES for the the forward problem \\
$\kappa_I$               & Number of iterations for unpreconditioned GMRES for the the inverse problem \\
$e_{rel}$                  & Relative error of the iterative solution with respect to the direct method solution \\
\hline 
\end{tabular}
\end{center}
}
\end{table}


{\bf Limitations:} First, our method is only applicable in low frequency problems. Second, here we are considering model problems in two dimensions. Third, this is a mostly experimental work. Our method has several parameters( e.g., overlap size, the type of AS, the rank in the RC preconditioner), which are currently chosen on an empirical way.


{\bf Related work:} 
The topic of domain decomposition for the solution of the partial differential equation formulation of the forward problem \eqref{eq:analytic_operator}, has been extensively studied, see \cite{BENAMOU199768,Marcus99,doi:10.1137/130917144,ey11,ERLANGGA2004409,Erlangga2008, Farhat2000, Nataf, Graham2017, 10.1007/978-3-642-35275-1_16, Mcinnes98additiveschwarz, CSTOLK2013240, ZEPEDANUNEZ2016347, doi:10.1137/15M104582X}. There is less work on the study of preconditioners in the solution of the integral equation form of the forward problem \cite{doi:10.1080/00207160802033335,BRUNO2004670,Moskow,doi:10.1137/140985147,doi:10.1137/16M1064660}. Notice, that the PDE formulation leads to sparse matrices, whereas, an integral equation formulation results in dense matrices. For this reason, it is not clear that an overlapping DD preconditioner will be effective for integral equations. For the inverse scattering problem, a few groups have used spectral strategies to speed-up the solution of the Gauss-Newton iteration, see \cite{Borges2017,CHAILLAT20124403,Hohage2001,doi:10.1137/16M1070475,Biros}. We are not aware of any work that uses traditional domain decomposition methods to  precondition the Hessian. To make the simulations simpler, we make some assumptions to simplify the problem, such as considering the problem in two dimensions and that our domain of integration is composed of point scatterers that are uniformly distributed in the unit square.

{\bf Article Outline:} In Section \ref{s:problem}, we describe the forward and inverse scattering problem formulations and how to obtain their numerical solutions. In Section \ref{s:method_forward}, we give a brief introduction about domain decomposition techniques used in this article,  we present and compare the domain decomposition preconditioners used for the forward problem, and we present the low rank correction procedure to obtain the RC preconditioner. In Section \ref{s:method_inverse}, we present our preconditioning strategy for the inverse scattering problem. Concluding remarks are made in Section \ref{s:conclusions}.

%% file: problem.tex
\section{The forward and inverse scattering problem}\label{s:problem}
In this section, we present the formulation of both the forward and inverse scattering problems with the assumptions used to simplify those problems. We then present the discrete system that we want to solve for both problems.

\subsection{Forward scattering problem}
The operator $\Fb$ in \eqref{eq:analytic_operator} is well-posed since the forward scattering problem has a unique and stable solution \cite{Colton}. To find the value of the scattered field we solve the equation
\begin{equation}
\Delta u^{scat}+k^2(1+q) u^{scat} = -k^2 q u^{inc}
\label{eq:Helmholtz}
\end{equation}
where $u^{scat}$ satisfies the Sommerfeld radiation condition, $u^{inc}$ is an incoming incident plane wave, and $k$ is the wavenumber.

Using Green's Formula for the Helmholtz equation (Chapter 2 of \cite{Colton}) and the Sommerfeld radiation condition, we obtain the integral form of Equation \eqref{eq:Helmholtz} that reads
\begin{equation}
u^{scat}(\xb)+k^2\int_{\mathbb{R}^2} G(k\|\xb-\yb\|)q(\yb)u^{scat}(\yb) d\yb= -k^2 \int_{\mathbb{R}^2} G(k\|\xb-\yb\|)q(\yb) u^{inc}(\yb) d\yb,
\label{eq:LipSch}
\end{equation}
the Lippman-Schwinger Equation, where $G(k\|\xb-\yb\|)=\frac{i}{4}H^{(1)}_0(k\|\xb-\yb\|)$ is the free-space Green's function for the two-dimensional Helmholtz equation and $H^{(1)}_0$ is the Hankel function of first kind of order zero.

We assume that $q(\xb)$ is a set of point scatterers distributed in a regular $\sqrt{N}\times \sqrt{N}$ grid of $\xb_i$ points in the square $\left[-0.5,0.5\right]^2$, so that 
\begin{equation}
q(\xb)=\sum_{i=1}^{N} q_i \delta(\xb-\xb_i),
\label{eq:domain_point}
\end{equation}
where $q_i$ is the charge of the point scatterer located at $\xb_i$. Then equation \eqref{eq:LipSch} becomes
\begin{equation}
u(\xb) + k^2 \sum_{i=1}^{N} q_i G(k\|\xb-\xb_i\|) = -k^2 \sum_{i=1}^{N} q_i G(k\|\xb-\xb_i\|).
\label{eq:domain_point_1}
\end{equation}
Taking $\xb=\xb_j$ and ignoring the self iterations, we obtain
\begin{equation}
u(\xb_j) + k^2 \sum_{\substack{i=1 \\ i\neq j}}^N q_i G(k\|\xb_j-\xb_i\|) = -k^2 \sum_{\substack{i=1 \\ i\neq j}}^N q_i G(k\|\xb_j-\xb_i\|).
\label{eq:domain_point_2}
\end{equation}

Using equation \eqref{eq:domain_point_2} on all the scatterer points in the domain, we obtain the system of equations
\begin{equation}
(\Ib+k^2\Gb\Qb)\ub^{scat} = -k^2 \Gb \Qb \ub^{inc},
\label{eq:discrete_forward}
\end{equation}
where $\Gb$ is the $N\times N$ matrix with elements
\begin{equation}\label{eq:green_source_source}
(\Gb)_{ij}=\begin{cases}
G(k\|\xb_i-\xb_j\|) &\mbox{if }  i\neq j \\
0 &\mbox{if } i=j, 
\end{cases}
\end{equation} 
 for $i,j=1,\ldots,N$, $\Ib$ is the $N\times N$ identity matrix, $\Qb$ is the $N\times N$ diagonal matrix with diagonal elements $(\Qb)_{ii}=q_i$, $\ub^{inc}$ is a vector with coordinates $(\ub^{inc})_i=u^{inc}(\xb_i)$ and the solution vector $\ub^{scat}$ is such that for each coordinate we have $(\ub^{scat})_i=u^{scat}(\xb_i)$. We define $\Ab=(\Ib+k^2\Gb\Qb)$. Our first goal is to define preconditioners for $\Ab$.

After solving the system \eqref{eq:discrete_forward} and finding $u^{scat}(\xb_i)$ at the scatterer positions, we can compute $u^{scat}(\xb)$ at any point using equation \eqref{eq:domain_point_1}. In particular, at the receivers we have
\begin{equation}
\Fb\left(q\right)\biggr\rvert_{\partial \mathcal{B}}=\ub^{\emph{meas}} =  -k^2 \Gb_r\Qb\ub^{\emph{inc}} - k^2\Gb_r\Qb\ub^{\emph{scat}},
\label{eq:Fvalue}
\end{equation}
where $\Gb_r$ is the $N_r\times N$ matrix with elements 
\begin{equation}
(\Gb_r)_{ij}=G(k\|\yb_i-\xb_j\|), 
\label{eq:green_rec_source}
\end{equation}
where $\yb_i$ are the coordinates of the $N_r$ receivers, $\ub^{\emph{meas}}$ is a vector with coordinates $(\ub^{\emph{meas}})_i=u^{meas}(\xb_i)$, and $u^{meas}(\xb_i)$ is the measured scattered field at the receivers.

\subsection{Inverse scattering problem}
Suppose we have an incoming plane wave $u^{inc} = \exp(ik\xb\cdot\theta)$ and we are given the vector $\db$ of data measurements of the scattered field at the receivers, we seek to find $\tilde{q}$ such that
\begin{equation}
\tilde{q}=\argmin_q \frac{1}{2}\biggr\rvert\biggr\rvert\db-\Fb\left(q\right)\biggr\rvert_{\partial \mathcal{B}}\biggr\rvert\biggr\rvert^2,
\label{eq:invprob_1}
\end{equation}
where the coordinates of $\db$ are $\db_j=u^{meas}(\xb_{r_j})$ and $\xb_{r_j}$ are the positions of the receivers, for $j=1,\ldots,N_r$.

We use the Gauss-Newton method to solve equation \eqref{eq:invprob_1} for $\tilde{q}$ \cite{Colton,Borges2017}. The Gauss-Newton step is given by
\begin{equation}
\Jb \delta q = \left( \db - \Fb\left(q^{(i)}\right)\biggr\rvert_{\partial \mathcal{B}}\right),
\label{eq:GNequation_0}
\end{equation}
where the matrix $\Jb$ is the discrete version of the Fr\'{e}chet derivative of $\Fb$ with respect to $q$. In our experiments, we consider scenarios in which the number of measurements is higher than the number of scatterer points, the system \eqref{eq:GNequation_0} is overdetermined. The expression for $\Jb$ given by
\begin{equation}
\Jb = -k^2 \Gb_r \Ub^{tot} +k^4 \Gb_r\Qb \Ab^{-1}\Gb\Ub^{tot},
\label{eq:Joperator}
\end{equation}
where $\Ub^{tot}$ is the $N\times N$ diagonal matrix with diagonal elements $(\Ub)^{tot}_{i,i}=u^{inc}(\xb_i)+u^{scat}(\xb_i)$.

In each iteration, the Tikhonov regularized Gauss-Newton step becomes
\begin{equation}
\Hb \delta q = \Jb^\ast \left( \db - \Fb\left(q^{(i)}\right)\biggr\rvert_{\partial \mathcal{B}} \right),
\label{eq:Hessian_prob_1}
\end{equation}
where $\Hb=(\Jb^\ast \Jb +\beta \Ib)$ is the Hessian of the problem, $\beta$ is the regularization parameter and $\Jb^\ast$ is the adjoint of $\Jb$. The Gauss-Newton method (without line search) is summarized in Algorithm \ref{alg:gnmethod}.

\begin{algorithm}
\caption{Gauss-Newton method for the inverse scattering problem.}
\label{alg:gnmethod}
\begin{algorithmic}[1]
\STATE{{\bf Input:} data $\db$, initial guess $q_0$, tolerances $\epsilon_1$,$\epsilon_2$, and maximum number of iterations $N_{it}$.}
\STATE{Set $q\coloneqq q_0$, $\delta q\coloneqq 0$ and $it\coloneqq 0$.}
\WHILE{$\|\db-\Fb(q)\|\geq\epsilon_1$ or $it<N_{it}$ or $\delta q\geq\epsilon_2$} 
\STATE{Solve $\Ab\ub^{\emph{scat}}=-k^2\Gb\Qb\ub^{\emph{inc}}$ using GMRES}
\STATE{Calculate $\Fb(q)\vert_{\partial \mathcal{B}}=-k^2 \Gb_r\Qb\ub^{\emph{inc}} - k^2 \Gb_r\Qb\ub^{\emph{scat}}$}
\STATE{Solve $\Hb\, \delta q=\Jb^\ast \left( \db-\Fb\left(q\right)\biggr\rvert_{\partial \mathcal{B}}\right)$ using GMRES}
\STATE{Update $q\leftarrow q+\delta q$}
\STATE{Update $it\leftarrow it+1$}
\ENDWHILE
\STATE{The approximate solution is $\tilde{q}\coloneqq q$.}
\end{algorithmic}
\end{algorithm}

The extension of the Gauss-Newton method using data from $N_d$ impinging incident plane waves is straightforward. Consider that $\db_j$ and $\Fb_j$ are the data measurements and forward operator referent to the plane wave with incident direction $\theta_j$, $j=1,\ldots,N_d$. We seek to find $\tilde{q}$ such that
\begin{equation}
\tilde{q}=\argmin_q \sum_{\substack{j=1}}^{N_d} \frac{1}{2}\biggr\rvert\biggr\rvert\db_j-\Fb_j\left(q\right)\biggr\rvert_{\partial \mathcal{B}}\biggr\rvert\biggr\rvert^2.
\label{eq:invprob_md}
\end{equation}
In this case, we have
\begin{equation}
\begin{bmatrix}
\Jb_1 \\
\vdots \\
\Jb_{N_d}
\end{bmatrix} \delta q = \begin{bmatrix} \db_1 - \Fb_1\left(q^{(i)}\right)  \\ \vdots \\ \db_{N_d} - \Fb_{N_d}\left(q^{(i)}\right) \end{bmatrix}
\end{equation}
where $\Jb_j$ is given by the formula \eqref{eq:Joperator} using $\Ub^{tot}_j$ instead of $\Ub^{tot}$, where $\Ub^{tot}_j$ is the $N\times N$ diagonal matrix with diagonal elements $u_j^{inc}(\xb_i)+u_j^{scat}(\xb_i)$, $u_j^{inc}$ is the incoming field with direction of propagation $\theta_j$ and $u_j^{scat}$ is its respective scattered field. If we set
\begin{equation}
\Jb=\begin{bmatrix}
\Jb_1 \\
\vdots \\
\Jb_{N_d}
\end{bmatrix}, \quad 
\db =\begin{bmatrix}
\db_1 \\
\vdots \\
\db_{N_d}
\end{bmatrix}, \quad \text{and} \quad
\Fb=\begin{bmatrix}
\Fb_1\left(q^{(i)}\right) \\
\vdots \\
\Fb_{N_d}\left(q^{(i)}\right)
\end{bmatrix},
\end{equation}
we obtain the system \eqref{eq:GNequation_0}, and consequently all our discussion on the Gauss-Newton method directly applies. From this point forward, we will suppress the notation indicating the direction of the incident plane wave when referring to the measurements and the forward operator. We consider that the operators and measurements are being used for multiple incoming waves.

%% file: fwd_prec.tex
\section{Overlapping domain decomposition preconditioning of the forward problem}\label{s:method_forward}

\begin{figure}
\centering
\begin{subfigure}{.4\linewidth}
\begin{tikzpicture}[scale=.7]
\draw[thick] (-2,-2) rectangle (2,2);
\draw[thick] (2,2) rectangle (6,6);
\draw[thick] (2,-2) rectangle (6,2);
\filldraw[fill=blue!40!white, draw=black] (-2,2) rectangle (2,6);
\node at (0,4) {\Large$\Omega_1^0$};
\node at (0,0) {\Large$\Omega_2^0$};
\node at (4,4) {\Large$\Omega_3^0$};
\node at (4,0) {\Large$\Omega_4^0$};
\end{tikzpicture}
\end{subfigure}
\begin{subfigure}{.4\linewidth}
\begin{tikzpicture}[scale=.7]
\filldraw[fill=blue!40!white] (-2,1.5) rectangle (2.5,6);
\draw[thick] (-2,-2) rectangle (2,2);
\draw[thick] (2,2) rectangle (6,6);
\draw[thick] (2,-2) rectangle (6,2);
\draw[thick] (-2,2) -- (-2,6);
\draw[thick] (-2,6) -- (2,6);
\draw[thick,dashed] (-2,1.5) rectangle (2.5,6);
\node at (0,4) {\Large$\Omega_i^\delta$};
\end{tikzpicture}
\end{subfigure}

\begin{subfigure}{.4\linewidth}
\begin{tikzpicture}[scale=.7]
\draw[thick] (-2,0) -- (6,0);
\draw[thick] (-2,2) -- (2,2);
\draw[dashed] (2,2) -- (2,0);
\node at (2,3) {\Large$R_i^0$};
\draw[thick,<->] (-2,-1) -- (2,-1);
\node at (0,-1.5) {\Large$\Omega_i^0$};
\end{tikzpicture}
\end{subfigure}
\begin{subfigure}{.4\linewidth}
\begin{tikzpicture}[scale=.7]
\draw[thick] (-2,0) -- (6,0);
\draw[thick] (-2,2) -- (1.5,2);
\draw[dashed] (1.5,2) -- (1.5,0);
\draw[dashed] (2.5,2) -- (2.5,0);
\draw[thick] (1.5,2) -- (2.5,2);
\node at (2,3) {\Large$R_i^\delta$};
\draw[thick,<->] (-2,-1) -- (2.5,-1);
\draw[thick,<->] (1.5,-.2) -- (2.5,-.2);
\node at (1.9,-0.5) {\large$\delta$};
\node at (0.5,-1.5) {\Large$\Omega_i^\delta$};
\end{tikzpicture}
\end{subfigure}

\begin{subfigure}{.4\linewidth}
\begin{tikzpicture}[scale=.7]
    \begin{axis}[grid=major,view={45}{45},
    xmin=-.5, xmax=.5,
    ymin=-.5, ymax=.5,
    zmin=  -.2, zmax=1.2]

\draw [color=gray, fill=blue!15] (-.5,0,1) -- (-.5,.5,1) -- (0,0.5,1) -- (0,0,1) -- (-.5,0,1);
\draw [color=gray, fill=blue!15] (-.5,-.5,0) -- (.5,-.5,0) -- (.5,.5,0) -- (0,.5,0) -- (0,0,0) -- (-.5,0,0) --  (-.5,-.5,0);
      \addplot3 [data cs=cart,surf,domain=-.5:.5, opacity=0.0]
      {0};      
    \end{axis}
  \end{tikzpicture}
\end{subfigure}
\begin{subfigure}{.4\linewidth}
\begin{tikzpicture}[scale=.7]
    \begin{axis}[grid=major,view={45}{45},
    xmin=-.5, xmax=.5,
    ymin=-.5, ymax=.5,
    zmin=  -.2, zmax=1.2]

\draw [color=gray, fill=blue!15] (-.5,-.5,0) -- (.5,-.5,0) -- (.5,.5,0) -- (.1,.5,0) -- (.1,-.1,0) -- (-.5,-.1,0) --  (-.5,-.5,0);
\draw [color=gray, fill=blue!15] (-.5,-.1,1) -- (-.5,.5,1) -- (.1,0.5,1) -- (.1,-.1,1) -- (-.5,-.1,1); 

      \addplot3 [data cs=cart,surf,domain=-.5:.5, opacity=0.0]
      {0};      
    \end{axis}
  \end{tikzpicture}
\end{subfigure}
\caption{In each row, from top to bottom: partition of the domain $\Omega$, cross-section of the restriction operator for $\Omega_0$ and the restriction operator for $\Omega_0$. The left column has a non overlapping partition and the right column an overlapping partition.}
\label{fig:DD_partition}
\end{figure}
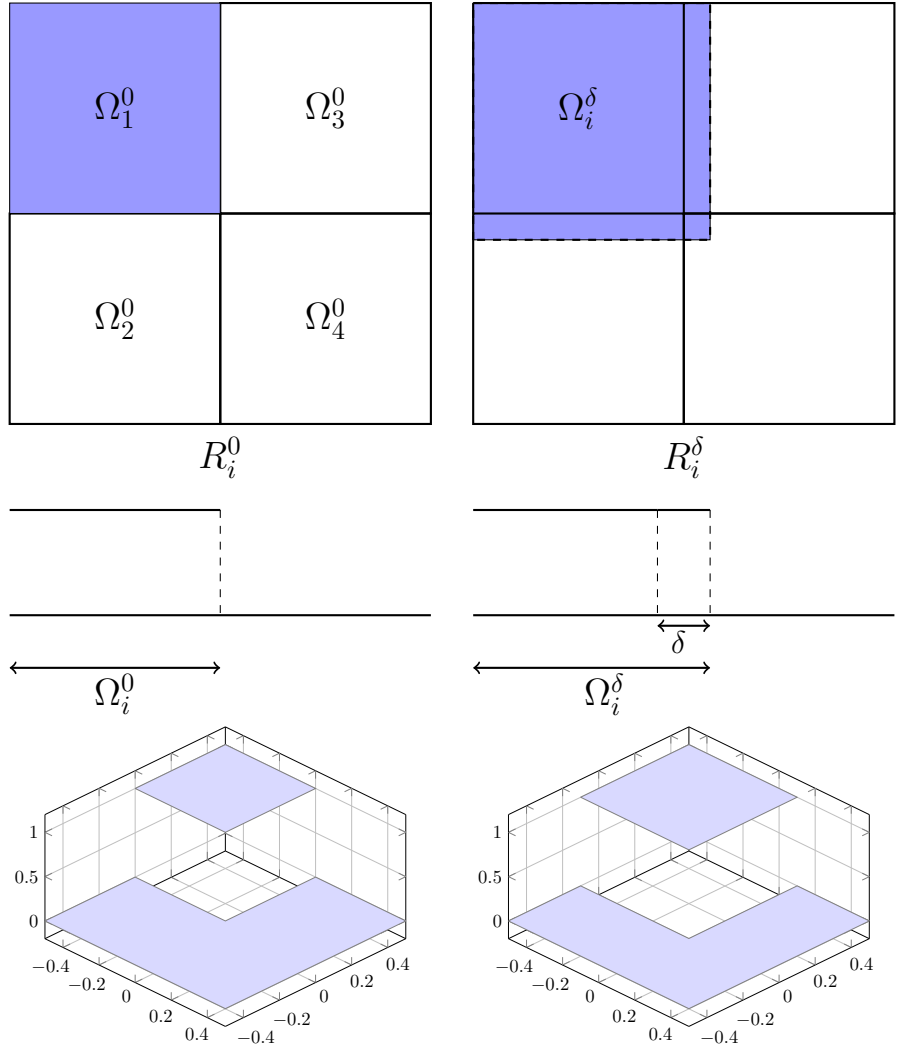

As we discussed, $q$ is composed of $N$ point scatterers distributed in a $\sqrt{N}\times \sqrt{N}$ regular grid in the domain $\Omega=\left[-.5,.5\right]^2$. We partition $\Omega$ into $N_{s}$ nonoverlapping subdomains. Without loss of generality we consider that the domain is partitioned in a perfect square number of same size squares as in Figure \ref{fig:DD_partition}. We have
\begin{equation}
\Omega=\bigcup_{i=1}^{N_{s}} \Omega_i^0.
\label{eq:partition_nonover}
\end{equation}
We define the overlapping partition of $\Omega$, as follows. Let $\Omega_i^\delta$ be the overlapping partition of $\Omega$, where $\Omega_i^\delta \supset \Omega_i^0$ is obtained by increasing the size of $\Omega_i^0$ by $\delta$, where $\delta$ is a measure of size of the overlap of the domains (it can be a percentage of the size of $\Omega_i^0$ or it can be the number of points in $\Omega_i^\delta$ not in $\Omega_i^0$). With this definition, we have:
\begin{equation}
\Omega=\bigcup_{i=1}^{N_{s}} \Omega_i^\delta.
\label{eq:partition_over}
\end{equation}

Associated with each $\Omega_i^0$, we define a restriction operator $\Rb_i^0$. In matrix terms, $\Rb_i^0$ is an $N\times N$ diagonal matrix whose diagonal elements are set to one if the corresponding point belongs to $\Omega_i^0$ and to zero otherwise. We define $\Rb_i^\delta$ in an analogous way, such that its diagonal elements are set to one if the corresponding point belongs to $\Omega_i^\delta$ and to zero otherwise  From the definition, we have
\begin{equation}
\sum_{i=1}^{N_{s}} \Rb_i^0(j,j) =1
\label{eq:restriction_op_over}
\end{equation}
for $j=1,\cdots,N$, where $\Rb_i^0(j,j)$ is the $j_\mathrm{th}$ diagonal element of $\Rb_i^0$, and
\begin{equation}
\sum_{i=1}^{N_{s}} \Rb_i^\delta(j,j) \geq 1
\label{eq:restriction_op_nonover}
\end{equation}
for $j=1,\cdots,N$, where $\Rb_i^\delta(j,j)$ is the $j_\mathrm{th}$ diagonal element of $\Rb_i^\delta$.

\subsection{Domain decomposition preconditioning of the forward problem} To simplify notation, we denote $\ub=\ub^{\emph{scat}}$ and $\bb=-k^2\Gb\Qb\ub^{\emph{inc}}$. We want to speed-up the convergence of GMRES for solving
\begin{equation}
\Ab\ub=\bb.
\label{eq:forward_lin_system}
\end{equation}
Using a domain decomposition based preconditioner, the system to be solved becomes
\begin{equation*}
 \tilde{\Ab}^{-1}\Ab \ub= \tilde{\Ab}^{-1} \bb,
\end{equation*}
where $\tilde{\Ab}^{-1}$ is the proposed preconditioner.

We start by defining the $N\times N$ matrices
\begin{equation}
\Ab_{ii} =\Rb_i^\delta \Ab \Rb_i^\delta,
\label{eq:Asubmatrix}
\end{equation}
for each subdomain $\Omega_i^\delta$, $i=1,\ldots,N_s$. Note that the $(j,m)-$entry on matrix $\Ab_{ii}$ is given by
\begin{equation}\label{eq:Aii_entry}
\Ab_{ii}(j,m)=
\begin{cases}
\Ab(j,m), & \mbox{if } \xb_j,\xb_m\in \Omega_i^\delta,  \\
0,              & \mbox{otherwise}.
\end{cases}
\end{equation}
Suppose $\Omega_i^\delta$ is composed by $\hat{n}$ points $\xb_{i_1},\ldots,\xb_{i_{\hat{n}}}$. The restriction of $\Ab_{ii}$ to $\Omega_i^{\delta}$ is defined as the $\hat{n}\times\hat{n}$ matrix $\hat{\Ab}_{ii}=\Ab_{ii}($$i_1$:$i_{\hat{n}}$,$i_1$:$i_{\hat{n}})$. Although the matrix $\Ab_{ii}$ is not invertible, $\hat{\Ab}_{ii}$ can be inverted. We define the inverse of the matrix $\Ab_{ii}$ as the matrix $\Ab_{ii}^{-1}$ with $(i_j,i_m)-$entry given by
\begin{equation}
\Ab_{ii}^{-1}(i_j,i_m) = 
\begin{cases}
\hat{\Ab}_{ii}^{-1}(j,m), & \mbox{if } \xb_{i_j},\xb_{i_m}\in \Omega_i^\delta,  \\
0,              & \mbox{otherwise}.
\end{cases}
\label{eq:rest_inv}
\end{equation}

In \cite{Marcus99}, the authors proposed several variants of the Additive Schwarz preconditioners and compared these variants for systems of equations obtained from using the finite element method to solve the Helmholtz Equation. In this subsection, we apply all the non weighted methods listed in \cite{Marcus99} to a system of equations obtained from the Lippman-Schwinger equation applied in point scatterers. All the methods proposed in \cite{Marcus99} are listed in Table \eqref{table:precs}. In the same paper, the authors commented that it is possible to use weighted versions of those methods. In the weighted versions, instead of having the relation \eqref{eq:restriction_op_nonover} for the Restriction operator, we have a linear interpolation and have $\sum_{i=1}^{N_{s}} \Rb_i^\delta(j,j) = 1$. We made some tests with the weighted versions and we did not observe any considerable improvement in the performance. For this reason, we only report results for the non-weighted version.

\begin{table}
\caption{List of preconditioners. The matrix $A_{ii}^{-1}$ is given by Equation \eqref{eq:rest_inv}.}\label{table:precs}
\begin{center}
\begin{tabular}{l|c|c}
\hline
Method name & Abbreviation & Preconditioner \\
\hline\hline
                  Additive Schwarz & AS & $\tilde{\Ab}^{-1}_{AS} = \sum_{i=1}^{N_{s}}\Rb_i^\delta \Ab_{ii}^{-1} \Rb_i^\delta$ \\
Restricted Additive Schwarz & RAS & $\tilde{\Ab}^{-1}_{RAS} = \sum_{i=1}^{N_{s}}\Rb_i^0 \Ab_{ii}^{-1} \Rb_i^\delta$ \\
Additive Harmonic Schwarz & AHS & $\tilde{\Ab}^{-1}_{AHS} = \sum_{i=1}^{N_{s}}\Rb_i^\delta \Ab_{ii}^{-1} \Rb_i^0$ \\
Symmetrized Restricted Additive Schwarz  & SRAS & $\tilde{\Ab}^{-1}_{SRAS} = \sum_{i=1}^{N_{s}}\Rb_i^0 \Ab_{ii}^{-1}\Rb_i^0$ \\
\hline 
\end{tabular}
\end{center}
\end{table}

Assuming that $\delta\ll N$, the construction of each matrix $\Ab_{ii}^{-1}$ requires $\mathcal{O}(N^3/N_s^3)$ operations. Since this step is parallelizable, all the subdomain matrices can be constructed simultaneously and consequently the construction of the preconditioner requires $\mathcal{O}(N^3/N_s^3)$ operations.

\subsection{Improving performance by adding a low rank correctioon} 
It is well-known that domain decomposition preconditioners perform better when combined with a coarse grid solve. This inspired us to consider a low-rank correction to the overlapping schemes. We use $\tilde{\Ab}^{-1}$ for the domain decomposition preconditioner and $\tilde{\Ab}_{RC}^{-1}$ for the RC preconditioner. 

To obtain the RC preconditioner, first, we approximately compute the singular value decomposition $\Ub\Sb\Vb^\ast=\left(\Ab-\tilde{\Ab}\right)$.This decomposition can be obtained using randomized SVD methods and requires only matrix-vector multiplications \cite{GunnarPNAS}. We set $\Ub_{N_\lambda}=\Ub(:,$$1$:$N_\lambda)$, $\Sb_{N_\lambda}=\Sb($$1$:$N_\lambda,$$1$:$N_\lambda)$ and $\Vb_{N_\lambda}=\Vb(:,$$1$:$N_\lambda)$ and construct
\begin{equation}
\tilde{\Ab}_{RC}^{-1}=(\tilde{\Ab}+\Ub_{N_\lambda} \Sb_{N_\lambda} \Vb_{N_\lambda}^\ast)^{-1}.
\label{eq:rc_nonwoodbury}
\end{equation}
We can use the Woodbury formula to express $\tilde{\Ab}_{RC}^{-1}$  as a function of $\tilde{\Ab}^{-1}$
\begin{equation}
\tilde{\Ab}_{RC}^{-1} = \tilde{\Ab}^{-1}+\tilde{\Ab}^{-1}\Ub_{N_\lambda}(\Sb_{N_\lambda}^{-1}-\Vb_{N_\lambda} \tilde{\Ab}^{-1}\Ub_{N_\lambda})^{-1}\Vb_{N_\lambda} \tilde{\Ab}^{-1}.
\label{eq:rc_woodbury}
\end{equation}

The matrix $\tilde{\Ab}_{RC}^{-1}$ can be used as a preconditioner for the solution of Equation \eqref{eq:forward_lin_system}. Typically  $\tilde{\Ab}_{RC}^{-1}$ is never computed, but we can rapidly apply it using \eqref{eq:rc_woodbury}.
%

{\bf Complexity analysis:} The complexity of this scheme is the same as the complexity of building the domain decomposition based preconditioner plus the complexity of obtaining the low rank correction. This complexity is proportional to the complexity of building the inverse for each of the $N_s$ subdomains, with a total of $\mathcal{O}(N^3/N_s^2)$. One of the complications to obtain the rank correction is that $\tilde{\Ab}$ is not readily available. This can be mitigated by the application of randomized SVD. Also, in the calculation of \eqref{eq:rc_nonwoodbury}, we also do not have $\tilde{\Ab}$, however, with the use of the Woodbury formula, we obtain equation \eqref{eq:rc_woodbury}, which requires only the already available DD preconditioner. The total complexity of obtaining the low-rank correction becomes
\begin{itemize}
\item Application of randomized SVD to obtain $N_\lambda$ largest singular values and associated singular vectors: $\mathcal{O}(N \log(N_\lambda))$ using randomized algorithms, or $\mathcal{O}(NN_\lambda)$ using classical algorithms; and 
\item Application of the  Woodbury formula: $\mathcal{O}(N_\lambda^3+NN_\lambda+N)$.
\end{itemize}
The total cost of the construction of $\Ab^{-1}_{RC}$ is $\mathcal{O}(N\log(N_\lambda)+N_\lambda^3+NN_\lambda+N+N^3/N_s^2)$, where the first term is for the randomized SVD, the next three terms come from the Woodbury formula and the last term comes from the pre-calculation of the domain decomposition preconditioner.

\subsection{Numerical Experiments for Forward Preconditioning} 
We present experiments that verify the effectiveness of the domain decomposition preconditioners in different scenarios for the forward scattering problem. We compare with GMRES without preconditioner and study the effect of the number of subdomains and overlapping on the preconditioners as well as the shape of the subdomains $\Omega_i^\delta$. In the fourth experiment, we study the scalability of the domain decomposition preconditioners as we increase the number of scatterers. Those experiments show that the RAS preconditioner performs best. We combine the RAS preconditioner with RC and study its performance for different values of $N_\lambda$. A list of the experiments for the forward problem with their descriptions and results is provided in Table \ref{table:fp_experiments}.

We define the relative error of the iterative solution with respect to the direct method solution $e_{rel}\coloneqq\|\ub_{\mathrm{GMRES}}-\ub_{\mathrm{LU}}\|/\|\ub_{\mathrm{LU}}\|$, where $\ub_{\mathrm{GMRES}}$ is the solution obtained by GMRES and $\ub_{\mathrm{LU}}$ is the solution obtained by the LU direct solver. The LU solution is obtained by solving \eqref{eq:forward_lin_system} using the backslash operator in MATLAB.

We report the number of iterations for different variants of the preconditioner. In Appendix \ref{s:appendB}, we report the residual and the solution error for the experiments.

We summarize our choice of parameters for GMRES bellow
\begin{itemize}
\item GMRES configuration with preconditioners: tolerance of $10^{-13}$, no restart, and maximum number of iterations $N-1$.
\item GMRES configuration with no preconditioners: residual tolerance of $10^{-11}$, no restart, and maximum number of iterations $N-1$.
\end{itemize}

The different choice for tolerance is to account for the fact that GMRES from MATLAB uses the preconditioned residual for termination. With this different choice we aim to guarantee that the actual residuals $\bb - \Ab u$ for both unpreconditioned and preconditioned GMRES have the same order magnitude.

\begin{table}
\caption{List of forward problem  experiments.}\label{table:fp_experiments}
{\small
\begin{center}
\begin{tabular}{|c|l|c|c|}
\hline
Experiment & Description & Tables & Figures\\
\hline\hline
\multirow{2}{*}{F.1} & GMRES using DD preconditioning on forward         & \multirow{2}{*}{\ref{table:fwd_exp1}, \ref{table:fwd_exp1_sup}} & \multirow{2}{*}{\ref{fig:domain_example}, \ref{fig:q4_scat_field},\ref{fig:eig_example_f1_1},\ref{fig:eig_example_f1_2}} \\
      &  problem								       &									& \\\hline
\multirow{2}{*}{F.2} & Comparison of different partition geometries on                         & \multirow{2}{*}{\ref{table:fwd_exp2}, \ref{table:fwd_exp2_sup}} & \multirow{2}{*}{\ref{fig:q4_scat_field}}  \\
      & DD forward preconditioning								   &									     &   \\\hline 
\multirow{2}{*}{F.3} & Comparison of different scatterers on DD forward                              & \multirow{2}{*}{\ref{table:fwd_exp3}, \ref{table:fwd_exp3_sup}} & \multirow{2}{*}{\ref{fig:domain_example_16}, \ref{fig:q16_scat_field}} \\
      & preconditioning											&									     &   \\\hline 
\multirow{2}{*}{F.4} & Scalability of the DD preconditioning on forward     & \multirow{2}{*}{\ref{table:fwd_exp4}, \ref{table:fwd_exp4_sup}} & \multirow{2}{*}{X} \\
      & problem											      &										&    \\\hline
F.5 & RC preconditioner on forward problem					      & \ref{table:fwd_exp5}, \ref{table:fwd_exp5_sup} & X \\\hline
\end{tabular}
\end{center}
}
\end{table}

\begin{figure}[h!]
  \centering
  \begin{subfigure}{.4\linewidth}
\includegraphics[width=1\textwidth]{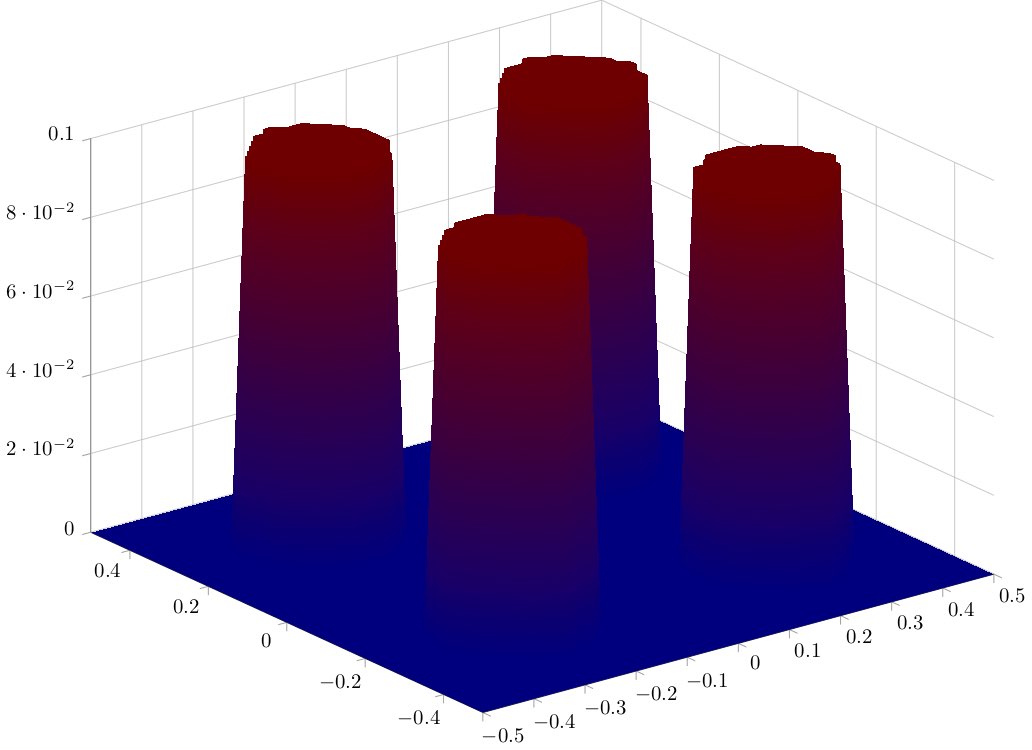}
\caption{isometric view of $q_4$}
\end{subfigure}
\begin{subfigure}{.4\linewidth}
\includegraphics[width=.8\textwidth]{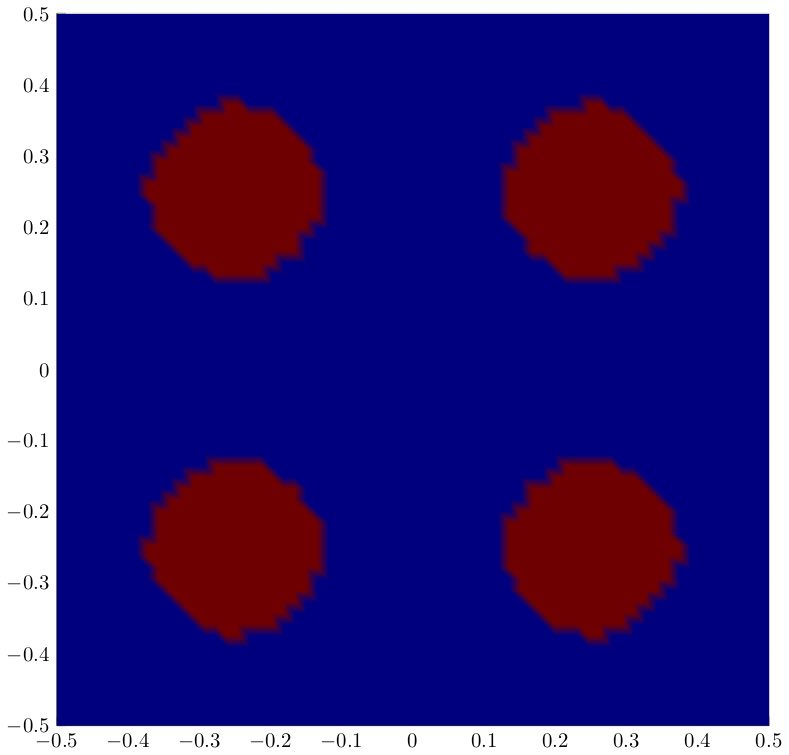}
\caption{top view of domain of $q_4$}
\end{subfigure}
\caption{The (a) isometric view and (b) top view of the domain $q_4$ used in the Experiments F.1-5, and Experiments I.1-4.}
\label{fig:domain_example}
\end{figure}

{\bf Experiment F.1 -- GMRES using DD preconditioning on forward problem:} This example aims to compare the performance of the preconditioners presented in Table \ref{table:precs} for solving Equation \eqref{eq:forward_lin_system}. We analyze the effects of the number of subdomains and size of the overlap at different wavenumbers. The following parameters are used: 
\begin{itemize}
\item Incoming wave: the incoming wave is given by $u^{inc}(\xb=(x,y))=\exp(ik x)$ with two frequencies, $k/(2\pi)=5$ and $20$;
\item Scatterer: we use a regular grid of $64^2$ point scatterers and their magnitudes are given by
$$
q_4(\xb=(x,y))\coloneqq
\begin{cases}
0.1,\quad \text{if} \quad \cos^2(2\pi x)+\cos^2(2\pi y)>0.5;  \\
0, \quad \quad \quad \quad \quad \quad \quad \quad  \text{otherwise}.
\end{cases}
$$
Figure \ref{fig:domain_example} shows a plot of $q_4$.
\item Domain decomposition: the number of subdomains is $N_s=4$, $16$ and $64$, and the overlap parameter is $\delta=1$ and $6$;
\end{itemize}

\begin{figure}[h!]
  \centering
  \begin{subfigure}{.4\linewidth}
\includegraphics[width=.8\textwidth]{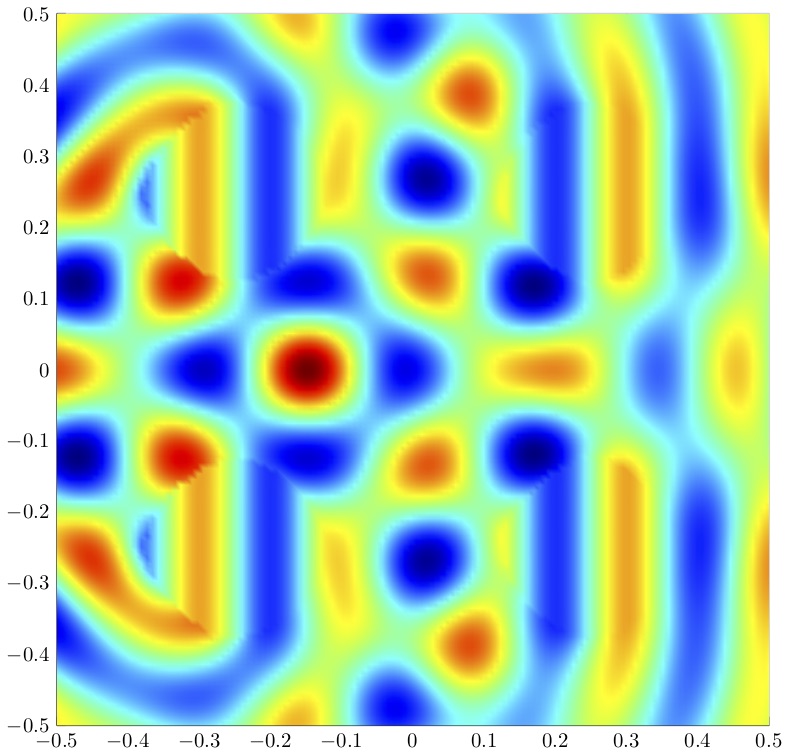}
\caption{$k/(2\pi)=5$}
\end{subfigure}
  \begin{subfigure}{.4\linewidth}
\includegraphics[width=.8\textwidth]{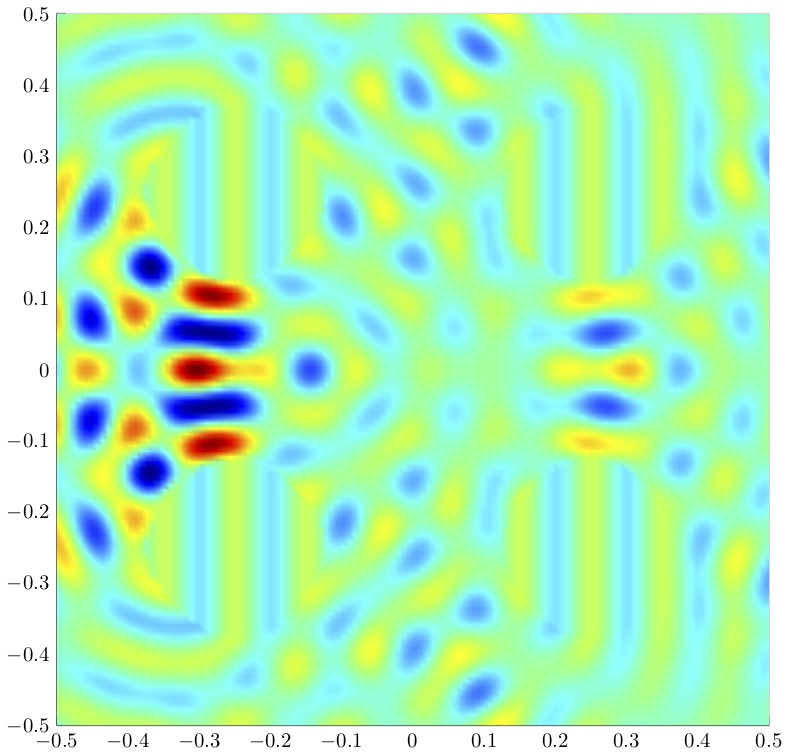}
\caption{$k/(2\pi)=10$}
\end{subfigure}
\begin{subfigure}{.4\linewidth}
\includegraphics[width=.8\textwidth]{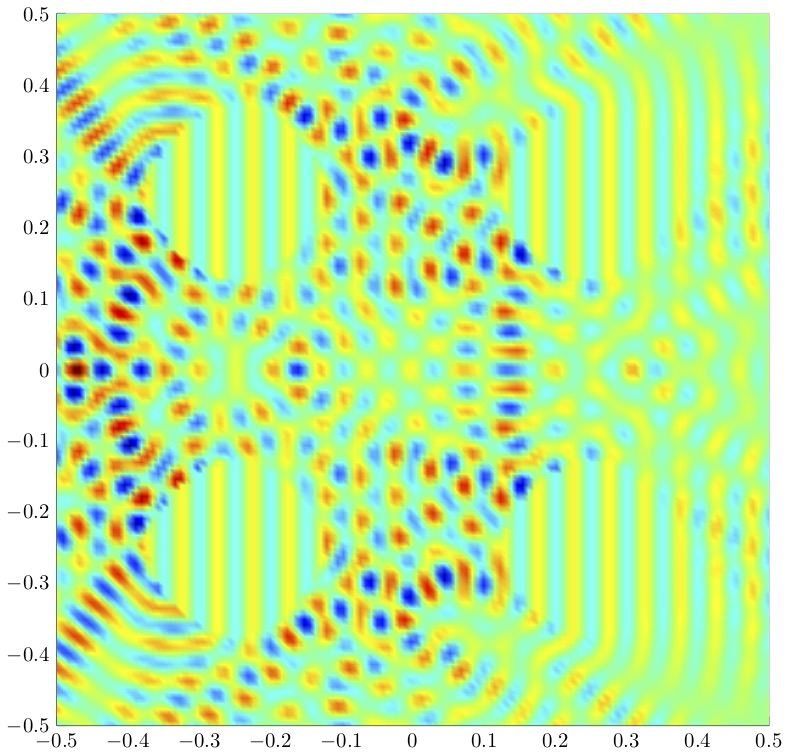}
\caption{$k/(2\pi)=20$}
\end{subfigure}
\begin{subfigure}{.4\linewidth}
\includegraphics[width=.8\textwidth]{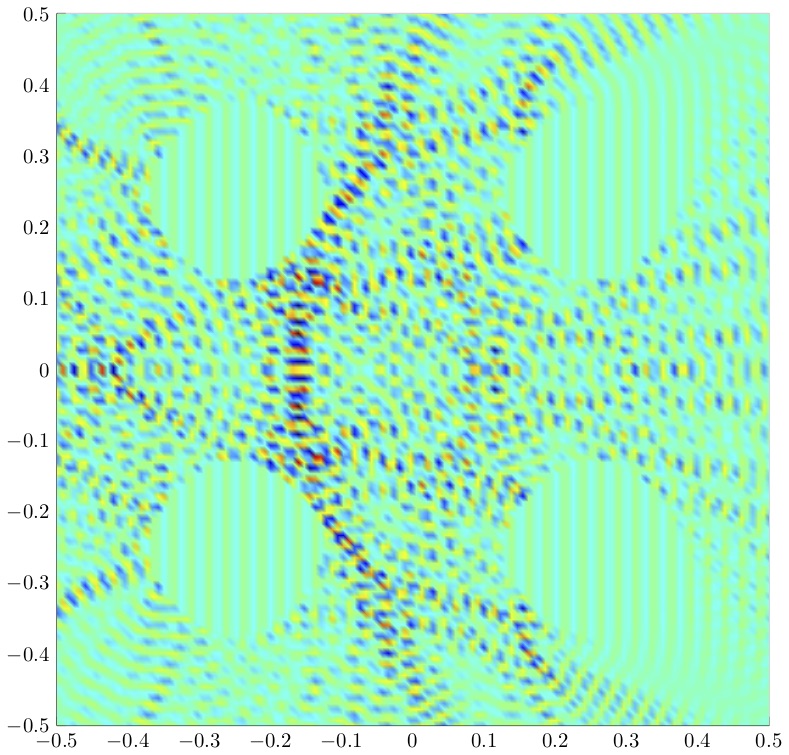}
\caption{$k/(2\pi)=40$}
\end{subfigure}  
\caption{The real part of the scattered field off of $q_4$ when the incident plane wave has incidence direction $\theta=(1,0)$ and wavenumber: (a) $k/(2\pi)=5$, (b) $k/(2\pi)=10$, (c) $k/(2\pi)=20$ and (d) $k/(2\pi)=40$.}
\label{fig:q4_scat_field}
\end{figure}

In Figure \ref{fig:q4_scat_field}, we plot the real part of $u^{scat}$ when the incident plane wave has incidence direction $\theta=(1,0)$ and wavenumber $k/(2\pi)=5$, $10$, $20$ and $40$.

In Table \ref{table:fwd_exp1}, we report the number of iterations necessary for the convergence of GMRES without preconditioner and with the precondtioners in Table \ref{table:precs}. In Table \ref{table:fwd_exp1_sup} of Appendix \ref{s:appendB}, we present the relative error $e_{rel}$ of the solution using GMRES with no preconditioners and with the domain decomposition based preconditioners.

In Figure \ref{fig:eig_example_f1_1} of Appendix \ref{s:appendB}, we plot the eigenvalues of (a) $\Ab$, (b) $\tilde{\Ab}^{-1}_{AS}\Ab$, (c) $\tilde{\Ab}^{-1}_{RAS}\Ab$ and (d) $\tilde{\Ab}^{-1}_{SRAS}\Ab$ when $N_s=16$ and $\delta=6$. The eigenvalues for $\tilde{\Ab}^{-1}_{RAS}\Ab$ are more clustered than the eigenvalues for $\Ab$, $\tilde{\Ab}^{-1}_{AS}\Ab$ and $\tilde{\Ab}^{-1}_{SRAS}\Ab$, justifying its best performance. In Figure \ref{fig:eig_example_f1_2} of Appendix \ref{s:appendB}, we present the eigenvalues of $\tilde{\Ab}^{-1}_{RAS}\Ab$ with (a)  $N_s=4$ and $\delta=1$, (b)  $N_s=4$ and $\delta=6$, (c)  $N_s=64$ and $\delta=1$ and (d)  $N_s=64$ and $\delta=6$. Note that the eigenvalues are more clustered for smaller $N_s$ and larger $\delta$.

\begin{table}
\center
\caption{(Experiment F.1) We present the number of iterations necessary for the convergence of GMRES without preconditioners (second column) and with domain decomposition preconditioning strategies (columns 5--8). The incoming plane wave has horizontal direction of propagation and wavenumbers $k/(2\pi)=5$ and $20$. We use a regular grid of $N=64^2$ point scatterers with magnitude given by the function $q_4$. The number of subdomains used is $N_s=4$, $16$, and $64$ and the overlap parameter is $\delta=1$ and $6$. For both cases, the number of GMRES iterations is reduced significantly, even in the case of having small subdomains. Notice that the full inverse of $\Ab$ requires inverting a $4096\times 4096$ matrix, whereas the preconditioner requires inverting matrices of size $70\times 70$. Also notice that the overlap significantly helps.}
\label{table:fwd_exp1}
\begin{tabular}{|c||c||c||c||*{4}{c|}}\hline
$k/2\pi$ & GMRES & $N_{s}$ & overlap  & AS & RAS & AHS & SRAS  \\\hline\hline
\multirow{6}{*}{5} & \multirow{6}{*}{97} &\multirow{2}{*}{4} & 1  & 18& 16& 15& 16\\\cline{4-8}
                           & &                          				        & 6    & 17& 15& 15& 15\\\cline{3-8}
			  & 				 &\multirow{2}{*}{16} & 1    & 54& 44& 44& 69\\\cline{4-8}

                           & &                          					& 6    & 30& 17& 17& 74\\\cline{3-8}
                           & 				&\multirow{2}{*}{64}  & 1    & 71& 44& 44& 68\\\cline{4-8}

                           & &                          					& 6    & 46& 18& 18& 74\\\hline
\multirow{6}{*}{20} & \multirow{6}{*}{293} &\multirow{2}{*}{4} & 1  & 43& 42& 42& 42\\\cline{4-8}

                              & &                          				& 6    & 43& 41& 41& 41\\\cline{3-8}
			     & 				&\multirow{2}{*}{16}   & 1 & 104& 74& 74& 130\\\cline{4-8}

                              & &                          				& 6    & 62& 41& 42& 143\\\cline{3-8}
                              & 				&\multirow{2}{*}{64}& 1 & 138& 76& 76& 129\\\cline{4-8}

                              & &                          				& 6    & 97& 41& 42& 143\\\hline
\end{tabular}
\end{table}

{\bf Summary:} The RAS and AHS have very similar performance and both work better than the AS and SRAS. The use of increasing overlap significantly improves the performance of the preconditioner. The increase of subdomains decreases the computational time for the construction of the preconditioners, but it also precludes an increase in the number of iterations for the convergence of the method.

{\bf Experiment F.2 -- Comparison of different partition geometries on DD forward preconditioning:} In this example we compare the choice of the geometry of the partition. We use two different geometries. In the first choice, called $\mathcal{G}_1$, the domain is subdivided into square domains of same size. In the second choice, called $\mathcal{G}_2$, the domain is divided in vertical bands. The following parameters are used: 
\begin{itemize}
\item Incoming wave: the incoming wave is given by $u^{inc}(\xb=(x,y))=\exp(ik x)$ with two frequencies, $k/(2\pi)=10$ and $40$;
\item Scatterer: we use a regular grid of $64^2$ point scatterers and their magnitudes are given by $q_4$; 
\item Domain decomposition: the number of subdomains is $N_s= 4$, $9$ and $16$; and the overlap parameter is $\delta=1$ and $3$;
\end{itemize}

In Table \ref{table:fwd_exp2}, we report the number of iterations necessary for the convergence of GMRES without preconditioner and using the AS, RAS and AHS preconditioners when using the two geometries for the partition of the domain.  In Table \ref{table:fwd_exp2_sup} of Appendix \ref{s:appendB}, we present the relative error $e_{rel}$.

\begin{table}[!htp]
\center
\caption{(Experiment F.2) We present the number of iterations necessary for the convergence of GMRES without preconditioners (second column) and with the  AS, RAS and AHS preconditioning strategies (columns 5--10). The incoming plane waves have wavenumbers $k/(2\pi)=10$ and $40$. We use a regular grid of $N=64^2$ point scatterers with magnitude given by the function $q_4$. $\mathcal{G}_1$ represents the partition composed of equal sized squares and $\mathcal{G}_2$ is the partition composed of vertical bands.The number of subdomains used is $N_s=4$, $9$, and $16$ and the overlap parameter is $\delta=1$ and $3$.}
\label{table:fwd_exp2}
\begin{tabular}{|c||c||c||c||*{6}{c|}}\hline
\multirow{2}{*}{$k/2\pi$} & \multirow{2}{*}{GMRES} & \multirow{2}{*}{$N_s$} & \multirow{2}{*}{$\delta$} & \multicolumn{2}{c|}{AS} & \multicolumn{2}{c|}{RAS} & \multicolumn{2}{c|}{AHS} \\\cline{5-10}
                                       &               &                                     &                                        & $\mathcal{G}_1$ & $\mathcal{G}_2$  & $\mathcal{G}_1$ & $\mathcal{G}_2$ & $\mathcal{G}_1$ & $\mathcal{G}_2$ \\\hline\hline
\multirow{6}{*}{10} & \multirow{6}{*}{164} & \multirow{2}{*}{4}   &   1  & 26 & 48 & 25 & 40 & 25 & 41 \\\cline{4-10}
                              &                                 &                               &   3  & 26 & 34 & 24 & 26 & 25 & 26 \\\cline{3-10}
                              &                                 & \multirow{2}{*}{9}   &   1  & 42 & 59 & 35 & 40 & 35 & 40 \\\cline{4-10}
                              &                                 &                               &   3  & 36 & 60 & 24 & 27 & 24 & 27 \\\cline{3-10}
                              &                                 & \multirow{2}{*}{16} &   1  & 72 & 99 & 57 & 63 & 58 & 63 \\\cline{4-10}
                              &                                 &                               &   3  & 60 & 60 & 32 & 36 & 32 & 36 \\\hline
\multirow{6}{*}{40} & \multirow{6}{*}{398} & \multirow{2}{*}{4}   &   1  & 71 & 152 & 67 & 136 & 67 & 138 \\\cline{4-10}
                              &                                 &                               &   3  & 71 & 85   & 67 & 68   & 67 & 68 \\\cline{3-10}
                              &                                 & \multirow{2}{*}{9}   &   1  & 123 & 176 & 125 & 135 & 124 & 139 \\\cline{4-10}
                              &                                 &                               &   3  & 92   & 155 & 77   & 79   & 77   & 79 \\\cline{3-10}
                              &                                 & \multirow{2}{*}{16} &   1  & 274 & 269 & 240 & 225 & 239 & 226 \\\cline{4-10}
                              &                                 &                               &   3  & 147 & 146 & 99   & 124 & 100 & 124 \\\hline
\end{tabular}
\end{table}

{\bf Summary:} For this case in particular, we were not able to experience any difference in performance between the two geometries used when both domains have almost the same amount of points. The methods behave similarly as in the previous example, with the RAS and AHS being faster and presenting similar results.

\begin{figure} 
  \centering
  \begin{subfigure}{.4\linewidth}
\includegraphics[width=1\textwidth]{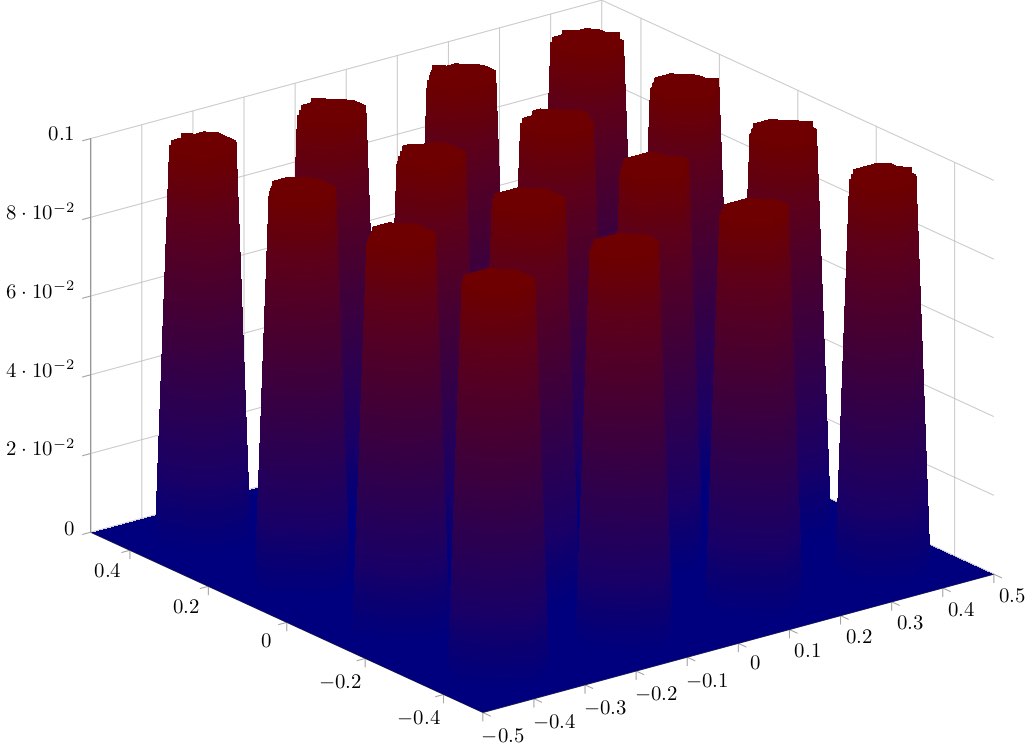}
\caption{isometric view of $q_{16}$}
\end{subfigure}
\begin{subfigure}{.4\linewidth}
\includegraphics[width=.8\textwidth]{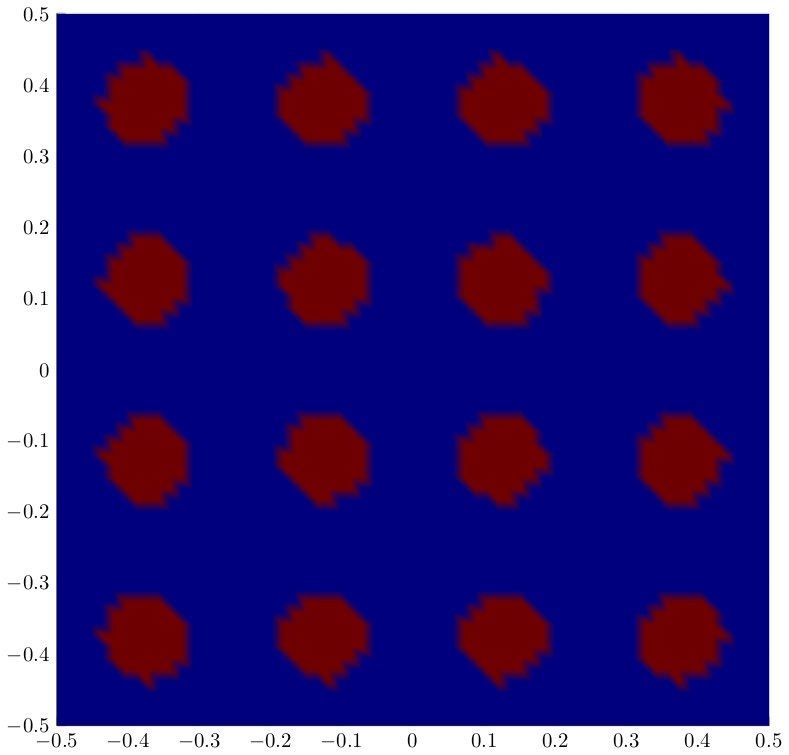}
\caption{top view of domain of $q_{16}$}
\end{subfigure}
\caption{The (a) isometric view and (b) top view of the domain $q_{16}$ used in the Experiment F.3.}
\label{fig:domain_example_16}
\end{figure}

{\bf Experiment F.3 -- Comparison of different scatterers on DD forward preconditioning:} In this example we compare the results for two different scatterers $q_4$ and $q_{16}$. The scatterer $q_{16}$ is given by
$$
q_{16}(\xb=(x,y))\coloneqq
\begin{cases}
0.1,\quad \text{if} \quad \cos^2(4\pi x)+\cos^2(4\pi y)>0.5;  \\
0, \quad \quad \quad \quad \quad \quad \quad \quad \text{otherwise}.
\end{cases}
$$
Figure \ref{fig:domain_example_16} shows a plot of $q_{16}$.
The following parameters are used: 
\begin{itemize}
\item Incoming wave: the incoming wave is given by $u^{inc}(\xb=(x,y))=\exp(ik x)$ with two frequencies, $k/(2\pi)=10$ and $40$;
\item Scatterers: we use a regular grid of $128^2$ scatterer points. The scatterers magnitudes are given by $q_4$ and $q_{16}$; 
\item Domain decomposition: the number of subdomains is $N_s= 4$ and $16$, and the overlap parameter is $\delta=1$ and $16$.
\end{itemize}

\begin{figure}[h!]
  \centering
  \begin{subfigure}{.4\linewidth}
\includegraphics[width=.8\textwidth]{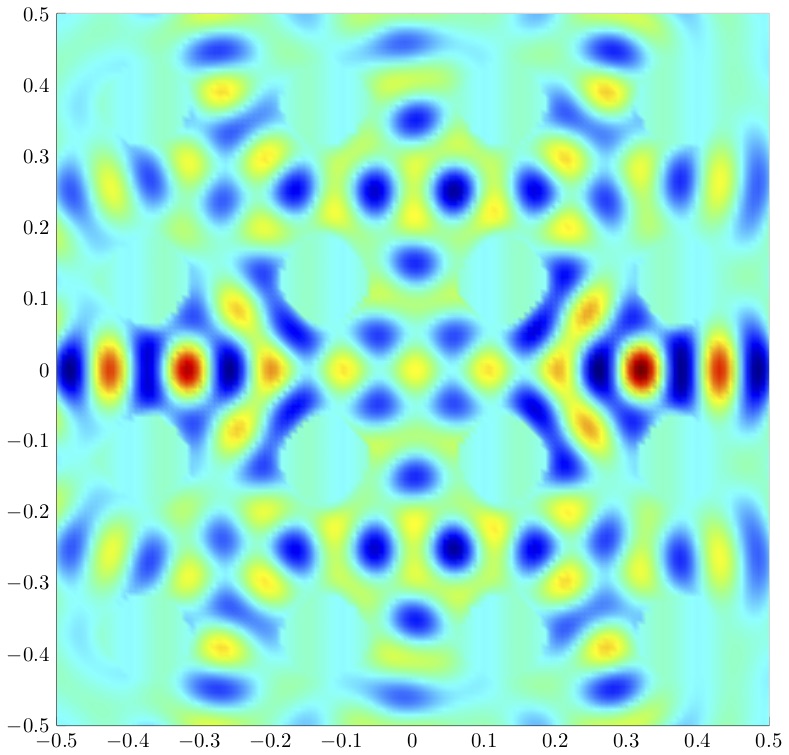}
\caption{$q_{16}$ and $k/(2\pi)=10$}
\end{subfigure}
\begin{subfigure}{.4\linewidth}
\includegraphics[width=.8\textwidth]{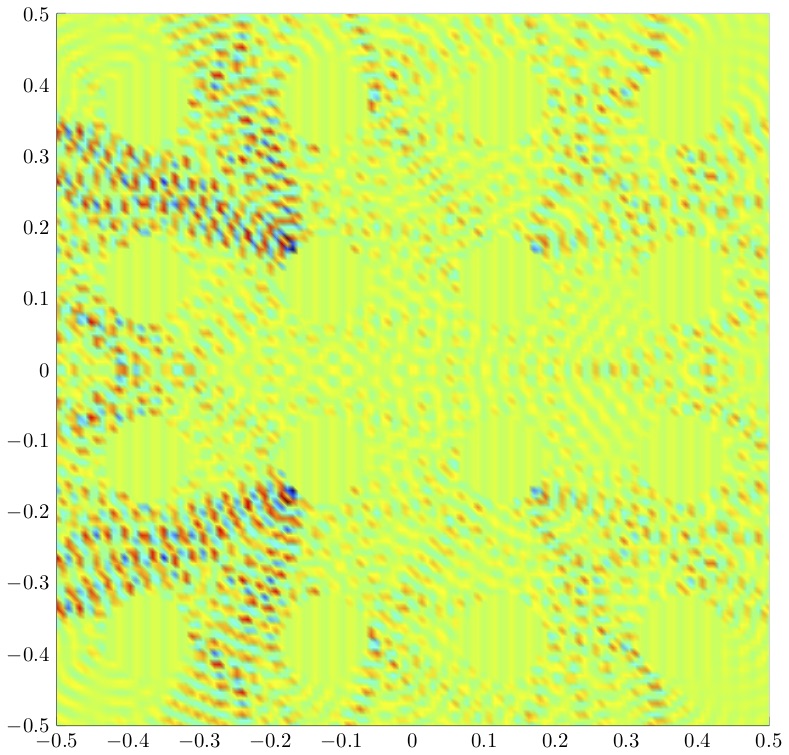}
\caption{$q_{16}$ and $k/(2\pi)=40$}
\end{subfigure}
\caption{The real part of the scattered field off of $q_{16}$ when the incident plane wave has incidence direction $\theta=(1,0)$ and wavenumber: (a) $k/(2\pi)=10$ and (b) $k/(2\pi)=40$.}
\label{fig:q16_scat_field}
\end{figure}

In Figure \ref{fig:q16_scat_field}, we show the plot of the real part of the scattered field off of $q_{16}$ when the incident plane wave has incidence direction $\theta=(1,0)$ and wavenumber $k/(2\pi)=10$ and $40$.

In Table \ref{table:fwd_exp3}, we report the number of iterations necessary for the convergence of GMRES without preconditioner and with the AS, RAS and AHS preconditioners when the scatterer is $q_4$ and $q_{16}$. In Table \ref{table:fwd_exp3_sup}, we present the relative error of the solution $e_{rel}$.

\begin{table}[!htp]	      
\caption{(Experiment F.3) We present the number of iterations necessary for the convergence of GMRES without preconditioner (columns 2--3) and with the AS, RAS and AHS preconditioners (columns 6--11). The incoming plane waves have wavenumbers $k/(2\pi)=10$ and $40$. We use $N=64^2$ scatterers points with magnitude given by the functions $q_4$ and $q_{16}$. The number of subdomains used is $N_s=4$ and $16$, and the overlap parameter is $\delta=1$ and $16$.}\label{table:fwd_exp3}
\center
\begin{tabular}{|c||c|c||c||c||*{6}{c|}}\hline
\multirow{2}{*}{$k/2\pi$} & \multicolumn{2}{c||}{GMRES} & \multirow{2}{*}{$N_s$} & \multirow{2}{*}{$\delta$} & \multicolumn{2}{c|}{AS} & \multicolumn{2}{c|}{RAS} & \multicolumn{2}{c|}{AHS} \\\cline{2-3}\cline{6-11}
                                       &  $q_4$       & $q_{16}$         &                                     &                                        & $q_4$ & $q_{16}$          & $q_4$ & $q_{16}$            & $q_4$ & $q_{16}$ \\\hline\hline
\multirow{4}{*}{10}          & \multirow{4}{*}{211}  & \multirow{4}{*}{259} &\multirow{2}{*}{4}   &   1                             & 23 & 45 & 21 & 42 & 21 & 42 \\\cline{5-11}
                                       &                                  &                                 &                              &   16                            & 23 & 34 & 21 & 22 & 21 & 22 \\\cline{4-11}
                                       &                                  &                                 &\multirow{2}{*}{16} &   1                             & 69 & 81 & 60 & 77 & 60 & 77 \\\cline{5-11}
                                       &                                  &                                 &                              &   16                           & 23 & 85 & 21 & 41 & 21 & 41 \\\hline
\multirow{4}{*}{40}          & \multirow{4}{*}{674}  &  \multirow{4}{*}{926}&\multirow{2}{*}{4}   &   1                             & 64 & 153 & 61 & 152 & 61 & 151 \\\cline{5-11}
                                       &                                  &                                 &                              &   16                            & 64 & 73 & 61 & 62 & 61 & 62 \\\cline{4-11}
                                       &                                  &                                 &\multirow{2}{*}{16} &   1                             & 134 & 279 & 100 & 273 & 100 & 271 \\\cline{5-11}
                                       &                                  &                                 &                              &   16                            & 65 & 207 & 61 & 152 & 61 & 152 \\\hline
\end{tabular}
\end{table}

{\bf Summary:} As expected, due to the increasing of the effect of multiple scattering, the number of iterations necessary for GMRES to converge is higher for the domain $q_{16}$. The RAS and AHS remain the best performing methods.

{\bf Experiment F.4 -- Scalability of the DD preconditioning on forward problem:} This example aims to check the scalability of the methods presented in Table \ref{table:precs} as we increase the number of scatterers. We fix the number of subdomains and increase the number of points in the grid. The following parameters are used: 
\begin{itemize}
\item Incoming wave: the incoming wave is given by $u^{inc}(\xb=(x,y))=\exp(ik x)$ with three frequencies, $k/(2\pi)=10$, $20$ and $40$;
\item Scatterer: we use regular grids of $N=64^2$, $128^2$ and $256^2$ scatterer points. The scatterer points magnitudes are given by $q_4$;
\item Domain decomposition: the number of subdomains is $N_s=16$, we use the overlap parameter $\delta=3$ for the grid with $64^2$ scatterer points, $\delta=6$ for $128^2$ scatterer points and $\delta=12$ for $256^2$ scatterer points;
\end{itemize}

In Table \ref{table:fwd_exp4}, we report the number of iterations for the convergence of GMRES without preconditioner and with the preconditioners AS, RAS and AHS. We present the relative error of the solution $e_{rel}$ in Table \ref{table:fwd_exp4_sup} on Appendix \ref{s:appendB}.

\begin{table}[!htp]	
\center
\caption{(Experiment F.4) We present the number of iterations for the convergence of the GMRES without preconditioner (column 3) and with the domain decomposition preconditioners in Table \ref{table:precs} (columns 5--8). The incoming plane wave has horizontal direction of propagation and frequencies $k/(2\pi)=10$, $20$ and $40$. We use regular grids of $N=64^2$, $128^2$ and $256^2$ point scatterers with magnitude given by the function $q_4$. The number of subdomains used is $N_s=16$ and the overlapping parameter is $\delta=3$, $6$ and $12$ for $N=64^2$, $N=128^2$ and $N=256^2$ respectively.}\label{table:fwd_exp4}
\begin{tabular}{|c||c||c||c||*{4}{c|}}\hline
$k/2\pi$ & $\sqrt{N}$ & GMRES & $\delta$ & AS & RAS & AHS & SRAS \\\hline\hline
\multirow{3}{*}{10} &                             64 & 164 & 3   & 60& 32& 32& 105  \\\cline{2-8}
                              & \multirow{1}{*}{128} & \multirow{1}{*}{211}  & 6   & 55& 24& 24& 116\\\cline{2-8}
                              & \multirow{1}{*}{256} & \multirow{1}{*}{277}  & 12 & 52& 21& 21& 121\\\hline
\multirow{3}{*}{20} &                             64 & 293& 3   & 78& 45& 45& 143\\\cline{2-8}
                              & \multirow{1}{*}{128} & \multirow{1}{*}{413}   & 6   & 71& 31& 31& 152\\\cline{2-8}
                              & \multirow{1}{*}{256} & \multirow{1}{*}{526}  & 12 & 65& 27& 27& 169\\\hline
\multirow{3}{*}{40} &                             64 & 398 & 3   & 147& 99& 100& 286\\\cline{2-8}
                              & \multirow{1}{*}{128} & \multirow{1}{*}{674}  & 6   & 112& 62& 61& 200\\\cline{2-8}
                              & \multirow{1}{*}{256} & \multirow{1}{*}{957}  & 12 & 101& 45& 45& 230\\\hline
\end{tabular}
\label{table:gmres_precs_radius_AS}
\end{table}

{\bf Summary:} The results show that to have the same accuracy with unpreconditioned GMRES, when keeping the number of subdomains constant, we need to increase the  overlap parameter as the number of scatterers in the grid increases. 

{\bf Experiment F.5 -- RC preconditioner on forward problem:} This example aims to show the performance of the RC preconditioner to solve Equation \eqref{eq:forward_lin_system} when the correction is applied in the RAS preconditioner. We analyze the effects of the number of subdomains, size of the overlap and number of singular values used for the correction at different wavenumbers. The following parameters are used: 
\begin{itemize}
\item Incoming wave: the incoming wave is given by $u^{inc}(\xb=(x,y))=\exp(ik x)$ with two frequencies, $k/(2\pi)=5$ and $20$;
\item Scatterer: we use a regular grid of $64^2$ scatterer points. The scatterer points magnitudes are given by $q_4$; 
\item Domain decomposition: number of subdomains is $N_s= 4$ and $16$, and the overlap parameter is $\delta=1$ and $8$;
\item RC preconditioner: we use the RAS preconditioner $\tilde{\Ab}^{-1}_{RAS}$ to construct the preconditioner $\tilde{\Ab}^{-1}_{RC}$. We choose $N_\lambda=20$, $40$, $60$ and $80$ for wavenumber $k/(2\pi)=5$ and $N_\lambda=40$, $80$, $120$ and $160$ for wavenumber $k/(2\pi)=20$.
\end{itemize}

In Table \ref{table:fwd_exp5}, we reportt the number of iterations for the convergence of GMRES without preconditioner, with the RAS preconditioner and with the RC preconditioner with different numbers of singular values for the correction.  In Table \ref{table:fwd_exp5_sup} in Appendix \ref{s:appendB}, we present the relative error of the solution $e_{rel}$.

\begin{table}[!htp]	    
\caption{(Experiment F.5) We present the number of iterations for the convergence of GMRES without preconditioner (column 2), with the RAS preconditioner (column 5) and with the RC preconditioner using $N_\lambda=20$, $40$, $60$ and $80$ at $k=5/(2\pi)$ and $N_\lambda=40$, $80$, $120$ and $160$ at $k/(2\pi)=20$ (columns 6--9). The incoming plane waves have wavenumbers $k/(2\pi)=5$ and $20$ with incidence direction $(1,0)$. We use a regular grid of $N=64^2$ point scatterers with magnitude given by the function $q_4$. The number of subdomains used is $N_s=4$ and $16$, and the overlap parameter is $\delta=1$ and $8$.}\label{table:fwd_exp5}
\center  
\begin{tabular}{|c||c||c|c||c||c|c|c|c|}\hline
\multirow{2}{*}{$k/(2\pi)$} & \multirow{2}{*}{GMRES} & \multirow{2}{*}{$N_s$} & \multirow{2}{*}{$\delta$} & \multirow{2}{*}{RAS} & \multicolumn{4}{c|}{RC -- $N_\lambda$}  \\ \cline{6-9}
                                         &				      &                                      &                                       &                                    &  20 & 40 & 60 & 80  \\\hline\hline
\multirow{4}{*}{ 5 } 		& \multirow{4}{*}{ 97 }         &\multirow{2}{*}{4}           &   1			            & 15   				& 9 & 5 & 3 & 3 \\\cline{4-9}
 &				      &                                      &   8          			    & 15				& 9& 5& 3& 3 \\\cline{3-9}
 &                                        &\multirow{2}{*}{16}         &   1          			     & 44				& 29& 24& 20& 17 \\\cline{4-9}
 &                                        &                                      &   8          			    & 15				& 9& 5& 3& 3 \\\hline \hline
 \multicolumn{5}{c||}{} & \multicolumn{4}{c|}{RC -- $N_\lambda$} \\\cline{6-9}
 \multicolumn{5}{c||}{}                         & 40 & 80 & 120 & 160  \\ \hline \hline                                       
\multirow{4}{*}{20 } &\multirow{4}{*}{ 293 }		&\multirow{2}{*}{4}  		&   1         				 & 42				& 18& 8& 4& 2  \\\cline{4-9}
 &                 	                 &                             		 &   8         			 & 40				& 18& 8& 4& 2 \\\cline{3-9}
 &        		                         &\multirow{2}{*}{16} 		&   1          			& 74				& 52& 38& 30& 22 \\\cline{4-9}
 &	                                  &                            		  &   8         			 & 41				& 18& 8& 4& 2 \\\hline  
 \end{tabular}
\end{table}

{\bf Summary:} The number of singular values used for the construction of the RC preconditioner is dependent on the number of subdomains and the amount of overlap used for the RAS preconditioner. It is also dependent on the wavenumber of the incoming wave. From the results, we have that as the wavenumber of the incoming wave increases,  a higher rank correction is required for obtaining better performance for the RC preconditioner.

\subsection{Conclusions} The preconditioning methods tested in this section provide improvements on the speed-up of convergence of the standard GMRES approach decreasing its computational cost. Among the domain decomposition preconditioners, the RAS and AHS readily outperform the AS and SRAS in all of our examples. As expected, the number of iterations decreases when the partition overlap is larger, and increases if the number of subdomains increases too much. The RC preconditioner improve even further the speed of convergence of the domain decomposition methods requiring a total computational complexity of $\mathcal{O}(N\log(N_\lambda)+N_\lambda^3+NN_\lambda+N+(N/N_s)^3)$. The methods presented are easily translated to higher dimensions and they are scalable with the increase of the number of scatterers.

%% file: inv_prec.tex
\section{Preconditioning of the inverse problem}\label{s:method_inverse}
To obtain the update of the domain at each step of the Gauss-Newton method, we must solve Equation \eqref{eq:Hessian_prob_1} using GMRES. Since $\Hb=\Jb^\ast\Jb+\beta\Ib$, and $\Jb$ involves $\Ab^{-1}$ \eqref{eq:Joperator}, the cost of a Hessian matrix-vector multiplication is that of two forward solves per direction. So it is imperative to reduce the number of GMRES iterations. To this end, we propose a preconditioner $\tilde{\Hb}_{RC}$ based on $\tilde{\Ab}_{RC}^{-1}$.

How can we use the $\tilde{\Hb}_{RC}$ as a preconditioner? It can be used as an H-matrix approximation and consequently apply it for solving the inverse problem step directly with very low accuracy; it can be used to construct high-rank approximations, and many others. We intend to explore some of these algorithmic variants in our future work. Our focus here is to show that $\tilde{\Hb}_{RC}$ is a good preconditioner using brute force factorization.

%

\subsection{Preconditioning using an approximation of the inverse forward operator}


To speed-up the iterative method for solving Equation \eqref{eq:Hessian_prob_1}, we intend to use a preconditioner that approximates $\Hb^{-1}$. Our first preconditioner is $\tilde{\Hb}=\tilde{\Jb}^\ast \tilde{\Jb}+\beta \Ib$, where 
\begin{equation*}
\tilde{\Jb} = -k^2 \Gb_r \Ub^{tot} + k^4 \Gb_r\Qb \tilde{\Ab}^{-1}\Gb\Ub^{tot}, \label{eq:frechet_der_one_dir_DD}
\end{equation*}
and $\tilde{\Ab}^{-1}$ is a domain decomposition preconditioner from Table \ref{table:precs}. Unfortunately, we found out that the domain decomposition preconditioners without rank correction do not work that well for the Hessian. However, they are effective when combined with the rank correction.



Next, we use the RC precondtidioner $\tilde{\Ab}_{RC}^{-1}$ as an approximation of $\Ab^{-1}$. The approximation of the matrix $\Jb$ is given by
\begin{equation*}
\tilde{\Jb} = -k^2 \Gb_r \Ub^{tot} +k^4 \Gb_r\Qb \tilde{\Ab}_{RC}^{-1}\Gb\Ub^{tot}. \label{eq:frechet_der_one_dir_RC}
\end{equation*}
The quality of the approximation of $\Jb$ by $\tilde{\Jb}$ depends on the number of eigenvalues $N_\lambda$ chosen for the rank correction.

Consider the singular value decomposition $\Ub\Sb\Vb^\ast=\left(\Ab-\tilde{\Ab}\right)$ and the submatrices $\Ub_{N_\lambda}=\Ub(:,$$1$:$N_\lambda)$, $\Sb_{N_\lambda}=\Sb($$1$:$N_\lambda,$$1$:$N_\lambda)$ and $\Vb_{N_\lambda}=\Vb(:,$$1$:$N_\lambda)$. According to \cite{GunnarPNAS}, setting a tolerance $\epsilon>0$, we can obtain $N_\lambda$ such that
\begin{equation}
\|\tilde{\Ab}-\Ab-\Ub_{N_\lambda}\Sb_{N_\lambda}\Vb_{N_\lambda}^\ast\|\leq\epsilon. \label{eq:MFdif}
\end{equation}

Since $\tilde{\Ab}^{-1}_{RC}=(\tilde{\Ab}-\Ub_{N_\lambda}\Sb_{N_\lambda}\Vb_{N_\lambda}^\ast)^{-1}$, we have
\begin{equation}
\tilde{\Ab}-\tilde{\Ab}_{RC}=\Ub_{N_\lambda}\Sb_{N_\lambda}\Vb_{N_\lambda}^\ast. \label{eq:invTFdif}
\end{equation}
From \eqref{eq:MFdif} and \eqref{eq:invTFdif}, we obtain
\begin{equation}
\|\tilde{\Ab}_{RC}-\Ab\|\leq\epsilon. \label{eq:FminusinvT}
\end{equation}

The norm of the difference between $\Jb$ and its approximation $\tilde{\Jb}_{RC}$ is given by
\begin{eqnarray}
\|\Jb-\tilde{\Jb}_{RC}\|&\leq& k^4\|\Gb_r\Qb(\Ab^{-1}-\tilde{\Ab}_{RC}^{-1})\Gb\Ub^{tot}\| \\ \nonumber
                         &\leq& k^4 \|\Gb_r\| \|\Qb\| \|\Ab^{-1}-\tilde{\Ab}_{RC}^{-1}\| \|\Gb\| \|\Ub^{tot}\|. \nonumber
\end{eqnarray}                               
We have that $\|\Gb\|$, $\|\Gb_d\|$ and $\|\Ub^{tot}\|$ are bounded by a constant depending on $k$ and $\|\Qb\|\leq \|q\|_{\infty}$. 

Since $\left(\Ab^{-1}-\tilde{\Ab}_{RC}^{-1}\right)=\tilde{\Ab}_{RC}^{-1}\left(\tilde{\Ab}_{RC}-\Ab\right)\Ab^{-1}$, and using properties of matrices norms, we get
\begin{equation}
\|\Ab^{-1}-\tilde{\Ab}_{RC}^{-1}\|\leq\|\tilde{\Ab}_{RC}^{-1}\| \|\Ab-\tilde{\Ab}_{RC}\| \|\Ab^{-1}\|\leq \tilde{C}(k,q,N_\lambda)\|\Ab-\tilde{\Ab}_{RC}\|, 
\label{eq:ft_bound}
\end{equation}
where $\tilde{C}(k,q,N_\lambda)$ is a constant that depends on $k$, $\|q\|_\infty$ and $N_\lambda$.

Using the bounds of the norm of the matrices and \eqref{eq:ft_bound}, we obtain
\begin{equation*}
\|\Jb-\tilde{\Jb}_{RC}\|\leq C(k,q,N_\lambda) \|\Ab-\tilde{\Ab}_{RC} \|,
\end{equation*}
where we combined the constants and reuse $C(k,q,N_\lambda)$ to denote the final constant.

We have the preconditioner $\tilde{\Hb}_{RC}=\tilde{\Jb}_{RC}^\ast\tilde{\Jb}_{RC}+\beta\Ib$. We refer to this preconditioner, as the {\bf Hessian rank corrected preconditioner} (HRC preconditioner). 



\subsection{Numerical Experiments for Inverse Preconditioning}
The aim of the first two numerical experiments is to test the effect of the overlap parameter and the number of subdomains in the choice of the number of singular values used for the construction of HRC preconditioner. In the third experiment, we check the scalability of the HRC preconditioner with increasing number of scatterers. In the fourth experiment, we compare the HRC preconditioner with a {\bf low-rank preconditioner} of the Gauss-Newton Hessian, which we term the LR preconditioner. In our last experiment, we use full aperture data from several incoming directions at multiple frequencies to obtain a full reconstruction of the scatterer using the recursive linearization algorithm (RLA) \cite{Chen}. A list of the experiments in this section with their description and results is provided in Table \ref{table:ip_experiments}. For a summary of the RLA, please see Appendix \ref{s:appendA}, and for a more detailed exposition see \cite{Chen}.

The scattered data $\db$ is measured at $N_r$ receivers at $R\left(\cos(2m\pi/N_r),\sin(2m\pi/N_r)\right)$, with $m=0,1,\ldots,N_r-1$. The data is generated by $N_d$ incoming incident waves with direction of propagation $\theta_j=\left(\cos(2j\pi/N_d),\sin(2j\pi/N_d)\right)$, for $j=0,1,\ldots,N_d-1$. For Examples I.1-4 the initial guess $\tilde{q}$ is chosen to be a tiny perturbation of $q$ given by
\begin{equation*}
\tilde{q}(x)=q(x)+\frac{10^{-2}\|q\|\mathcal{N}(0,1)}{kN^2},
\end{equation*}
where $\mathcal{N}(0,1)$ is a uniformly distributed random number in $\left[0,1\right]$. With $\tilde{q}$ we are able to calculate $\Fb(\tilde{q})$. For Example I.5, the initial guess is $\tilde{q}\equiv 0$.

To deal with the ill-posedness of this system, we scale and regularize the Gauss-Newton Hessian as follows: we estimate the maximum singular value $\sigma_{max}(\Jb)$ and set $\Hb = \sigma_{max}^{-2}\Jb^\ast\Jb+\beta\Ib$. For Examples I.1-4, we use $\beta=10^{-6}$. In this way, we control the conditioning of $\Hb$. For Example I.5, we use a different parameter that depends on the wavenumber, more details follow in the specific example.

We use GMRES with no restart to solve Equation \eqref{eq:Hessian_prob_1}. We also tried CG, however, the number of iterations to obtain the same accuracy on the solution with no preconditioner is much higher than GMRES. We believe that this behavior is justified by the conditioning of the matrix, since GMRES minimizes the residual using Givens rotations which is more stable than the Gram-Schmidt process used by CG. The behavior of the two methods becomes very similar when using the HRC preconditioner with increasing correction, therefore, we do not report the results of the experiments using CG here.

In all experiments, we use the RAS domain decomposition preconditioner $\tilde{\Ab}^{-1}_{RAS}$ to construct the RC preconditioner $\tilde{\Ab}^{-1}_{RC}$. The RAS preconditioner was chosen because it presented the best results among the domain decomposition preconditioners in the forward problem experiments.

For this section we define the relative error of the iterative solution with respect to the direct method solution $e_{rel}\coloneqq\|\delta q_{\mathrm{GMRES}}-\delta q_{\mathrm{LU}}\|/\|\delta q_{\mathrm{LU}}\|$, where $\delta q_{\mathrm{GMRES}}$ is the solution obtained by GMRES and $\delta q_{\mathrm{LU}}$ is the solution obtained by the LU direct solver. The LU solution is obtained by solving \eqref{eq:Hessian_prob_1} using the backslash in MATLAB.

\begin{table}
\caption{List of inverse problem experiments.}\label{table:ip_experiments}
{\small
\begin{center}
\begin{tabular}{|c|l|c|c|}
\hline
Experiment & Description & Tables & Figures\\
\hline\hline
\multirow{2}{*}{I.1} & Influence of the overlapping on inverse                    & \multirow{2}{*}{\ref{tab:error_exp_1_2}}   & \multirow{2}{*}{\ref{fig:errorH_inv_ol}, \ref{fig:error_inv_ol_ft}} \\
     &preconditioning										      &									&	\\\hline
\multirow{2}{*}{I.2} & Influence of the number of subdomains on              & \multirow{2}{*}{\ref{tab:error_exp_1_2}}  & \multirow{2}{*}{\ref{fig:errorH_inv_sd}, \ref{fig:error_inv_sd_ft}} \\
			     & inverse preconditioning					      & 							      & \\\hline 
I.3 & Scalability of the inverse preconditioning 		           				& \ref{tab:inv_scal_it}, \ref{tab:error_tf}, \ref{tab:inv_scal_error} & X \\\hline
I.4 & Comparison with low-rank inverse preconditioner                                                    & X  & \ref{fig:Hcomparison_I4}, \ref{fig:Heig_example_I4_k5}, \ref{fig:Heig_example_I4_k20} \\\hline 
I.5 & Nonlinear inverse problem 									& \ref{tab:inv_full} & \ref{fig:domain_example_I4}, \ref{fig:inv_full_sol}, \ref{fig:inv_full_bench}\\\hline
\end{tabular}
\end{center}
}
\end{table}

{\bf Experiment I.1 -- Influence of the overlapping on inverse preconditioning:} We analyze the performance of the HRC preconditioner when we change the size of the overlapping when keeping the number of subdomains constant. The following parameters are used: 
\begin{itemize}
\item Incoming waves: the incoming waves are given by the wave $u^{inc}(\xb=(x,y))=\exp(ik \xb\cdot\theta_j)$ with $k/(2\pi)=5$, $10$ and $20$, and $\theta_j$ prescribed as in the beginning of the subsection with $j=0,\ldots,7$;
\item Receivers: the receivers are located at $\xb_r=(x_r,y_r)=0.8(\cos(2\pi r/N_r), \sin(2\pi r/N_r))$, with $r=1,\ldots,N_r$ and $N_r=2,000$;
\item Scatterer: the scatterers points are distributed in an uniform grid of points with $N=64^2$ scatterers and their magnitude is given by the function $q_4$;
\item Domain decomposition: the number of subdomains is fixed at $N_s=16$ and the overlap parameter is $\delta=3$, $6$ and $9$;
\item HRC preconditioner: we use for the low-rank correction $N_\lambda=10+10m$, $m=0,\ldots,M_\lambda$, where $M_\lambda$ differs for each problem and it is the maximum number needed to obtain the prescribed error.
\end{itemize}

\begin{figure}
\footnotesize
\centering
\begin{subfigure}{.45\textwidth}
\begin{tikzpicture}[scale=0.7]
\begin{semilogyaxis}[xmin=10, xmax=160,domain=1:160,ymin=1, ymax=600,domain=1:600,
    xlabel={$N_\lambda$},
    ylabel={Iterations},
    grid=major,
    legend pos=north east,
    legend entries={$\delta=3$,$\delta=6$,$\delta=9$},
]
\addplot coordinates {
(10,545) (20,355) (30,265) (40,210) (50,130) (60,57)
(70,39) (80,33) (90,29) (100,29) (110,26) (120,17)
(130,14) (140,10) (150,9) (160,9)
};
\addplot coordinates {
(10,468) (20,202) (30,123) (40,51) (50,24) (60,13)
(70,10) (80,9) (90,8)
};
\addplot coordinates {
(10,322) (20,226) (30,162) (40,42) (50,10) (60,6)
};
\end{semilogyaxis}
\end{tikzpicture}
\caption{Iterations at $k/(2\pi)=5$}\label{fig:errorH_inv_ol_1}
\end{subfigure}
\begin{subfigure}{.45\textwidth}
\begin{tikzpicture}[scale=0.7]
\begin{semilogyaxis}[xmin=10, xmax=200,domain=1:200,ymin=1, ymax=800,domain=1:800,
    xlabel={$N_\lambda$},
    ylabel={Iterations},
    grid=major,
    legend pos=north east,
    legend entries={$\delta=3$,$\delta=6$,$\delta=9$},
]
\addplot coordinates {
(10,793) (20,282) (30,240) (40,145) (50,148) (60,90)
(70,84) (80,46) (90,35) (100,31) (110,27) (120,20)
(130,18) (140,15) (150,14) (160,9) (170,8) (180,7)
(190,6) (200,6)
};
\addplot coordinates {
(10,667) (20,261) (30,124) (40,48) (50,29) (60,17)
(70,15) (80,14) (90,11) (100,10) (110,9) (120,6)
};
\addplot coordinates {
(10,601) (20,422) (30,165) (40,110) (50,57) (60,25)
(70,11) (80,6) (90,4)
};
\end{semilogyaxis}
\end{tikzpicture}
\caption{Iterations at $k/(2\pi)=10$}\label{fig:errorH_inv_ol_2}
\end{subfigure}
\begin{subfigure}{.45\textwidth}
\begin{tikzpicture}[scale=0.7]
\begin{semilogyaxis}[xmin=10, xmax=280,domain=1:280,ymin=1, ymax=1300,domain=1:1300,
    xlabel={$N_\lambda$},
    ylabel={Iterations},
    grid=major,
    legend pos=north east,
    legend entries={$\delta=3$,$\delta=6$,$\delta=9$},
]
\addplot coordinates {
(10,1265) (20,1181) (30,1164) (40,979) (50,683) (60,350)
(70,329) (80,305) (90,297) (100,255) (110,254) (120,182)
(130,176) (140,165) (150,144) (160,123) (170,112) (180,105)
(190,93) (200,46) (210,36) (220,26) (230,21) (240,19)
(250,17) (260,16) (270,14) (280,12)
};
\addplot coordinates {
(10,1192) (20,967) (30,849) (40,787) (50,526) (60,214)
(70,108) (80,43) (90,22) (100,21) (110,19) (120,15)
(130,12) (140,11) (150,11) (160,11) (170,12) (180,11)
};
\addplot coordinates {
(10,1251) (20,985) (30,996) (40,804) (50,577) (60,318)
(70,260) (80,159) (90,119) (100,33) (110,20) (120,16)
(130,10) 
};
\end{semilogyaxis}
\end{tikzpicture}
\caption{Iterations at $k/(2\pi)=20$}\label{fig:errorH_inv_ol_3}
\end{subfigure}
\caption{(Experiment I.1) Plot of the number of iterations necessary to obtain $e_{rel}\approx\mathcal{O}(10^{-4})$ using the HRC preconditioner for the GMRES with different $N_\lambda$ at wavenumbers: (a) $k/(2\pi)=5$, (b) $k/(2\pi)=10$ and (c) $k/(2\pi)=20$. Each line represents the number of iterations obtained when the preconditioner is obtained using $\delta=3$, $6$ and $9$. The number of subdomains is $N_s=16$. }\label{fig:errorH_inv_ol}
\end{figure}
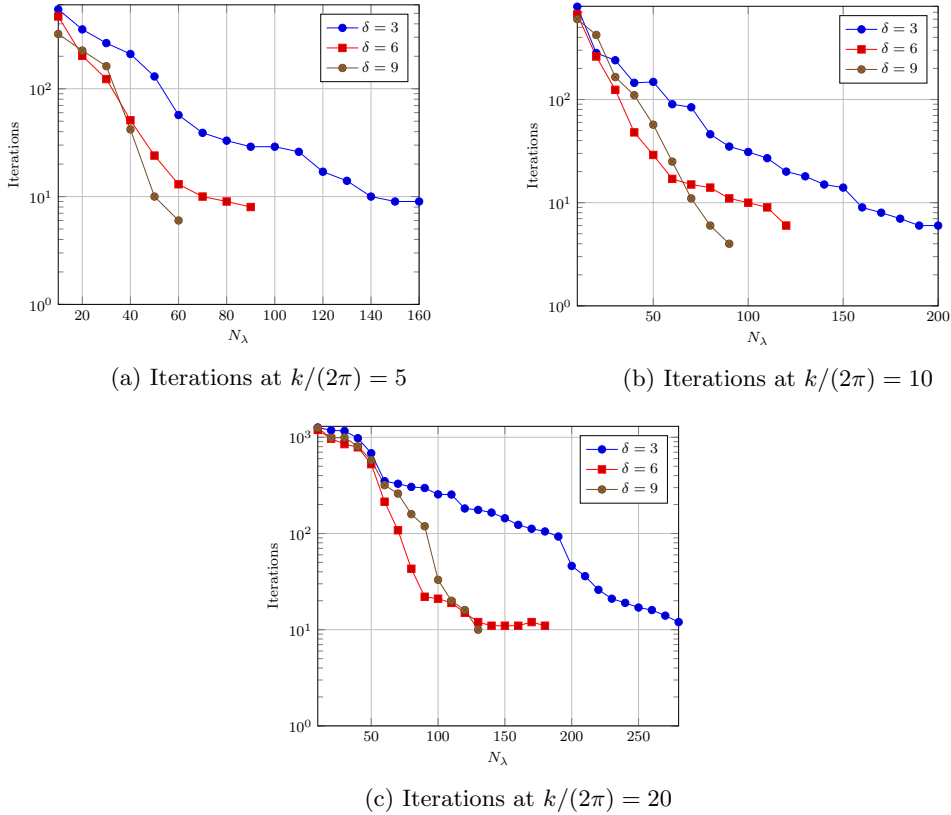

In Figure \ref{fig:errorH_inv_ol}, we present the number of iterations necessary for the HRC preconditioned GMRES to obtain $e_{rel}\approx \mathcal{O}(10^{-4})$ for different $N_\lambda$ at: (a) $k/(2\pi)=5$, (b) $k/(2\pi)=10$ and (c) $k/(2\pi)=20$. In each figure, we have three lines and each line represents the number of iterations necessary for convergence when using the HRC preconditioner obtained with different overlap parameter $\delta=3$, $6$ and $9$. In Table \ref{tab:error_exp_1_2}, we present the number of iterations with its respective error in the solution using GMRES without preconditioner  at wavenumbers $k/(2\pi)=5$, $10$ and $20$.

In Figure \ref{fig:error_inv_ol_ft} of Appendix \ref{s:appendB}, we present respectively $\|\tilde{\Ab}^{-1}_{RC}-\Ab^{-1}\|/\|\Ab^{-1}\|$ and $\|\tilde{\Ab}_{RC}-\Ab\|/\|\Ab\|$ for different $N_\lambda$ at $k/(2\pi)=5$ ((a),(b)), $10$ ((c),(d)) and $20$ ((e),(f)). Each line represents the error of the approximation for $\delta=3$, $6$ and $9$.

\begin{table}
\center
\caption{(Experiment I.1-2) Number of iterations and relative error of the solution using GMRES without preconditioner, $e_{rel}$, to solve Equation \eqref{eq:Hessian_prob_1} at frequencies $k=5(2\pi)$, $10(2\pi)$ and $20(2\pi)$.}
\begin{tabular}{|c||c|c|}\cline{2-3}\cline{2-3}
\multicolumn{1}{c|}{} & \multicolumn{2}{|c|}{GMRES}\\\hline
$k/(2\pi)$ & Iterations & Error\\\hline\hline
$5$  	 & 255 & 2.4e-4\\\hline
$10$ & 118 & 2.9e-4\\\hline
$20$ & 237 & 8.3e-4\\\hline
\end{tabular}
\label{tab:error_exp_1_2}
\end{table}

{\bf Summary:} The value of the parameter $N_\lambda$ needs to be larger at higher frequencies to obtain the same accuracy in the approximation of the inverse of the forward operator. This is expected due to the fact that the singular values of $\Ab$ decay faster at lower frequencies than at higher frequencies. We also notice that, in accordance with intuition, if we use a larger overlap parameter to create $\Ab^{-1}_{RAS}$ than we need smaller values for the parameter $N_\lambda$ to obtain a prescribed fixed accuracy. 

{\bf Experiment I.2 -- Influence of the number of subdomains on inverse preconditioning:} We analyze the performance of the HRC preconditioner when we change the number of subdomains when we have a constant overlap parameter. The following parameters are used: 
\begin{itemize}
\item Incoming waves: the incoming waves are given by the wave $u^{inc}(\xb=(x,y))=\exp(ik \xb\cdot\theta_j)$ with $k/(2\pi)=5$, $10$ and $20$, and $\theta_j$ prescribed as in the beginning of the subsection with $j=0,\ldots,7$;
\item Receivers: the receivers are located at $\xb_r=(x_r,y_r)=0.8(\cos(2\pi r/N_r), \sin(2\pi r/N_r))$, with $r=1,\ldots,N_r$, and $N_r=2,000$;
\item Scatterer: the scatterers points are distributed in an uniform grid of points with $N=64^2$ scatterers and their magnitude is given by the function $q_4$;
\item Domain decomposition: the number of subdomains is $N_s=4$, $16$, $25$, $36$ and $64$, and the overlap parameter is $\delta=8$;
\item HRC preconditioner: we use for the low-rank correction $N_\lambda=10+10m$, $m=0,\ldots,M_\lambda$, where $M_\lambda$ differs for each problem and it is the maximum number needed to obtain the prescribed error.
\end{itemize}

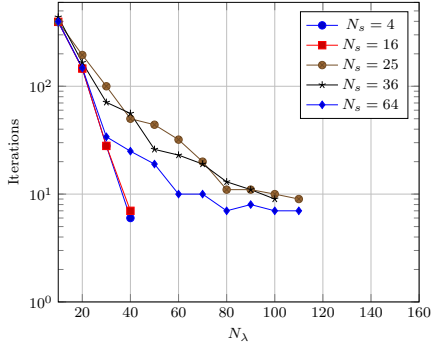
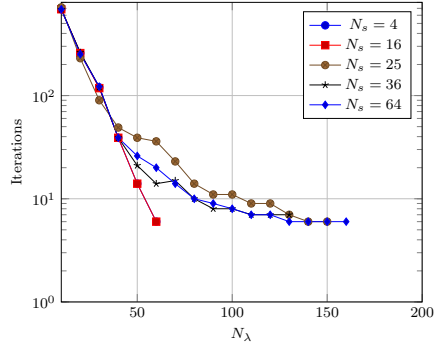
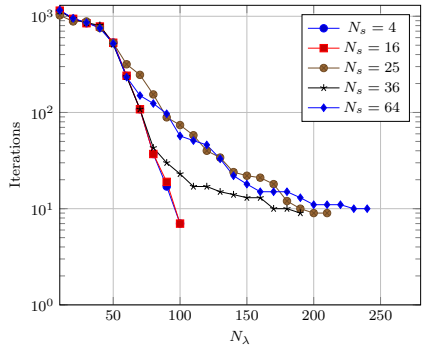
\begin{figure}
\footnotesize
\centering
\begin{subfigure}{.45\textwidth}
\begin{tikzpicture}[scale=0.7]
\begin{semilogyaxis}[xmin=10, xmax=160,domain=1:160,ymin=1, ymax=600,domain=1:600,
    xlabel={$N_\lambda$},
    ylabel={Iterations},
    grid=major,
    legend pos=north east,
    legend entries={$N_s=4$,$N_s=16$,$N_s=25$,$N_s=36$,$N_s=64$},
]
\addplot coordinates {
(10,393) (20,146) (30,28) (40,6)
};
\addplot coordinates {
(10,396) (20,146) (30,28) (40,7)
};
\addplot coordinates {
(10,396) (20,195) (30,100) (40,50) (50,44) (60,32)
(70,20) (80,11) (90,11) (100,10) (110,9) 
};
\addplot coordinates {
(10,438) (20,166) (30,71) (40,56) (50,26) (60,23)
(70,19) (80,13) (90,11) (100,9) 
};
\addplot coordinates {
(10,398) (20,148) (30,34) (40,25) (50,19) (60,10)
(70,10) (80,7) (90,8) (100,7) (110,7) 
};
\end{semilogyaxis}
\end{tikzpicture}
\caption{Iterations at $k/(2\pi)=5$}\label{fig:error_inv_sd_1}
\end{subfigure}
\begin{subfigure}{.45\textwidth}
\begin{tikzpicture}[scale=0.7]
\begin{semilogyaxis}[xmin=10, xmax=200,domain=1:200,ymin=1, ymax=800,domain=1:800,
    xlabel={$N_\lambda$},
    ylabel={Iterations},
    grid=major,
    legend pos=north east,
    legend entries={$N_s=4$,$N_s=16$,$N_s=25$,$N_s=36$,$N_s=64$},
]
\addplot coordinates {
(10,686) (20,258) (30,120) (40,39) (50,14) (60,6)
};
\addplot coordinates {
(10,686) (20,258) (30,118) (40,39) (50,14) (60,6)
};
\addplot coordinates {
(10,716) (20,230) (30,90) (40,49) (50,39) (60,36)
(70,23) (80,14) (90,11) (100,11) (110,9) (120,9)
(130,7) (140,6) (150,6)
};
\addplot coordinates {
(10,683) (20,261) (30,121) (40,39) (50,21) (60,14)
(70,15) (80,10) (90,8) (100,8) (110,7) (120,7)
(130,7)
};
\addplot coordinates {
(10,685) (20,252) (30,122) (40,39) (50,26) (60,20)
(70,14) (80,10) (90,9) (100,8) (110,7) (120,7)
(130,6) (140,6) (150,6) (160,6)
};
\end{semilogyaxis}
\end{tikzpicture}
\caption{Iterations at $k/(2\pi)=10$}\label{fig:error_inv_sd_2}
\end{subfigure}
\begin{subfigure}{.45\textwidth}
\begin{tikzpicture}[scale=0.7]
\begin{semilogyaxis}[xmin=10, xmax=280,domain=1:280,ymin=1, ymax=1300,domain=1:1300,
    xlabel={$N_\lambda$},
    ylabel={Iterations},
    grid=major,
    legend pos=north east,
    legend entries={$N_s=4$,$N_s=16$,$N_s=25$,$N_s=36$,$N_s=64$},
]
\addplot coordinates {
(10,1144) (20,941) (30,847) (40,783) (50,533) (60,240)
(70,108) (80,37) (90,17) (100,7)
};
\addplot coordinates {
(10,1144) (20,941) (30,847) (40,782) (50,531) (60,241)
(70,108) (80,37) (90,19) (100,7) 
};
\addplot coordinates {
(10,1022) (20,884) (30,883) (40,755) (50,528) (60,317)
(70,246) (80,154) (90,89) (100,74) (110,58) (120,40)
(130,34) (140,24) (150,22) (160,21) (170,18) (180,12)
(190,10) (200,9) (210,9) 
};
\addplot coordinates {
(10,1129) (20,947) (30,853) (40,807) (50,509) (60,235)
(70,108) (80,43) (90,30) (100,23) (110,17) (120,17)
(130,15) (140,14) (150,13) (160,13) (170,10) (180,10)
(190,9)
};
\addplot coordinates {
(10,1154) (20,952) (30,860) (40,752) (50,522) (60,233)
(70,150) (80,124) (90,97) (100,57) (110,51) (120,46)
(130,33) (140,22) (150,18) (160,15) (170,15) (180,15)
(190,13) (200,11) (210,11) (220,11) (230,10) (240,10)
};
\end{semilogyaxis}
\end{tikzpicture}
\caption{Iterations at $k/(2\pi)=20$}\label{fig:error_inv_sd_3}
\end{subfigure}
\caption{(Experiment I.2) Plot of the number of iterations necessary to obtain $e_{rel}\approx\mathcal{O}(10^{-4})$ using the HRC preconditioner for the GMRES for $N= 64^2$ with different $N_\lambda$ at wavenumbers: (a) $k/(2\pi)=5$, (b) $k/(2\pi)=10$ and (c) $k/(2\pi)=20$. Each line represents the number of iterations obtained using the HRC preconditioner using $N_s=4$, $16$, $25$, $36$ and $64$ subdomains. The overlap parameter is $\delta=8$. }\label{fig:errorH_inv_sd}
\end{figure}

In Figure \ref{fig:errorH_inv_sd}, we present the number of iterations necessary for the HRC preconditioned GMRES to obtain $e_{rel}\approx \mathcal{O}(10^{-4})$ for different $N_\lambda$ at different wavenumbers: (a) $k/(2\pi)=5$, (b) $k/(2\pi)=10$ and (c) $k/(2\pi)=20$. In each figure, we have five lines and each line represents the number of iterations necessary for convergence when using the HRC preconditioner obtained with different number of subdomains.

In Figure \ref{fig:error_inv_sd_ft} of Appendix \ref{s:appendB}, we present respectively $\|\tilde{\Ab}^{-1}_{RC}-\Ab^{-1}\|/\|\Ab^{-1}\|$ and $\|\tilde{\Ab}_{RC}-\Ab\|/\|\Ab\|$ for different $N_\lambda$ at $k/(2\pi)=5$ ((a),(b)), $10$ ((c),(d)) and $20$ ((e),(f)). Each line represents the error of the approximation for $N_s=4$, $16$, $25$, $36$ and $64$.

{\bf Summary:} We note that as we increase the number of subdomains used, we require a larger $N_\lambda$ parameter to obtain better results and eventually it deteriorates when we have a very large number of subdomains. 

{\bf Experiment I.3 -- Scalability of the inverse preconditioning:} This example aims to check the scalability of the HRC preconditioner. We fix the number of subdomains and increase the number of points in the grid and the overlap parameter. The following parameters are used: 
\begin{itemize}
\item Incoming waves: the incoming waves are given by the wave $u^{inc}(\xb=(x,y))=\exp(ik \xb\cdot\theta_j)$ with $k/(2\pi)=5$ and $20$, and $\theta_j$ prescribed as in the beginning of the subsection with $j=0,\ldots,7$;
\item Receivers: the receivers are located at $\xb_r=(x_r,y_r)=0.8(\cos(2\pi r/N_r), \sin(2\pi r/N_r))$, with $r=1,\ldots,N_r$, and  $N_r=10,000$;
\item Scatterer: we use a regular grid with $N=64^2$, $128^2$ and $256^2$ scatterers and their magnitude is given by the function $q_4$;
\item Domain decomposition: the number of subdomains is $N_s=16$ constant, we use the overlap parameter $\delta=3$ for the grid with $64^2$ scatterer points, $\delta=6$ for $128^2$ scatterer points and $\delta=12$ for $256^2$ scatterer points;
\item HRC preconditioner: we choose $N_\lambda=20$, $40$, $60$ and $80$ for wavenumber $k/(2\pi)=5$ and $N_\lambda=40$, $90$ and $120$ for wavenumber $k/(2\pi)=20$.
\end{itemize}

We present in Table \ref{tab:inv_scal_it} the number of GMRES iterations necessary for the method to converge to the relative error $e_{rel}$ with order of magnitude $\mathcal{O}(\Phi)$, with $\Phi=10^{-2}$, $10^{-3}$, $10^{-4}$ and $10^{-5}$. We present the respective $e_{rel}$ in Table \ref{tab:inv_scal_error} of Appendix \ref{s:appendB}. We report the values of $\|\tilde{\Ab}^{-1}_{RC}-\Ab^{-1}\|/\|\Ab^{-1}\|$ and $\|\tilde{\Ab}_{RC}-\Ab\|/\|\Ab\|$ for the experiments in Table \ref{tab:error_tf}.

\begin{table}
\caption{(Experiment I.3) Number of iterations necessary for the convergence of GMRES  with the prescribed  order of magnitude to the solution of the Equation \eqref{eq:Hessian_prob_1} using the HRC preconditioner at (a) $k/(2\pi)=5$ and (b) $k/(2\pi)=20$. The prescribed order of magnitude of the error is $\mathcal{O}(\Phi)$, with $\Phi=10^{-2}$, $10^{-3}$, $10^{-4}$ and $10^{-5}$. The total number of scatterers in the domain are $N=32^2$, $64^2$ and $128^2$. The number of subdomains is $N_s=16$ and the overlap parameter $\delta=3$, $6$ and $12$ respectively for the domains with $N=32^2$, $64^2$ and $128^2$ scatterers. The low-rank approximation parameter is $N_\lambda=20$, $40$, $60$ and $80$ at $k=5(2\pi)$ and $N_\lambda=40$, $90$ and $120$ at $k=20(2\pi)$.}
\label{tab:inv_scal_it}
\begin{subtable}{\textwidth}
\center
\caption{Iterations for $k/(2\pi)=5$}
\begin{tabular}{|c||c||*{6}{c|}}\cline{4-7}\cline{4-7}
\multicolumn{3}{c|}{} & \multicolumn{4}{c|}{$N_\lambda$}\\\hline
$N$ & $\mathcal{O}(Error)$ & GMRES & 20 & 40 & 60 & 80 \\\hline\hline
\multirow{4}{*}{$32^2$}   & $10^{-2}$ & 221 & 205 & 67   & 14 & 3\\
                                        & $10^{-3}$ & 248 & 231 & 81   & 15 & 4\\
                                        & $10^{-4}$ & 263 & 261 & 91   & 17 & 4\\
                                        & $10^{-5}$ & 274 & 286 & 102 & 20 & 5\\\hline
\multirow{4}{*}{$64^2$}   & $10^{-2}$ & 198 & 159 & 37 & 11  & 8 \\
                                        & $10^{-3}$ & 238 & 183 & 45 & 12 & 9\\
                                        & $10^{-4}$ & 255 & 198 & 51 & 13 & 9\\
                                        & $10^{-5}$ & 266 & 217 & 58 & 15 & 11\\\hline
\multirow{4}{*}{$128^2$} & $10^{-2}$ & 195 & 110 & 21 &   9 & 6\\
                                        & $10^{-3}$ & 237 & 134 & 25 & 11 & 7\\
                                        & $10^{-4}$ & 256 & 152 & 26 & 13 & 8\\
                                        & $10^{-5}$ & 268 & 169 & 29 & 14 & 9\\\hline
\end{tabular}
\end{subtable}
\begin{subtable}{\textwidth}
\center
\caption{Iterations for $k/(2\pi)=20$}
\begin{tabular}{|c||c||*{5}{c|}}\cline{4-6}\cline{4-6}
\multicolumn{3}{c|}{} & \multicolumn{3}{c|}{$N_\lambda$}\\\hline
$N$ & $\mathcal{O}(Error)$ & GMRES & 40 & 90 & 120  \\\hline\hline
\multirow{4}{*}{$32^2$}   & $10^{-2}$ &  112  &  112  & 33 & 2\\
                                        & $10^{-3}$ &  158 &  131  & 40 & 3\\
                                        & $10^{-4}$ &  208 &  150  & 44 & 4\\
                                        & $10^{-5}$ &  250 &  166  & 46 & 4\\\hline
\multirow{4}{*}{$64^2$}   & $10^{-2}$ &  122 &   657 & 15 & 11\\
                                        & $10^{-3}$ &  182 &   754 & 19 & 14\\
                                        & $10^{-4}$ &  236 &   786 & 22 & 16\\
                                        & $10^{-5}$ &  285 &   805 & 24 & 17\\\hline
\multirow{4}{*}{$128^2$} & $10^{-2}$ &  123 &   544 & 10 & 9\\
                                        & $10^{-3}$ &  191 &   630 & 12 & 12\\
                                        & $10^{-4}$ &  266 &   694 & 14 & 13\\
                                        & $10^{-5}$ &  330 &   732 & 16 & 13\\\hline
\end{tabular}
\end{subtable}
\end{table}

{\bf Summary:} The results show that the method is fully scalable with the increase of the number of points in the domain, requiring approximately the same number of iterations to obtain the same accuracy for increasing domain size.

{\bf Experiment I.4 -- Comparison with low-rank inverse preconditioner:} we compare the HRC preconditioner with a low-rank preconditioner (called here LR preconditioner) obtained by inverting the regularized low-rank approximation of the operator $\Hb$. A similar version of this preconditioner was previously presented in \cite{Hohage2001}. 

First, we describe how to construct the LR preconditioner $\tilde{\Hb}^{-1}_{LR}$. We compute the singular value decomposition of $\Jb^{\ast}\Jb=\Ub\Sb\Ub^\ast$. In practice, the SVD is constructed using randomized projections. Here, we use the exact SVD and truncate it to the target rank. Next, approximate $\Jb^{\ast}\Jb\approx \Ub_{N_\lambda}\Sb_{N_\lambda}\Ub_{N_\lambda}^\ast$, where $\Ub_{N_\lambda}=\Ub(:,$$1$:$N_\lambda)$, $\Sb_{N_\lambda}=\Sb($$1$:$N_\lambda,$$1$:$N_\lambda)$ and $\Vb_{N_\lambda}=\Vb(:,$$1$:$N_\lambda)$. Finally, set the LR preconditioner as
\begin{equation}
\tilde{\Hb}^{-1}_{LR}=\Ub_{N_\lambda}(\Sb_{N_\lambda}+\beta \Ib)^{-1}\Ub_{N_\lambda}^{\ast} + \beta^{-1}(\Ib-\Ub_{N_\lambda}\Ub_{N_\lambda}^{\ast}).\label{eq:lr_prec}
\end{equation}

The following parameters are used in this experiment: 
\begin{itemize}
\item Incoming waves: the incoming waves are given by $u^{inc}(\xb=(x,y))=\exp(ik \xb\cdot\theta_j)$ with three frequencies, $k/(2\pi)=5$ $10$ and $20$, and $\theta_j$ prescribed as in the beginning of the subsection with $j=0,\ldots,7$;
\item Receivers: the receivers are located at $\xb_r=(x_r,y_r)=0.8(\cos(2\pi r/N_r), \sin(2\pi r/N_r))$, with $r=1,\ldots,N_r$, and $N_r=2000$;
\item Scatterer: we use a regular grid with $N=64^2$ scatterers and their magnitude is given by the function $q_4$;
\item Domain decomposition: the number of subdomains is $N_s=16$ and the overlap parameter is $\delta=4$ and $8$;
\item HRC and LR preconditioners: we use for the low-rank correction $N_\lambda=10+10m$, $m=0,\ldots,M_\lambda$, where $M_\lambda$ differs for each problem and it is the maximum number needed to obtain the relative error $e_{rel}\approx \mathcal{O}(10^{-4})$.
\end{itemize}

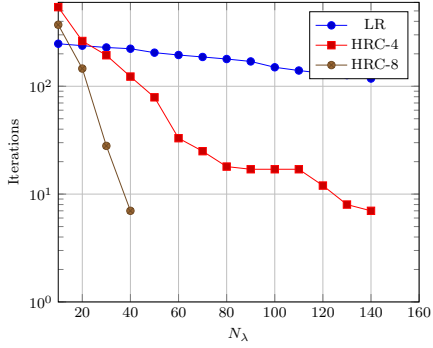
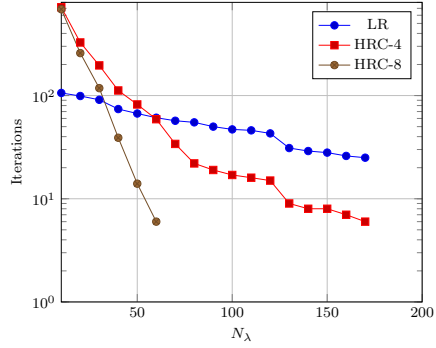
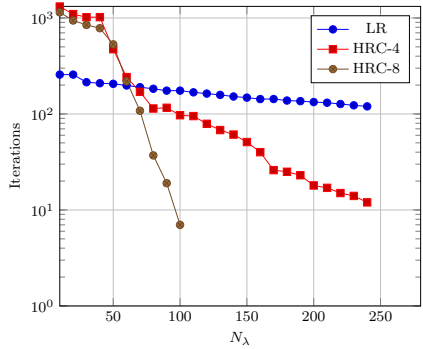
\begin{figure}
\footnotesize
\centering
\begin{subfigure}{.45\textwidth}
\begin{tikzpicture}[scale=0.7]
\begin{semilogyaxis}[xmin=10, xmax=160,domain=1:160,ymin=1, ymax=600,domain=1:600,
    xlabel={$N_\lambda$},
    ylabel={Iterations},
    grid=major,
    legend pos=north east,
    legend entries={LR,HRC-4,HRC-8},
]
\addplot coordinates {
(10,248) (20,238) (30,229) (40,223) (50,205) (60,195)
(70,187) (80,179) (90,170) (100,150) (110,140) (120,132) 
(130,127) (140,118)
};
\addplot coordinates {
(10,542) (20,263) (30,194) (40,123) (50,79) (60,33)
(70,25) (80,18) (90,17) (100,17) (110,17) (120,12) 
(130,8) (140,7)
};
\addplot coordinates {
(10,372) (20,146) (30,28) (40,7)
};
\end{semilogyaxis}
\end{tikzpicture}
\caption{Iterations at $k/(2\pi)=5$}\label{fig:Hcomparison_I4:a}
\end{subfigure}
\begin{subfigure}{.45\textwidth}
\begin{tikzpicture}[scale=0.7]
\begin{semilogyaxis}[xmin=10, xmax=200,domain=1:200,ymin=1, ymax=800,domain=1:800,
    xlabel={$N_\lambda$},
    ylabel={Iterations},
    grid=major,
    legend pos=north east,
    legend entries={LR,HRC-4,HRC-8},
]
\addplot coordinates {
(10,106) (20,99) (30,91) (40,74) (50,67) (60,61)
(70,57) (80,55) (90,50) (100,47) (110,46) (120,43) 
(130,31) (140,29) (150,28) (160,26) (170,25)
};
\addplot coordinates {
(10,717) (20,327) (30,196) (40,112) (50,82) (60,59)
(70,34) (80,22) (90,19) (100,17) (110,16) (120,15) 
(130,9) (140,8) (150,8) (160,7) (170,6)
};
\addplot coordinates {
(10,686) (20,258) (30,118) (40,39) (50,14) (60,6)
};
\end{semilogyaxis}
\end{tikzpicture}
\caption{Iterations at $k/(2\pi)=10$}\label{fig:Hcomparison_I4:b}
\end{subfigure}
\begin{subfigure}{.45\textwidth}
\begin{tikzpicture}[scale=0.7]
\begin{semilogyaxis}[xmin=10, xmax=280,domain=1:280,ymin=1, ymax=1320,domain=1:1320,
    xlabel={$N_\lambda$},
    ylabel={Iterations},
    grid=major,
    legend pos=north east,
    legend entries={LR,HRC-4,HRC-8},
]
\addplot coordinates {
(10,257) (20,257) (30,214) (40,209) (50,206) (60,198)
(70,190) (80,183) (90,175) (100,175) (110,168) (120,163) 
(130,158) (140,152) (150,148) (160,143) (170,143) (180,138)
(190,136) (200,133) (210,131) (220,127) (230,123) (240,120)
};
\addplot coordinates {
(10,1315) (20,1097) (30,1018) (40,1020) (50,473) (60,242)
(70,170) (80,114) (90,116) (100,97) (110,95) (120,79) 
(130,68) (140,61) (150,51) (160,40) (170,26) (180,25)
(190,23) (200,18) (210,17) (220,15) (230,14) (240,12)
};
\addplot coordinates {
(10,1144) (20,941) (30,847) (40,782) (50,531) (60,219)
(70,108) (80,37) (90,19) (100,7) 
};
\end{semilogyaxis}
\end{tikzpicture}
\caption{Iterations at $k/(2\pi)=20$}\label{fig:Hcomparison_I4:c}
\end{subfigure}
\caption{(Experiment I.4) Plot of the number of iterations necessary to obtain $e_{rel}\approx\mathcal{O}(10^{-4})$ using LR preconditioned and HRC preconditioned GMRES using $N= 64^2$ with different $N_\lambda$. That is, we use a rank of $N_\lambda$ for both HRC and LR preconditioners. We perform the tests at wavenumbers: (a) $k/(2\pi)=5$, (b) $k/(2\pi)=10$ and (c) $k/(2\pi)=20$. In each picture, the blue line represents the number of iterations using the LR preconditioner, the red and brown lines represents the number of iterations using the HRC preconditioner with $\delta=4$ and $8$, respectively (HRC-4 for $\delta=4$, and HRC-8 for $\delta=8$).  The number of subdomains used for the HRC preconditioner is fixed, $N_s=16$.}\label{fig:Hcomparison_I4}
\end{figure}

In Figure \ref{fig:Hcomparison_I4}, we present the number of iterations necessary for the preconditioned GMRES to obtain relative error $e_{rel}\approx \mathcal{O}(10^{-4})$ for different $N_\lambda$ at different wavenumbers: (a) $k/(2\pi)=5$, (b) $k/(2\pi)=10$ and (c)$k/(2\pi)=20$. The blue line represents the number of iterations necessary using the LR preconditioner, the red line represents the number of iterations necessary using the HRC preconditioner with $\delta=4$ (in the legend of the figure as HRC-4) and the brown line represents the number of iterations necessary using the HRC preconditioner with $\delta=8$ (in the legend of the figure as HRC-8).

In Figure \ref{fig:Heig_example_I4_k5}, we report the singular values for $\Hb$, $\tilde{\Hb}_{LR}$ with $N_\lambda=140$, $\tilde{\Hb}_{HRC}$ with $N_\lambda=140$, $N_s=16$ and $\delta=4$, and $\tilde{\Hb}_{HRC}$ with $N_\lambda=40$, $N_s=16$ and $\delta=8$ when the incoming incident waves have wavenumber $k/(2\pi)=5$. In Figure \ref{fig:Heig_example_I4_k20}, we report the singular values for $\Hb$, $\tilde{\Hb}_{LR}$ with $N_\lambda=240$, $\tilde{\Hb}_{HRC}$ with $N_\lambda=240$, $N_s=16$ and $\delta=4$, and $\tilde{\Hb}_{HRC}$ with $N_\lambda=100$, $N_s=16$ and $\delta=8$ when the incoming incident waves have wavenumber $k/(2\pi)=20$

{\bf Summary:} When using small $N_\lambda$ the number of iterations for convergence of GMRES is larger when using the HRC preconditioner than when using the LR preconditioner. This behavior is reversed when using larger $N_\lambda$ and the HRC preconditioner performance is much better than of the LR preconditioner.

\begin{figure}[h!]
  \centering
  \begin{subfigure}{.45\linewidth}
\includegraphics[width=\textwidth]{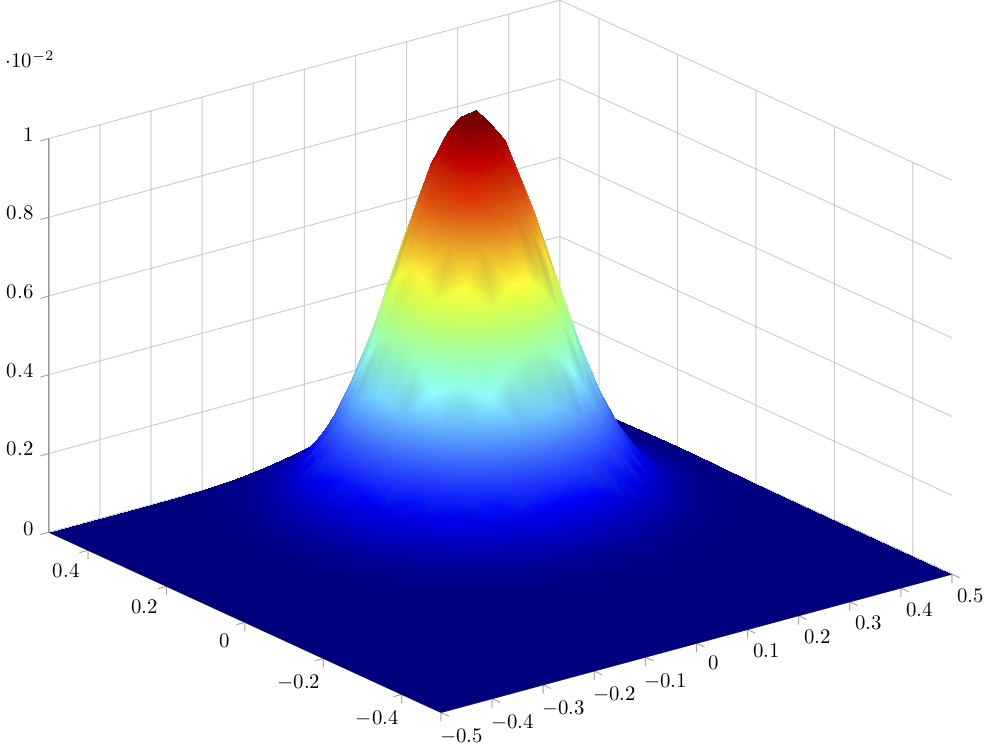}
\caption{isometric view of $q_b$}
\end{subfigure}
\begin{subfigure}{.45\linewidth}
\includegraphics[width=.8\textwidth]{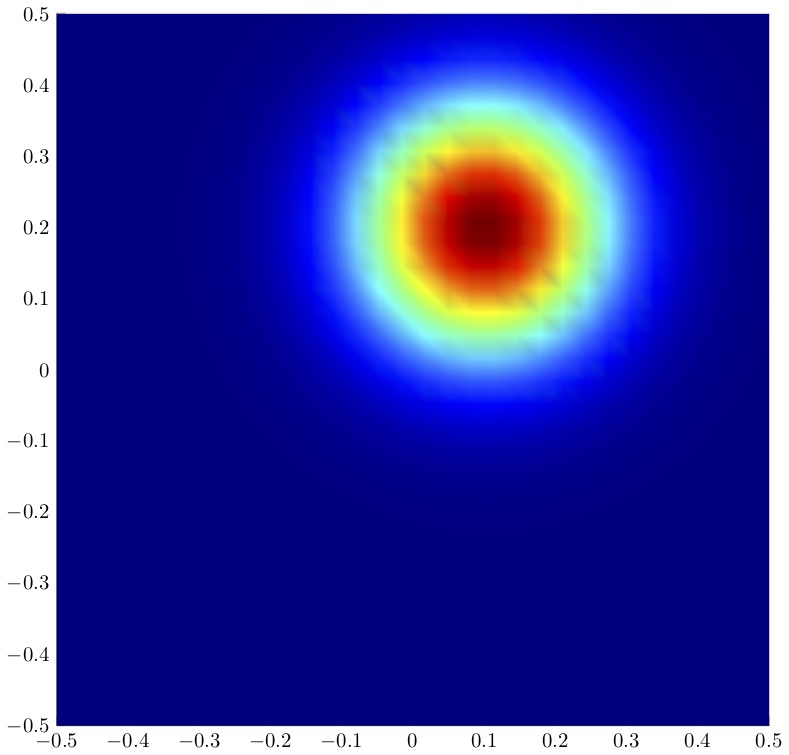}
\caption{top view of domain of $q_b$}
\end{subfigure}
\caption{The (a) isometric view and (b) top view of the domain $q_b$ used in the Experiment I.4.}
\label{fig:domain_example_I4}
\end{figure}

{\bf Experiment I.5 -- Nonlinear inverse problem:} In  this example, we compare the full reconstruction of a scatterer using GMRES with no preconditioner, GMRES with the LR preconditioner, and GMRES with the HRC preconditioner. We use the recursive linearization algorithm (RLA) to reconstruct the scatterer, given scattered data generated by incoming waves with multiple frequencies. The scatterer magnitudes are given by
\begin{equation*}
q_b(\xb=(x,y))=0.01\exp\left(-\left((x-0.1)^2+(y-0.2)^2\right)/0.03\right),
\end{equation*}
and can be seen in Figure \ref{fig:domain_example_I4}. We have chosen this function because it is very easy to simulate and obtain a very accurate reconstruction of it using a relatively low frequency amount of data, in comparison to $q_4$. 

The following parameters are used to simulate this experiment: 
\begin{itemize}
\item Incoming waves: the incoming waves are given by $u^{inc}(\xb=(x,y))=\exp(ik_\ell \xb \cdot \theta_j)$ with $k_\ell=1+0.25\ell$, for $\ell=1,\ldots,37$, and $\theta_j$ prescribed as in the beginning of the subsection with $j=0,\ldots,7$;
\item Receivers: the receivers are located at $\xb_r=(x_r,y_r)=0.8(\cos(2\pi r/N_r),\sin(2\pi r/N_r))$, with $r=1,\ldots,N_r$, and $N_r=2000$;
\item Scatterer: we use a regular grid with $N=32^2$ scatterers and their magnitude is given by the function $q_b$;
\item Domain decomposition: the number of subdomains is fixed set to $N_s=16$ and the overlap parameter is $\delta=4$;
\item HRC and LR preconditioners: to keep the number of iterations low, the choice of the parameter $N_\lambda$ must depend on the wavenumber $k$. We decide to use the function $N_\lambda(k)= \lceil 40k/9+140/9\rceil$. With this function, $N_\lambda(1)=20$ and $N_\lambda(10)=60$;
\item GMRES: the tolerance of the residual is $10^{-7}$ and the maximum number of iterations is $1,000$, with no restarts being used;
\item Gauss-Newton: the stopping criteria parameters for the Gauss-Newton method are the maximum number of iterations equal to $50$, the norm of the update $\delta q$ must be less than $10^{-3}/k$ and the norm of the objective functional must be less $10^{-4}/k$;
\item Regularization: we choose $\beta=10^{-0.9k-3.7}$, so that at $\beta(1)\approx2.5 \times10^{-5}$ and $\beta(10)=2\times 10^{-13}$; and
\item Initial guess: the initial guess is the regular grid scatterer points with the magnitude given by the function identically zero in the domain.
\end{itemize}

\begin{figure}
\center
\begin{subfigure}{.3\textwidth}
\includegraphics[width=.9\textwidth]{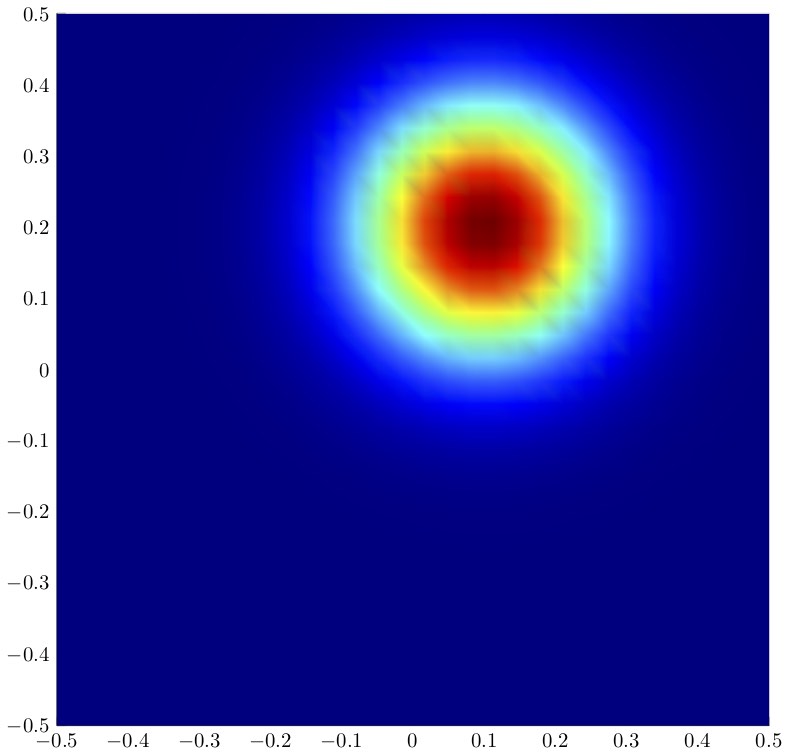}
\caption{GMRES }\label{fig:inv_full_sol_1}
\end{subfigure}
\begin{subfigure}{.3\textwidth}
\includegraphics[width=.9\textwidth]{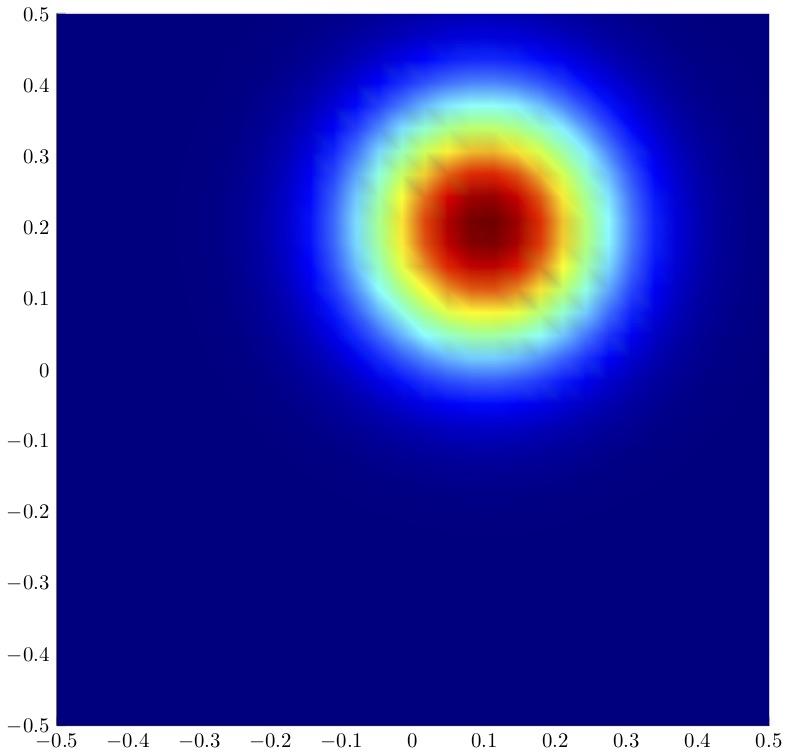}
\caption{LR Precconditioner}\label{fig:inv_full_sol_2}
\end{subfigure}
\begin{subfigure}{.3\textwidth}
\includegraphics[width=.9\textwidth]{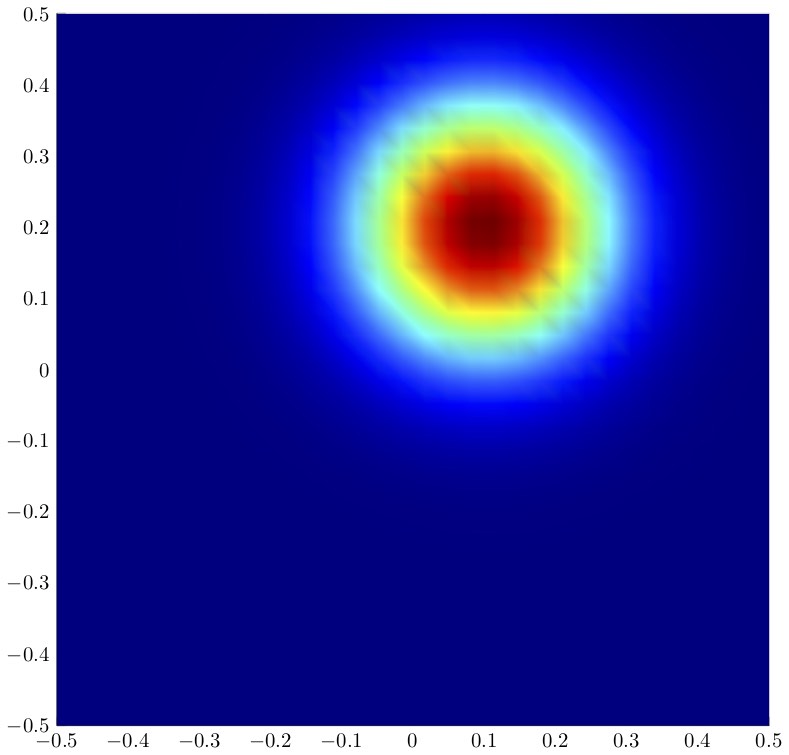}
\caption{HRC Preconditioner}\label{fig:inv_full_sol_3}
\end{subfigure}
\caption{(Experiment I.5) Reconstructions of the scatterer $q_b$ obtained by the RLA using: (a) GMRES without preconditioner, (b) LR preconditioned GMRES and (c) HRC preconditioned GMRES. }\label{fig:inv_full_sol}
\end{figure}

\begin{figure}
\center
\begin{subfigure}{.45\textwidth}
\input{err_exI4.tex}
\caption{Error of the solution}\label{fig:inv_full_bench_1}
\end{subfigure}
\begin{subfigure}{.45\textwidth}
\input{iter_GN_exI4.tex}
\caption{Gauss-Newton iterations}\label{fig:inv_full_bench_2}
\end{subfigure}

\begin{subfigure}{.45\textwidth}
\input{iter_exI4.tex}
\caption{GMRES iterations}\label{fig:inv_full_bench_3}
\end{subfigure}
\caption{(Experiment I.5) Regarding the full reconstruction of the scatterer $q_b$ using RLA, we present at each wavenumber $k$: (a) the relative error of the reconstruction with respect to $q_b$, (b) the number of iterations necessary for the convergence of the Gauss-Newton method, and (c) the total number of GMRES iterations used. In each figure, the curve with $\Box$ marks has the values for the solution using GMRES with no preconditioner, the curve with the $\circ$ marks has the values for the solution using LR preconditioned GMRES and the curve with the $\times$ marks has the values for the solution using the HRC preconditioned GMRES.}\label{fig:inv_full_bench}
\end{figure}
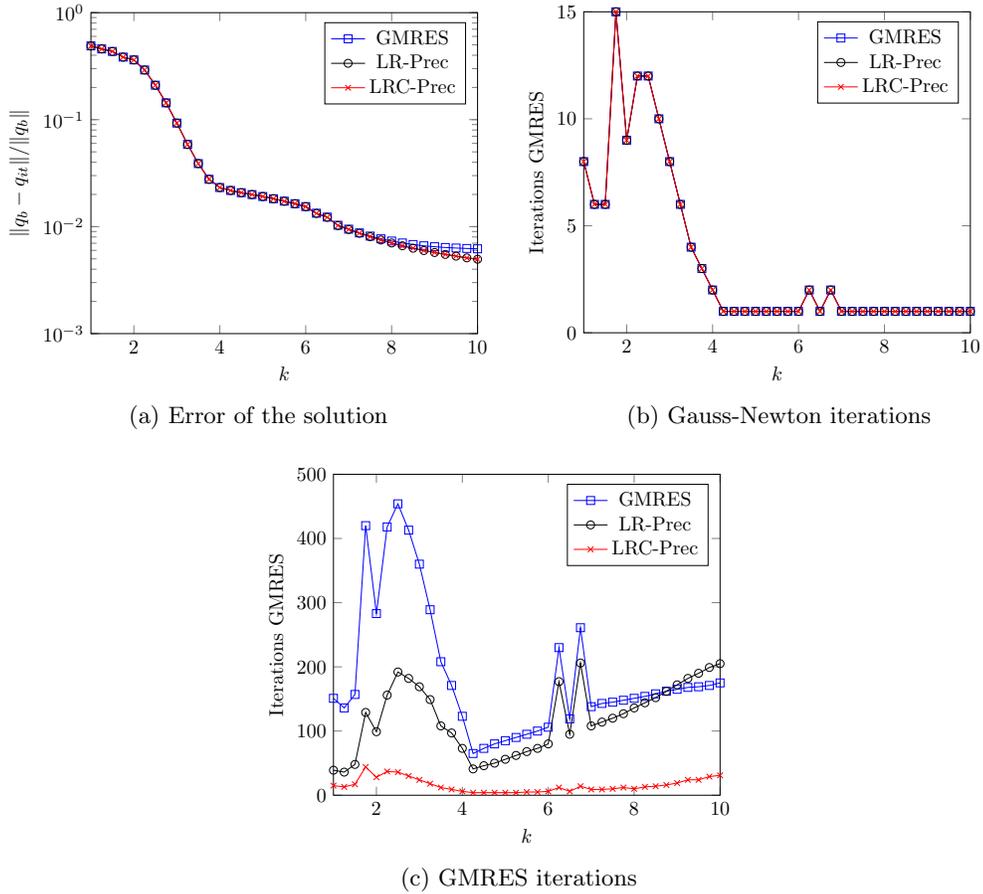

The reconstructions obtained using GMRES without preconditioner, the LR preconditioned GMRES and the HRC preconditioned GMRES can be seen in Figures \ref{fig:inv_full_sol_1}, \ref{fig:inv_full_sol_2} and \ref{fig:inv_full_sol_3} respectively. In Figure \ref{fig:inv_full_bench}, we present at each wavenumber: (a) the relative error between the reconstruction and $q_b$, (b) the number of iterations of the Gauss-Newton method and (c) the total number of GMRES iterations used.

In Table \ref{tab:inv_full}, we present at the wavenumbers $k=1$, $2.5$, $5$, $7.5$ and $10$ the number of GMRES iterations with and without the preconditioners in the columns labeled ``Step". The total number of iterations from wavenumber $1$ up to the wavenumber $k$ is presented in the column labeled ``Total".

\begin{table}
\center
\caption{(Experiment I.5) Total number of iterations of GMRES at each frequency with and without preconditioner. In the column labeled ``Step", we present the sum of the number of iterations of GMRES for the solution of the Gauss-Newton method at the respective wavenumber in the column $k$. In the column labeled ``Total",  we present the number of iterations of GMRES used for the RLA from the wavenumber 1 up to the respective  wavenumber in the column $k$.}
\label{tab:inv_full}
\begin{tabular}{|c||c|c|c|c|c|c|c|}\cline{2-3}\cline{5-8}
 \multicolumn{1}{c}{}       & \multicolumn{2}{|c|}{GMRES} & \multicolumn{1}{c}{} & \multicolumn{2}{|c|}{LR Preconditioner} & \multicolumn{2}{c|}{HRC Preconditioner}\\\hline
$k$  & Step & Total & $N_\lambda$ & Step & Total & Step & Total \\\hline\hline
  1      & 151 &   151 & 20 & 39   & 39     & 15 &   15  \\
  2.5   & 454 & 2019 & 27 & 192 & 699   & 36 & 190  \\
  5      &   85 & 3886 & 38 & 56   & 1670 & 4   & 305  \\
  7.5   & 145 & 5313 & 49 & 120 & 2773 & 10 & 385  \\
10      & 175 & 6934 & 60 & 205 & 4442 & 31 & 577  \\\hline
\end{tabular}
\end{table}

{\bf Summary:} The total number of GMRES iterations when using the HRC preconditioner is ten times smaller than the total number of GMRES iterations without the preconditioner. The preconditioner is very effective to be used with the RLA and even though we can experience a small increase in the number  of iterations with the increase of the wavenumber, this can be remedy by using a more aggressive choice of $N_\lambda$.

\subsection{Conclusion on the preconditioning for the inverse problem}
With the right choice of the number of singular values, we can construct the RC preconditioner to approximates the inverse of the forward operator and consequently use this approximation to construct the HRC preconditioner to decrease the number of iterations for the solution of the system \eqref{eq:Hessian_prob_1}. For larger values of $N$, additional approximation for HRC is needed, for example using H-matrix methods, to approximate the entries of $\tilde{\Hb}$ and then factorize it.

%% file: err_exI4.tex
%
%
\begin{tikzpicture}[scale=0.75]

\begin{axis}[%
xmin=1,
xmax=10,
ymode=log,
ymin=0.001,
ymax=1,
yminorticks=true,
    xlabel={$k$},
    ylabel={$\|q_b-q_{it}\|/\|q_b\|$},
    legend pos={north east},
    legend entries={GMRES,LR-Prec,LRC-Prec}
]
\addplot [color=blue,solid,mark=square,mark options={solid}]
  table[row sep=crcr]{%
1	0.488866571750757\\
1.25	0.458348426927888\\
1.5	0.434043818687475\\
1.75	0.384336012299078\\
2	0.361807094272891\\
2.25	0.291181379799955\\
2.5	0.209786039845685\\
2.75	0.14350546059028\\
3	0.0929324293108572\\
3.25	0.0586813078197453\\
3.5	0.0388601656898878\\
3.75	0.0277028233162365\\
4	0.0231997219787256\\
4.25	0.021778806354784\\
4.5	0.0207613284554842\\
4.75	0.0199403330020418\\
5	0.0191237256926422\\
5.25	0.0182377060915822\\
5.5	0.0173013110736338\\
5.75	0.0163854565101861\\
6	0.0153977568403244\\
6.25	0.0133453698591004\\
6.5	0.0122776273423279\\
6.75	0.0103032718520063\\
7	0.00946820736282391\\
7.25	0.00873569601105775\\
7.5	0.00815762280863731\\
7.75	0.00769607807367152\\
8	0.0073318303327827\\
8.25	0.00703640467799568\\
8.5	0.0068081604873903\\
8.75	0.00661893334217832\\
9	0.00648879533980654\\
9.25	0.00638221160699841\\
9.5	0.00630934618354211\\
9.75	0.0062510545228166\\
10	0.006205509875242\\
};
\addplot [color=black,solid,mark=o,mark options={solid}]
  table[row sep=crcr]{%
1	0.488791051716843\\
1.25	0.458322628854076\\
1.5	0.434028748098967\\
1.75	0.384327002678395\\
2	0.361800869489528\\
2.25	0.291180369258296\\
2.5	0.209788114298219\\
2.75	0.143508551035161\\
3	0.0929348073882157\\
3.25	0.0586833126165759\\
3.5	0.038861203491186\\
3.75	0.0277022225094813\\
4	0.0231973613441648\\
4.25	0.0217751278035193\\
4.5	0.0207575711014082\\
4.75	0.0199367666361409\\
5	0.019120085702449\\
5.25	0.018233374970924\\
5.5	0.0172934722825461\\
5.75	0.0163683575849941\\
6	0.0153653175477292\\
6.25	0.0133012090308472\\
6.5	0.0122226625751319\\
6.75	0.0102372781098486\\
7	0.0094129190078927\\
7.25	0.00867093096655685\\
7.5	0.00805077337239861\\
7.75	0.00750299983120757\\
8	0.00702991348288894\\
8.25	0.0066241054452338\\
8.5	0.00628521976931505\\
8.75	0.00597439765572433\\
9	0.00572843137570342\\
9.25	0.00549761480013898\\
9.5	0.00530017833726263\\
9.75	0.00510786140937514\\
10	0.00495579550791586\\
};
\addplot [color=red,solid,mark=x,mark options={solid}]
  table[row sep=crcr]{%
1	0.488791053869173\\
1.25	0.458322632192429\\
1.5	0.434028751858764\\
1.75	0.38432701100686\\
2	0.361800877471549\\
2.25	0.291180373728929\\
2.5	0.20978811707489\\
2.75	0.143508552760583\\
3	0.0929348086479022\\
3.25	0.0586833133961306\\
3.5	0.0388612039872558\\
3.75	0.0277022227860154\\
4	0.023197361568835\\
4.25	0.0217751279029258\\
4.5	0.0207575708863877\\
4.75	0.0199367650855672\\
5	0.0191200820234609\\
5.25	0.018233369209602\\
5.5	0.0172934606170349\\
5.75	0.0163683390142519\\
6	0.0153652858399172\\
6.25	0.0133011616110686\\
6.5	0.0122226080361298\\
6.75	0.0102372080478179\\
7	0.00941283114805669\\
7.25	0.0086708068910869\\
7.5	0.00805054624399711\\
7.75	0.00750257118574151\\
8	0.0070293817662354\\
8.25	0.00662344424350433\\
8.5	0.00628439732194091\\
8.75	0.00597337712408504\\
9	0.00572711976047705\\
9.25	0.00549606681349777\\
9.5	0.00529805125431803\\
9.75	0.00510483424516721\\
10	0.0049510879849076\\
};
\end{axis}
\end{tikzpicture}%

%% file: iter_GN_exI4.tex
%
%
\begin{tikzpicture}[scale=0.75]

\begin{axis}[%
xmin=1,
xmax=10,
ymin=0,
ymax=15,
    xlabel={$k$},
    ylabel={Iterations GMRES},
    legend pos=north east,
    legend entries={GMRES,LR-Prec,LRC-Prec}
]
\addplot [color=blue,solid,mark=square,mark options={solid}]
  table[row sep=crcr]{%
1	8\\
1.25	6\\
1.5	6\\
1.75	15\\
2	9\\
2.25	12\\
2.5	12\\
2.75	10\\
3	8\\
3.25	6\\
3.5	4\\
3.75	3\\
4	2\\
4.25	1\\
4.5	1\\
4.75	1\\
5	1\\
5.25	1\\
5.5	1\\
5.75	1\\
6	1\\
6.25	2\\
6.5	1\\
6.75	2\\
7	1\\
7.25	1\\
7.5	1\\
7.75	1\\
8	1\\
8.25	1\\
8.5	1\\
8.75	1\\
9	1\\
9.25	1\\
9.5	1\\
9.75	1\\
10	1\\
};
\addplot [color=black,solid,mark=o,mark options={solid}]
  table[row sep=crcr]{%
1	8\\
1.25	6\\
1.5	6\\
1.75	15\\
2	9\\
2.25	12\\
2.5	12\\
2.75	10\\
3	8\\
3.25	6\\
3.5	4\\
3.75	3\\
4	2\\
4.25	1\\
4.5	1\\
4.75	1\\
5	1\\
5.25	1\\
5.5	1\\
5.75	1\\
6	1\\
6.25	2\\
6.5	1\\
6.75	2\\
7	1\\
7.25	1\\
7.5	1\\
7.75	1\\
8	1\\
8.25	1\\
8.5	1\\
8.75	1\\
9	1\\
9.25	1\\
9.5	1\\
9.75	1\\
10	1\\
};
\addplot [color=red,solid,mark=x,mark options={solid}]
  table[row sep=crcr]{%
1	8\\
1.25	6\\
1.5	6\\
1.75	15\\
2	9\\
2.25	12\\
2.5	12\\
2.75	10\\
3	8\\
3.25	6\\
3.5	4\\
3.75	3\\
4	2\\
4.25	1\\
4.5	1\\
4.75	1\\
5	1\\
5.25	1\\
5.5	1\\
5.75	1\\
6	1\\
6.25	2\\
6.5	1\\
6.75	2\\
7	1\\
7.25	1\\
7.5	1\\
7.75	1\\
8	1\\
8.25	1\\
8.5	1\\
8.75	1\\
9	1\\
9.25	1\\
9.5	1\\
9.75	1\\
10	1\\
};
\end{axis}
\end{tikzpicture}%

%% file: iter_exI4.tex
%
%
\begin{tikzpicture}[scale=0.75]

\begin{axis}[%
xmin=1,
xmax=10,
ymin=0,
ymax=500,
    xlabel={$k$},
    ylabel={Iterations GMRES},
    legend pos=north east,
    legend entries={GMRES,LR-Prec,LRC-Prec}
]
\addplot [color=blue,solid,mark=square,mark options={solid}]
  table[row sep=crcr]{%
1	151\\
1.25	136\\
1.5	157\\
1.75	420\\
2	283\\
2.25	418\\
2.5	454\\
2.75	413\\
3	360\\
3.25	289\\
3.5	208\\
3.75	171\\
4	123\\
4.25	65\\
4.5	73\\
4.75	80\\
5	85\\
5.25	90\\
5.5	95\\
5.75	100\\
6	106\\
6.25	230\\
6.5	119\\
6.75	261\\
7	138\\
7.25	143\\
7.5	145\\
7.75	148\\
8	151\\
8.25	154\\
8.5	158\\
8.75	162\\
9	165\\
9.25	168\\
9.5	169\\
9.75	171\\
10	175\\
};
\addplot [color=black,solid,mark=o,mark options={solid}]
  table[row sep=crcr]{%
1	39\\
1.25	36\\
1.5	48\\
1.75	129\\
2	99\\
2.25	156\\
2.5	192\\
2.75	182\\
3	169\\
3.25	149\\
3.5	108\\
3.75	97\\
4	73\\
4.25	41\\
4.5	46\\
4.75	50\\
5	56\\
5.25	62\\
5.5	68\\
5.75	73\\
6	80\\
6.25	177\\
6.5	95\\
6.75	206\\
7	108\\
7.25	114\\
7.5	120\\
7.75	127\\
8	136\\
8.25	144\\
8.5	152\\
8.75	162\\
9	172\\
9.25	182\\
9.5	190\\
9.75	199\\
10	205\\
};
\addplot [color=red,solid,mark=x,mark options={solid}]
  table[row sep=crcr]{%
1	15\\
1.25	13\\
1.5	17\\
1.75	44\\
2	28\\
2.25	37\\
2.5	36\\
2.75	30\\
3	24\\
3.25	18\\
3.5	12\\
3.75	9\\
4	6\\
4.25	4\\
4.5	4\\
4.75	4\\
5	4\\
5.25	4\\
5.5	5\\
5.75	5\\
6	6\\
6.25	12\\
6.5	6\\
6.75	14\\
7	9\\
7.25	9\\
7.5	10\\
7.75	12\\
8	10\\
8.25	13\\
8.5	14\\
8.75	16\\
9	19\\
9.25	24\\
9.5	24\\
9.75	29\\
10	31\\
};
\end{axis}
\end{tikzpicture}%

%% file: conclusions.tex
\section{Conclusions}\label{s:conclusions}
We have presented preconditioning strategies for the integral forms of both the forward and inverse acoustic scattering problems in two dimensions.

For the forward problem, initially, we extended to the integral equation case the domain decomposition based preconditioning strategies for PDEs: Additive Schwarz, Restricted Additive Schwarz, Additive Harmonic Schwarz and Symmetric Restricted Additive Schwarz. We presented examples comparing the methods using different number of subdomains and size of overlap parameter at different frequencies. The main conclusion for this part is that the convergence of GMRES using the RAS and AHS precondioners is faster than with the other preconditioners. As expected, the convergence is better when using larger overlap, and as the number of subdomains in the partition increases the convergence speed-up deteriorates. 

For the inverse problem, we used the forward problem RC preconditioner to construct the HRC preconditioner. Examples are presented to show the behavior of the HRC preconditioner using different number of subdomains and different size of overlap parameter for the partition of the domain. As we noted with the forward problem preconditioners, the convergence is faster when using larger overlap in the partition of the domain and as the number of subdomains in the partition increases the convergence speed-up worsens significantly. An example showing the scalability of the method for domains with increasing number of points at different frequencies is also presented. Finally, in the last example of the section, the reconstruction of a set of scatterers points using the recursive linearization algorithm is presented. The convergence of the method using the HRC preconditioner is better than the GMRES solver without preconditioner and using the LR preconditioner, which is considered the state-of-the-art for this problem.

The preconditioning strategies presented here are a viable alternative to speed-up the solution of the forward and inverse scattering problems specially when the size of the domain is extremely large, due to their scalability. They also can be easily adapted to three dimensions and to other problems such as electromagnetics.

In the future, we intend to expand the techniques in this article to the continuous case and use them to solve the multifrequency inverse scattering problem for penetrable media for large scale problems in two and three dimensions.

%% file: appendix.tex
\appendix
\section{Recursive Linearization Algorithm}\label{s:appendA}
At low frequencies, the inverse scattering problem is uniquely solvable; however, it presents poor stability, meaning that it is difficult to obtain high resolution of the contrast function. On the other hand, at higher frequencies, the objective function presents multiple minima but is very stable. This trade-off between frequency and stability of the problem forms the basis of the {\bf Recursive Linearization Algorithm}. The {\bf RLA} uses standard frequency continuation to solve a sequence of inverse single-frequency scattering problems at increasing frequencies, using the solution of each problem as the initial guess for the subsequent problem. The summarized description is in Algorithm \ref{alg:rla}.

\begin{algorithm}
\caption{Recursive Linearization Algorithm with Gauss-Newton method (RLA).}
\label{alg:rla}
\begin{algorithmic}[1]
\STATE{{\bf Input:} data $\db(k_j)$ for $j=1,\ldots,Q$ with $k_1  < \dots < k_Q$, initial guess $q_0$, tolerances $\epsilon_1(k_j)$,$\epsilon_2(k_j)$ and maximum number of iterations $N_{it}$.}
\FOR{$j=1,\ldots,Q$}
\STATE{Set $q\coloneqq q_{j-1}$, $\delta q\coloneqq 0$ and $it\coloneqq 0$.}
\WHILE{$\|\db(k_j)-\Fb\left(q\right)\|\geq\epsilon_1(k_j)$ and $it<N_{it}$ and $\delta q\geq\epsilon_2(k_j)$} 
\STATE{Solve $\Hb\, \delta q= \Jb^\ast \left(\db(k_j)-\Fb\left(q\right)\right)$}
\STATE{Update $q\leftarrow q+\delta q$}
\STATE{Update $it\leftarrow it+1$}
\ENDWHILE
\STATE{Set $q_j\coloneqq q$.}
\ENDFOR
\end{algorithmic}
\end{algorithm}

\section{Supplemental numerical results}\label{s:appendB}
We present supplemental results regarding the experiments in the forward and inverse scattering problems. A list of the experiments, their related results and the description of the results is in Table \ref{table:experiments:appendix}.

\begin{table}[H]
\caption{List of supplemental results for the numerical experiments.}\label{table:experiments:appendix}
{\small
\begin{center}
\begin{tabular}{|c||c||l|}
\hline
Experiment & Results & Description of Results \\
\hline\hline
\multirow{5}{*}{F.1} & Table \ref{table:fwd_exp1_sup}     & Relative error $e_{rel}$ for the simulations in Experimentt F.1. \\\cline{2-3}
                               & \multirow{2}{*}{Figure \ref{fig:eig_example_f1_1}} & Eigenvalues in the complex plane of $\Ab$, $\tilde{\Ab}^{-1}_{AS}\Ab$, $\tilde{\Ab}^{-1}_{RAS}\Ab$ and $\tilde{\Ab}^{-1}_{SRAS}\Ab$ \\
                               &							     & when $N_s=16$, $\delta=6$ and $k/(2\pi)=20$. \\\cline{2-3} 
                               & \multirow{2}{*}{Figure \ref{fig:eig_example_f1_2}} & Eigenvalues in the complex plane of $\tilde{\Ab}^{-1}_{RAS}\Ab$ with $N_s=4$ and $16$, \\
                               &							    & $\delta=1$ and $6$, and $k/(2\pi)=20$. \\\hline                               
F.2                          & Table \ref{table:fwd_exp2_sup}     & Relative error $e_{rel}$ for the simulations in Experiments F.2. \\\hline 
F.3 			      & Table \ref{table:fwd_exp3_sup}     & Relative error $e_{rel}$ for the simulations in Experiments F.3. \\\hline
F.4 			      & Table \ref{table:fwd_exp4_sup}     & Relative error $e_{rel}$ for the simulations in Experiments F.4. \\\hline
F.5 			      & Table \ref{table:fwd_exp5_sup}     & Relative error $e_{rel}$ for the simulations in Experiments F.5. \\\hline\hline
\multirow{2}{*}{I.1}  & \multirow{2}{*}{Figure \ref{fig:error_inv_ol_ft}}        & Plots of $\|\tilde{\Ab}^{-1}_{RC}-\Ab^{-1}\|/\|\tilde{\Ab}^{-1}\|$ and $\|\tilde{\Ab}_{RC}-\Ab\|/\|\Ab\|$ at wavenumbers  \\
			      &							     & $k/(2\pi)=5$, $10$ and $20$, with $N_s=16$ and $\delta=3$, $6$ and $9$.\\\hline
\multirow{2}{*}{I.2} & \multirow{2}{*}{Figure \ref{fig:error_inv_sd_ft}}       & Plots of $\|\tilde{\Ab}^{-1}_{RC}-\Ab^{-1}\|/\|\Ab^{-1}\|$ and $\|\tilde{\Ab}_{RC}-\Ab\|/\|\Ab\|$ at wavenumbers \\
			      & 						      & $k/(2\pi)=5$, $10$, and $20$, with $N_s=4$, $16$, $25$, $36$ and $64$, and $\delta=8$. \\\hline 
\multirow{3}{*}{I.3}  & \multirow{2}{*}{Table \ref{tab:error_tf}} & Plots of $\|\tilde{\Ab}^{-1}_{RC}-\Ab^{-1}\|/\|\Ab^{-1}\|$ and $\|\tilde{\Ab}_{RC}-\Ab\|/\|\Ab\|$ for $N=32^2$, $64^2$ \\
			      &								&  and $128^2$ scatterers when $k/(2\pi)=5$ and $20$, with different $N_\lambda$. \\\cline{2-3}
                               & \multirow{1}{*}{Table \ref{tab:inv_scal_error}} & Relative error $e_{rel}$ for the simulations in Experiments I.3. \\\hline
\multirow{4}{*}{I.4}  & \multirow{2}{*}{Figure \ref{fig:Heig_example_I4_k5}} & Singular values of $\Hb$, $\Hb_{LR}$ and $\Hb_{HRC-4}$ using $N_\lambda=140$, and \\
			      &								& $\Hb_{HRC-8}$ using $N_\lambda=40$ when $k/(2\pi)=5$. \\\cline{2-3}
			      & \multirow{2}{*}{Figure \ref{fig:Heig_example_I4_k20}} &  Singular values of $\Hb$, $\Hb_{LR}$ and $\Hb_{HRC-4}$ using $N_\lambda=240$, and \\
			      &								& $\Hb_{HRC-8}$ using $N_\lambda=100$ when $k/(2\pi)=20$. \\\hline
\end{tabular}
\end{center}
}
\end{table}


\begin{table}
\center
\caption{(Experiment F.1 -- Supplemental results) We present the relative error $e_{rel}$ of the GMRES solution with the domain decomposition preconditioners and without using preconditioners. The incoming plane wave has horizontal direction of propagation and frequencies $k/(2\pi)=5$ and $20$. We use a regular grid of $N=64^2$ point scatterers with magnitude given by the function $q_4$. The number of subdomains is $N_s=4$, $16$, and $64$ and the overlap parameter is $\delta=1$ and $6$.}\label{table:fwd_exp1_sup}
\begin{tabular}{|c||c||c||c||*{4}{c|}}\hline
$k/2\pi$ & GMRES & $N_{s}$ & overlap  & AS & RAS & AHS & SRAS \\\hline\hline
\multirow{6}{*}{5} & \multirow{6}{*}{1.6e-08} &\multirow{2}{*}{4} & 1 &	4.3e-10& 2.2e-10& 3.3e-09& 2.2e-10\\\cline{4-8}
                           & &                          						& 6    & 	2.7e-09& 2.3e-09& 2.1e-09& 2.3e-09\\\cline{3-8}
			  & 					 &\multirow{2}{*}{16} & 1    & 	2.4e-09& 5.9e-09& 3.4e-09& 1.1e-08\\\cline{4-8}

                           & &                          						& 6    & 	9.1e-10& 3.4e-09& 3.8e-09& 1.3e-08\\\cline{3-8}
                           & 					&\multirow{2}{*}{64}  & 1    & 	1.4e-08& 1.6e-08& 4.2e-09& 2.3e-08\\\cline{4-8}
                           & &                          						& 6    & 	4.4e-10& 4.0e-10& 2.4e-09& 8.0e-09\\\hline
\multirow{6}{*}{20} & \multirow{6}{*}{1.9e-10} &\multirow{2}{*}{4} & 1  & 	1.8e-08& 4.9e-09& 4.2e-09& 4.9e-09\\\cline{4-8}

                           &  &                          					& 6   & 	7.9e-08& 2.0e-08& 8.3e-09& 2.0e-08\\\cline{3-8}
			  &  						&\multirow{2}{*}{16} & 1 & 2.1e-07& 1.8e-07& 3.3e-08& 2.0e-07\\\cline{4-8}
                           &  &                          					& 6   & 	2.7e-08& 3.0e-08& 8.4e-09& 3.6e-07\\\cline{3-8}
                           &  						&\multirow{2}{*}{64}& 1 & 2.3e-06& 2.1e-07& 1.5e-07& 3.9e-07\\\cline{4-8}
                           &  &                          					& 6    & 	3.7e-08& 8.6e-09& 1.5e-08& 3.2e-07\\\hline

\end{tabular}
\end{table}

\begin{figure}
  \centering
  \begin{subfigure}{.45\linewidth}
\includegraphics[width=.9\textwidth]{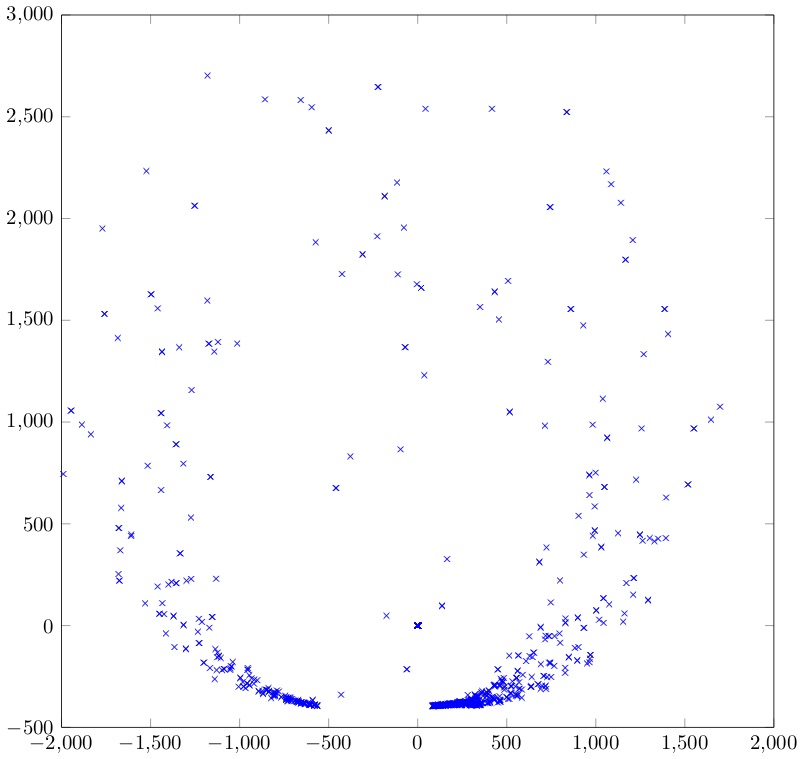}
\caption{Eigenvalues of $\Ab$}
\end{subfigure}
\begin{subfigure}{.45\linewidth}
\includegraphics[width=.9\textwidth]{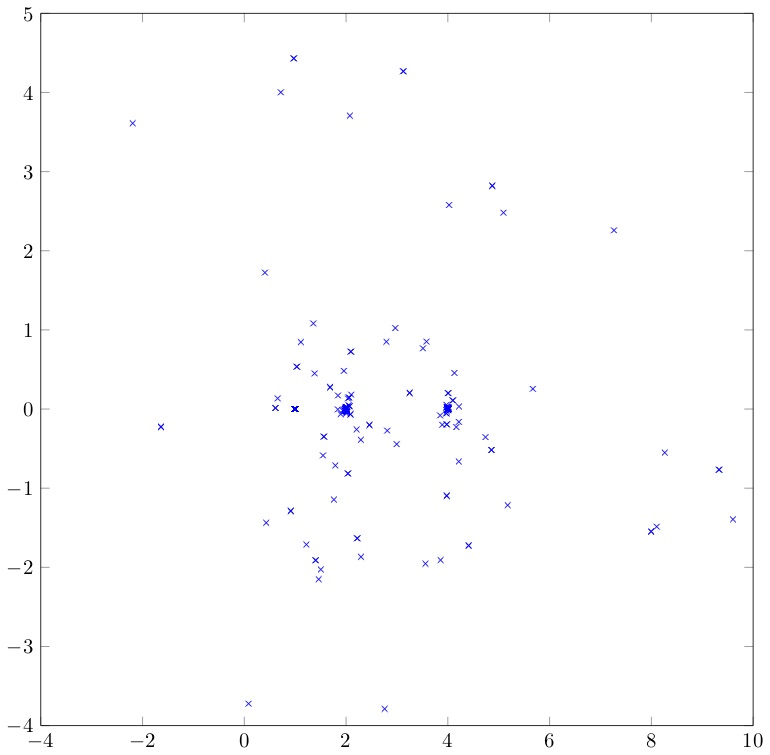}
\caption{Eigenvalues of $\tilde{\Ab}^{-1}_{AS}\Ab$}
\end{subfigure}

  \begin{subfigure}{.45\linewidth}
\includegraphics[width=.9\textwidth]{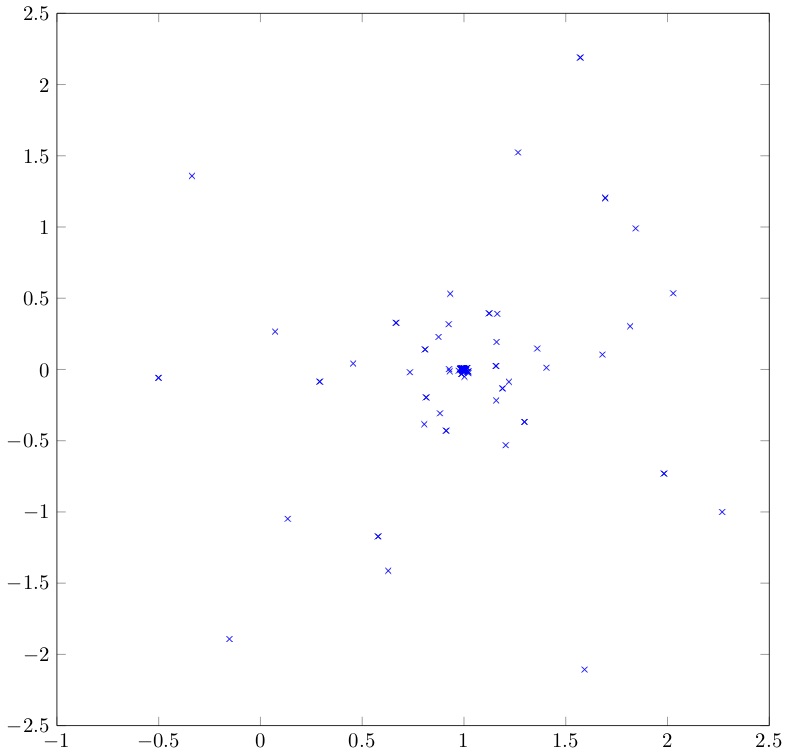}
\caption{Eigenvalues of $\tilde{\Ab}^{-1}_{RAS}\Ab$}
\end{subfigure}
\begin{subfigure}{.45\linewidth}
\includegraphics[width=.9\textwidth]{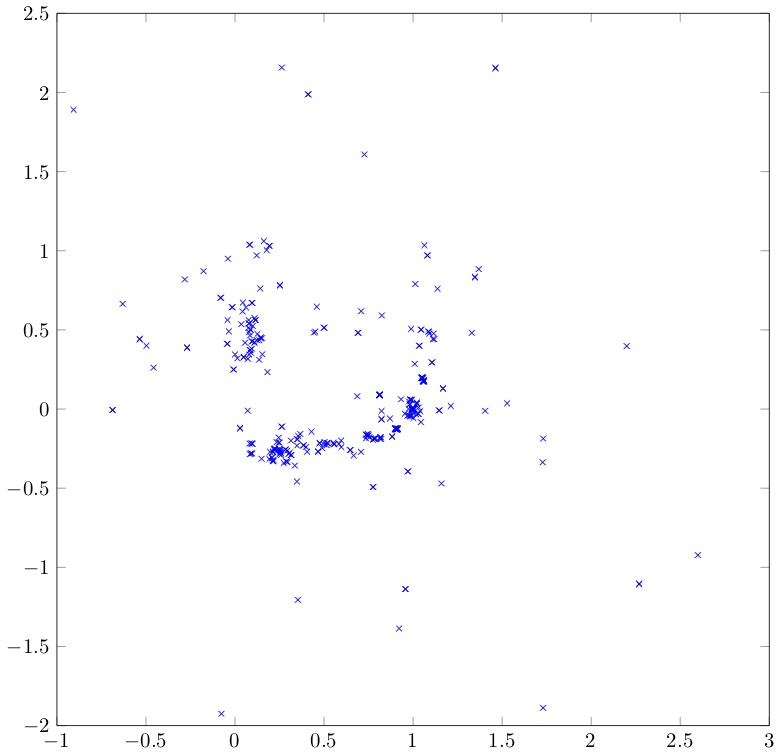}
\caption{Eigenvalues of $\tilde{\Ab}^{-1}_{SRAS}\Ab$}
\end{subfigure}
\caption{(Experiment F.1 -- Supplemental results) Plot in the complex plane of the eigenvalues of (a) $\Ab$, (b) $\tilde{\Ab}^{-1}_{AS}\Ab$, (c) $\tilde{\Ab}^{-1}_{RAS}\Ab$ and (d) $\tilde{\Ab}^{-1}_{SRAS}\Ab$, at $k/(2\pi)=20$, using the domain decomposition parameters $N_s=16$ and $\delta=6$.}
\label{fig:eig_example_f1_1}
\end{figure}

\begin{figure}
  \centering
  \begin{subfigure}{.45\linewidth} 
\includegraphics[width=.9\textwidth]{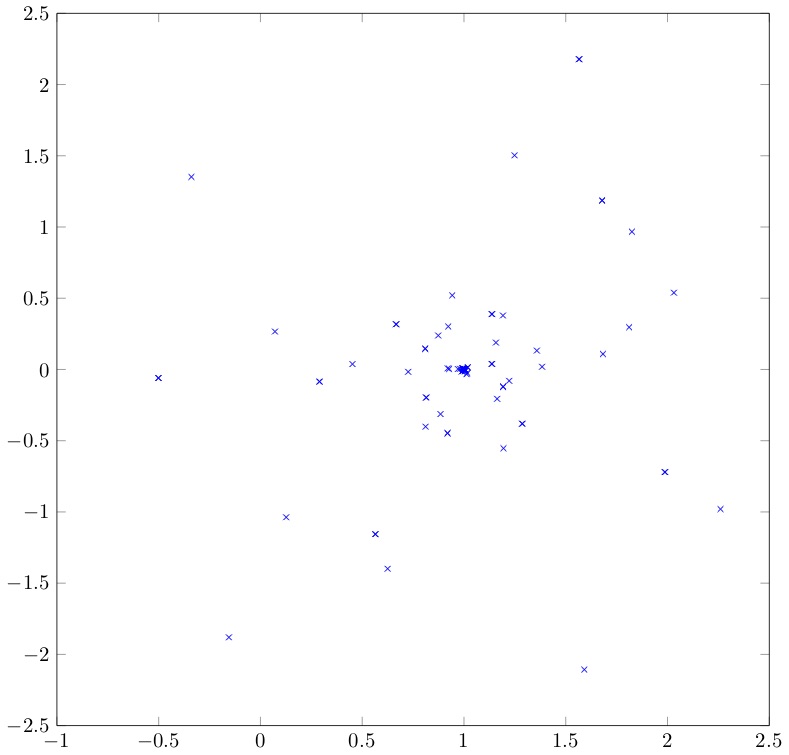}
\caption{Eigenvalues for $N_s =4$ and $\delta=1$}
\end{subfigure}
\begin{subfigure}{.45\linewidth}
\includegraphics[width=.9\textwidth]{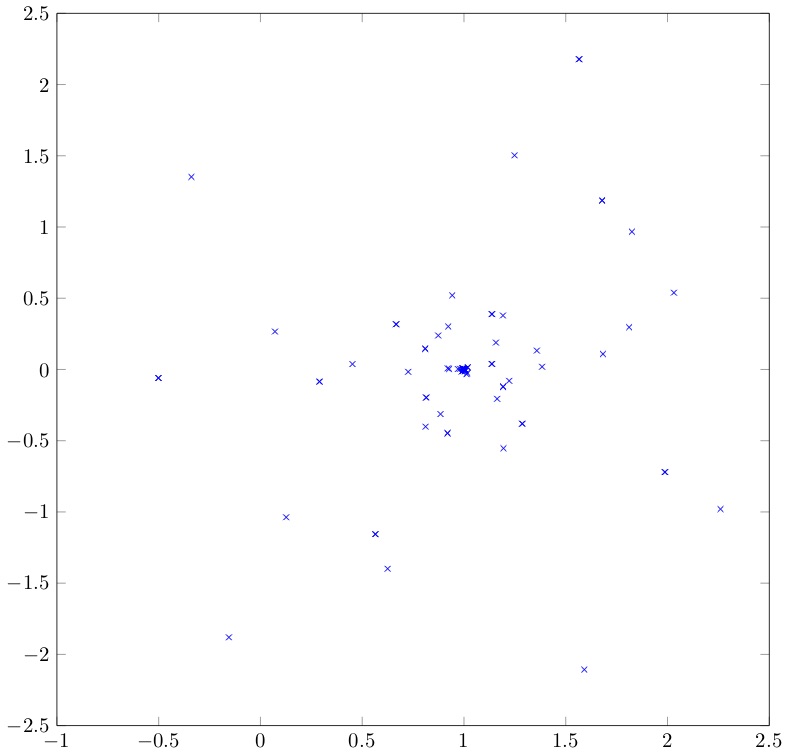}
\caption{Eigenvalues for $N_s =4$ and $\delta=6$}
\end{subfigure}

  \begin{subfigure}{.45\linewidth}
  \includegraphics[width=.9\textwidth]{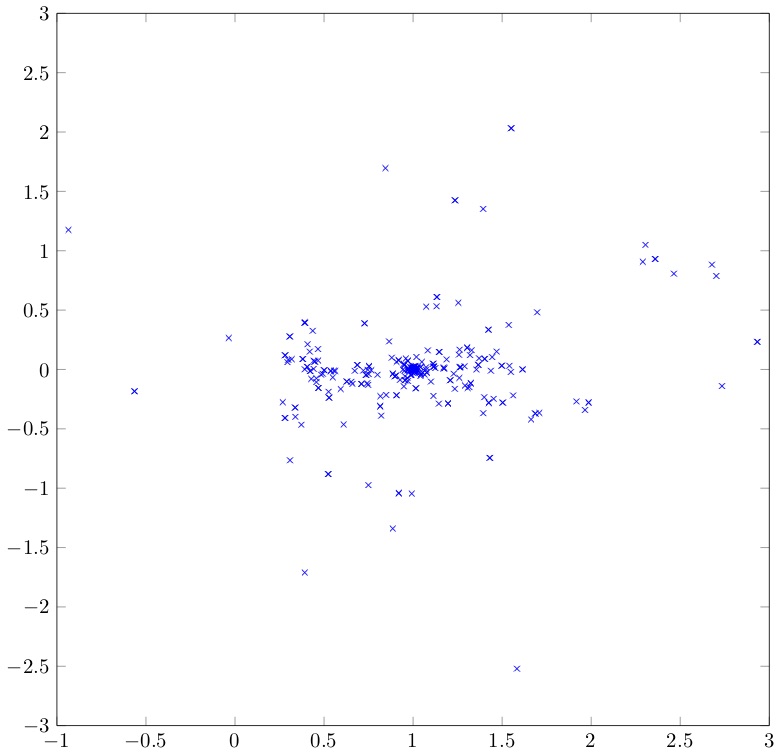}
\caption{Eigenvalues for $N_s =64$ and $\delta=1$}
\end{subfigure}
\begin{subfigure}{.45\linewidth}
\includegraphics[width=.9\textwidth]{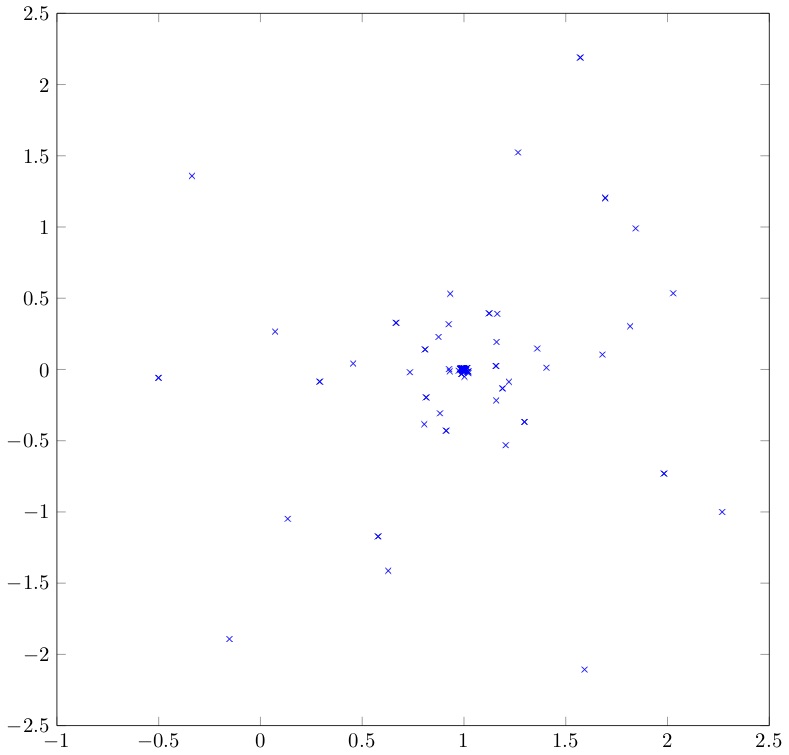}
\caption{Eigenvalues for $N_s =64$ and $\delta=6$}
\end{subfigure}
\caption{(Experiment F.1 -- Supplemental results) Plot in the complex plane of the eigenvalues of $\tilde{\Ab}^{-1}_{RAS}\Ab$ at $k/(2\pi)=20$, using the domain decomposition parameters (a) $N_s=4$ and $\delta=1$, (b) $N_s=4$ and $\delta=6$, (c) $N_s=64$ and $\delta=1$, and (d) $N_s=64$ and $\delta=6$.}
\label{fig:eig_example_f1_2}
\end{figure}


\begin{table}
\center
\caption{(Experiment F.2 -- Supplemental results) We present the relative error $e_{rel}$ of the solution obtained by the iterative method without using preconditioners  and using the preconditioners AS, RAS and AHS. The incoming plane waves have wavenumbers $k/(2\pi)=10$ and $40$. We use a regular grid of $N=64^2$ point scatterers with magnitude given by the function $q_4$. $\mathcal{G}_1$ represents the partition composed of equal sized squares and $\mathcal{G}_2$ is the partition composed of vertical bands.}
\label{table:fwd_exp2_sup}
\begin{tabular}{|c||c||c||c||*{6}{c|}}\hline
\multirow{2}{*}{$k/2\pi$} & \multirow{2}{*}{GMRES} & \multirow{2}{*}{$N_s$} & \multirow{2}{*}{$\delta$} & \multicolumn{2}{c|}{AS} & \multicolumn{2}{c|}{RAS} & \multicolumn{2}{c|}{AHS} \\\cline{5-10}
                                       &               &                                     &                                        & $\mathcal{G}_1$ & $\mathcal{G}_2$  & $\mathcal{G}_1$ & $\mathcal{G}_2$ & $\mathcal{G}_1$ & $\mathcal{G}_2$ \\\hline\hline
\multirow{6}{*}{10} & \multirow{6}{*}{6.1e-10} & \multirow{2}{*}{4}   &   1  & 4.5e-8   & 1.4e-8 & 1.5e-9 & 1.3e-7 & 1.7e-8 & 1.1e-7 \\\cline{4-10}
                              &                                       &                               &   3  & 2.7e-8   & 2.7e-8 & 3.8e-8 & 2.2e-8 & 5.1e-9 & 3.9e-8 \\\cline{3-10}
                              &                                       & \multirow{2}{*}{9}   &   1  & 7.1e-10 & 2.2e-8 & 2.0e-8 & 7.4e-8 & 2.7e-8 & 9.6e-8 \\\cline{4-10}
                              &                                       &                               &   3  & 7.6e-9   & 2.7e-8 & 1.7e-8 & 1.5e-8 & 2.0e-8 & 1.1e-8 \\\cline{3-10}
                              &                                       & \multirow{2}{*}{16} &   1  & 2.2e-8   & 1.1e-7 & 1.3e-7 & 6.7e-8 & 4.4e-8 & 7.4e-8 \\\cline{4-10}
                              &                                       &                               &   3  & 3.2e-8   & 4.9e-8 & 2.3e-8 & 4.2e-8 & 1.5e-8 & 1.3e-8 \\\hline
\multirow{6}{*}{40} & \multirow{6}{*}{7.3e-11} & \multirow{2}{*}{4}   &   1  & 2.0e-8   & 2.8e-7 & 7.7e-8 & 2.9e-8 & 8.4e-8 & 1.9e-7 \\\cline{4-10}
                              &                                       &                               &   3  & 2.5e-8   & 9.8e-8 & 7.4e-8 & 5.5e-7 & 1.2e-7 & 7.2e-7 \\\cline{3-10}
                              &                                       & \multirow{2}{*}{9}   &   1  & 7.8e-7   & 8.2e-7 & 5.7e-7 & 1.4e-7 & 1.0e-6 & 1.2e-6 \\\cline{4-10}
                              &                                       &                               &   3  & 2.3e-7   & 2.3e-7 & 5.5e-7 & 4.5e-7 & 1.2e-6 & 1.7e-7 \\\cline{3-10}
                              &                                       & \multirow{2}{*}{16} &   1  & 8.3e-7   & 1.7e-7 & 2.6e-7 & 1.1e-6 & 3.2e-7 & 6.5e-8 \\\cline{4-10}
                              &                                       &                               &   3  & 6.7e-7   & 1.2e-6 & 2.0e-7 & 5.8e-7 & 6.8e-7 & 3.4e-7 \\\hline
\end{tabular}
\end{table}


\begin{table}
\center
\caption{(Experiment F.3 -- Supplemental results) We present the relative error of the solution $e_{rel}$ obtained by the iterative method without using preconditioners and using the preconditioners AS, RAS and AHS. The incoming plane waves have wavenumbers $k/(2\pi)=10$ and $40$. We use a regular grid of $N=64^2$ point scatterers with magnitude given by the functions $q_4$ and $q_{16}$. The number of subdomains is $N_s=4$ and $16$, and the overlap parameter is $\delta=1$ and $8$.}\label{table:fwd_exp3_sup}
\begin{tabular}{|c||c|c||c||c||*{6}{c|}}\hline
\multirow{2}{*}{$k/2\pi$} & \multicolumn{2}{c||}{GMRES} & \multirow{2}{*}{$N_s$} & \multirow{2}{*}{$\delta$} & \multicolumn{2}{c|}{AS} & \multicolumn{2}{c|}{RAS} & \multicolumn{2}{c|}{AHS} \\\cline{2-3}\cline{6-11}
                                       &  $q_4$       & $q_{16}$         &                                     &                                        & $q_4$ & $q_{16}$          & $q_4$ & $q_{16}$            & $q_4$ & $q_{16}$ \\\hline\hline
\multirow{4}{*}{10}          & \multirow{4}{*}{2.4e-9}  & \multirow{4}{*}{7.4e-10} &\multirow{2}{*}{4}   &   1                  & 7.4e-8 & 3.3e-9 & 1.7e-8 & 3.2e-8 & 1.7e-8 & 3.5e-8 \\\cline{5-11}
                                       &                                  &                                 &                              &   16                           & 5.4e-7 & 9.7e-9 & 1.7e-8 & 1.1e-8 & 1.1e-8 & 1.5e-8 \\\cline{4-11}
                                       &                                  &                                 &\multirow{2}{*}{16} &   1                             & 2.2e-8 & 4.3e-8 & 2.7e-8 & 2.0e-7 & 7.5e-8 & 3.1e-7 \\\cline{5-11}
                                       &                                  &                                 &                              &   16                           & 3.2e-7 & 1.1e-7 & 1.7e-8 & 9.8e-8 & 3.7e-8 & 1.4e-7 \\\hline
\multirow{4}{*}{40}          & \multirow{4}{*}{3.3e-10}  &  \multirow{4}{*}{3.9e-10}&\multirow{2}{*}{4}   &   1                 & 3.6e-6 & 4.3e-7 & 8.0e-7 & 6.4e-7 & 9.1e-7 & 1.4e-7 \\\cline{5-11}
                                       &                                  &                                 &                              &   16                           & 2.6e-7 & 6.3e-7 & 9.2e-7 & 6.6e-7 & 1.2e-6 & 1.5e-6 \\\cline{4-11}
                                       &                                  &                                 &\multirow{2}{*}{16} &   1                             & 4.6e-6 & 3.9e-7 & 1.4e-5 & 1.6e-6 & 3.4e-6 & 1.1e-6 \\\cline{5-11}
                                       &                                  &                                 &                              &   16                           & 9.8e-7 & 4.5e-7 & 7.7e-7 & 1.8e-6 & 1.2e-6 & 4.5e-5 \\\hline
\end{tabular}
\label{table:gmres_dif_scat_error}
\end{table}


\begin{table}
%
\center
\caption{(Experiment F.4 -- Supplemental results) We present the relative error $e_{rel}$ obtained by the iterative method without using preconditioners  and using the preconditioners AS, RAS and AHS. The incoming plane wave has horizontal direction of propagation and frequencies $k/(2\pi)=10$, $20$ and $40$. We use regular grids of $N=64^2$, $128^2$ and $256^2$ point scatterers with magnitude given by the function $q_4$. The number of subdomains is $N_s=16$ and the overlap parameter is $\delta=3$, $6$ and $12$ for $N=64^2$, $N=128^2$ and $N=256^2$ respectively.}\label{table:fwd_exp4_sup}
\begin{tabular}{|c||c||c||c||*{4}{c|}}\hline
$k/2\pi$ & $\sqrt{N}$ & GMRES & $\delta$ & AS & RAS & AHS & SRAS \\\hline\hline
\multirow{3}{*}{10} &                             64 & 6.1e-10& 3&3.2e-08& 2.3e-08& 1.5e-08& 7.1e-08 \\\cline{2-8}
                              & \multirow{1}{*}{128} & \multirow{1}{*}{2.4e-09}  & 6   & 9.5e-08& 1.0e-07& 7.3e-08& 2.3e-07\\\cline{2-8}
                              & \multirow{1}{*}{256} & \multirow{1}{*}{1.1e-08}  & 12 & 2.6e-06& 6.9e-06& 5.0e-06& 8.8e-07\\\hline
\multirow{3}{*}{20} &                             64 & 1.9e-10  & 3  & 4.9e-08& 3.5e-08& 6.5e-09& 1.2e-06 \\\cline{2-8}
                              & \multirow{1}{*}{128} & \multirow{1}{*}{1.7e-09}  & 6   & 2.9e-07& 1.2e-06& 5.8e-07& 6.1e-06\\\cline{2-8}
                              & \multirow{1}{*}{256} & \multirow{1}{*}{6.0e-09}  & 12 & 2.0e-06& 1.6e-05& 1.7e-05& 2.7e-05\\\hline
\multirow{3}{*}{40} &                             64 & 7.3e-11  & 3   & 6.7e-07& 2.0e-07& 6.8e-07& 1.1e-06 \\\cline{2-8}
                              & \multirow{1}{*}{128} & \multirow{1}{*}{3.6e-10}  & 6   & 2.7e-06& 4.5e-07& 4.2e-07& 7.3e-06\\\cline{2-8}
                              & \multirow{1}{*}{256} & \multirow{1}{*}{1.5e-09}  & 12 & 1.5e-04& 1.2e-04& 1.3e-04& 1.6e-04\\\hline
\end{tabular}
\center
\end{table}


\begin{table}
\caption{(Experiment F.5 -- Supplemental results) We present the relative error $e_{rel}$ for GMRES, RAS and the RC preconditioner using $N_\lambda=20$, $40$, $60$ and $80$ at $k=5$ and $N_\lambda=40$, $80$, $120$ and $160$ at $k=20$ . The incoming plane waves have wavenumbers $k/(2\pi)=5$ and $20$ with incidence direction $(1,0)$. We use a regular grid of $N=64^2$ point scatterers with magnitude given by the function $q_4$. The number of subdomains is $N_s=4$ and $16$, and the overlap parameter is $\delta=1$ and $8$.}\label{table:fwd_exp5_sup}
\begin{subtable}{\textwidth}
\caption{Relative error $e_{rel}$ of the scattered field at $k/(2\pi)=5$}
\center
\begin{tabular}{|c||c|c||c||c|c|c|c|}\hline
\multirow{2}{*}{GMRES}       & \multirow{2}{*}{$N_s$} & \multirow{2}{*}{$\delta$} & \multirow{2}{*}{RAS} & \multicolumn{4}{c|}{RC -- $N_\lambda$}  \\ \cline{5-8}
                                             &                                      &                                       &                                    &   20 & 40 & 60 & 80  \\\hline\hline
\multirow{4}{*}{ 4.0e-10}    &\multirow{2}{*}{4}	&   1       & 1.7e-09& 5.1e-12& 6.1e-12& 5.4e-12& 6.5e-12  \\\cline{3-8}
	                                  &                              &   8          & 2.0e-09& 4.5e-12& 4.1e-12& 4.7e-12& 5.3e-12 \\\cline{2-8}
					 &\multirow{2}{*}{16} &   1          & 6.4e-09& 5.6e-10& 2.0e-11& 5.2e-12& 5.8e-12 \\\cline{3-8}
	                                 &                              &   8          & 2.4e-09& 5.1e-12& 4.1e-12& 4.4e-12& 6.0e-12 \\\hline 
\end{tabular}
\end{subtable}

\begin{subtable}{\textwidth}
\caption{Relative error $e_{rel}$ of the scattered field at $k/(2\pi)=20$}
\center
\begin{tabular}{|c||c|c||c||c|c|c|c|}\hline
\multirow{2}{*}{GMRES}       & \multirow{2}{*}{$N_s$} & \multirow{2}{*}{$\delta$} & \multirow{2}{*}{RAS} & \multicolumn{4}{c|}{RC -- $N_\lambda$}  \\ \cline{5-8}                                       
                                             &                                      &                                       &                                    & 40 & 80 & 120 & 160  \\ \hline 
\multirow{4}{*}{ 2.5e-10}   &\multirow{2}{*}{4}  &   1   		& 4.9e-09& 4.2e-09& 2.3e-11& 2.2e-11& 2.1e-11 \\\cline{3-8}
        	        		                  &                              &   8          & 1.1e-08& 4.2e-09& 2.3e-11& 2.4e-11& 2.7e-11 \\\cline{2-8}
        		                          &\multirow{2}{*}{16} &   1          & 1.8e-07& 2.6e-07& 1.8e-07& 6.5e-10& 8.9e-11\\\cline{3-8}
	                                  &                              &   8          & 1.7e-08& 4.2e-09& 2.4e-11& 2.1e-11& 1.9e-11  \\\hline
                                       
\end{tabular}
\end{subtable}
\end{table}


\begin{figure}
\footnotesize
\centering
\begin{subfigure}{.45\textwidth}
\begin{tikzpicture}[scale=0.7]
\begin{semilogyaxis}[xmin=10, xmax=160,domain=1:160,ymin=1e-4, ymax=5e+0,domain=1e-4:5e+0,
    xlabel={$N_\lambda$},
    ylabel={$\frac{\|\tilde{\Ab}^{-1}_{RC}-\Ab^{-1}\|}{\|\Ab^{-1}\|}$},
    grid=major,
    legend pos=north east,
    legend entries={$\delta=3$,$\delta=6$,$\delta=9$},
]
\addplot coordinates {
(10,5.208301e-01) (20,3.505532e-01) (30,2.908219e-01) (40,1.727415e-01) (50,9.841513e-02) (60,7.055136e-02)
(70,5.770865e-02) (80,4.726983e-02) (90,3.931140e-02) (100,3.393519e-02) (110,2.768438e-02) (120,1.929340e-02)
(130,1.660990e-02) (140,1.251251e-02) (150,9.269606e-03) (160,7.096203e-03)
};
\addplot coordinates {
(10,4.125958e-01) (20,1.374619e-01) (30,7.148799e-02) (40,4.552846e-02) (50,2.231503e-02) (60,1.466071e-02)
(70,1.037704e-02) (80,8.608461e-03) (90,5.712878e-03)
};
\addplot coordinates {
(10,4.030965e-01) (20,2.410227e-01) (30,1.162914e-01) (40,4.245485e-02) (50,9.893943e-03) (60,2.955278e-03)
};
\end{semilogyaxis}
\end{tikzpicture}
\caption{$\frac{\|\tilde{\Ab}^{-1}_{RC}-\Ab^{-1}\|}{\|\Ab^{-1}\|}$ at $k/(2\pi)=5$}\label{fig:error_inv_ol_ft_1}
\end{subfigure}
\begin{subfigure}{.45\textwidth}
\begin{tikzpicture}[scale=0.7]
\begin{semilogyaxis}[xmin=10, xmax=160,domain=1:160,ymin=1e-4, ymax=1e+0,domain=1e-4:1e+0,
    xlabel={$N_\lambda$},
    ylabel={$\frac{\|\tilde{\Ab}_{RC}-\Ab\|}{\|\Ab\|}$},
    grid=major,
    legend pos=north east,
    legend entries={$\delta=3$,$\delta=6$,$\delta=9$},
]
\addplot coordinates {
(10,3.383569e-01) (20,1.429379e-01) (30,7.822588e-02) (40,4.119527e-02) (50,2.553344e-02) (60,1.728831e-02)
(70,1.212183e-02) (80,9.206517e-03) (90,6.942970e-03) (100,5.172318e-03) (110,3.894971e-03) (120,2.961039e-03)
(130,2.144714e-03) (140,1.635361e-03) (150,1.264931e-03) (160,9.329382e-04)
};
\addplot coordinates {
(10,2.899069e-01) (20,5.684602e-02) (30,1.637049e-02) (40,8.270498e-03) (50,4.626358e-03) (60,3.047981e-03)
(70,2.039835e-03) (80,1.388437e-03) (90,9.271018e-04)
};
\addplot coordinates {
(10,2.322861e-01) (20,6.511837e-02) (30,2.932501e-02) (40,1.029972e-02) (50,2.253141e-03) (60,6.448004e-04)
};
\end{semilogyaxis}
\end{tikzpicture}
\caption{$\frac{\|\tilde{\Ab}_{RC}-\Ab\|}{\|\Ab\|}$ at $k/(2\pi)=5$}\label{fig:error_inv_ol_ft_2}
\end{subfigure}

\begin{subfigure}{.45\textwidth}
\begin{tikzpicture}[scale=0.7]
\begin{semilogyaxis}[xmin=10, xmax=200,domain=1:200,ymin=1e-4, ymax=5e+0,domain=1e-4:5e+0,
    xlabel={$N_\lambda$},
    ylabel={$\frac{\|\tilde{\Ab}^{-1}_{RC}-\Ab^{-1}\|}{\|\Ab^{-1}\|}$},
    grid=major,
    legend pos=north east,
    legend entries={$\delta=3$,$\delta=6$,$\delta=9$},
]
\addplot coordinates {
(10,9.347583e-01) (20,1.708244e+00) (30,1.160982e+00) (40,1.762622e+00) (50,1.918648e+00) (60,4.060002e-01)
(70,5.203111e-01) (80,1.755629e-01) (90,9.659678e-02) (100,1.036847e-01) (110,1.328356e-01) (120,5.215775e-02)
(130,3.359188e-02) (140,2.671753e-02) (150,2.757707e-02) (160,2.280855e-02) (170,1.191053e-02) (180,6.974879e-03)
(190,7.683879e-03) (200,7.296567e-03)
};
\addplot coordinates {
(10,9.273228e-01) (20,7.478443e-01) (30,3.985020e-01) (40,1.118984e-01) (50,1.442538e-01) (60,5.139933e-02)
(70,5.008204e-02) (80,2.960799e-02) (90,1.407799e-02) (100,1.295516e-02) (110,2.015281e-02) (120,5.296076e-03)
};
\addplot coordinates {
(10,1.106471e+00) (20,6.513435e-01) (30,1.136197e+00) (40,5.059565e-01) (50,1.743176e-01) (60,9.006806e-02)
(70,2.141996e-02) (80,1.446481e-02) (90,1.110764e-03) 
};
\end{semilogyaxis}
\end{tikzpicture}
\caption{$\frac{\|\tilde{\Ab}^{-1}_{RC}-\Ab^{-1}\|}{\|\Ab^{-1}\|}$ at $k/(2\pi)=10$}\label{fig:error_inv_ol_ft_3}
\end{subfigure}
\begin{subfigure}{.45\textwidth}
\begin{tikzpicture}[scale=0.7]
\begin{semilogyaxis}[xmin=10, xmax=200,domain=1:200,ymin=1e-4, ymax=1e+0,domain=1e-4:1e+0,
    xlabel={$N_\lambda$},
    ylabel={$\frac{\|\tilde{\Ab}_{RC}-\Ab\|}{\|\Ab\|}$},
    grid=major,
    legend pos=north east,
    legend entries={$\delta=3$,$\delta=6$,$\delta=9$},
]
\addplot coordinates {
(10,5.025465e-01) (20,2.841442e-01) (30,1.482301e-01) (40,1.026126e-01) (50,7.702646e-02) (60,5.461248e-02)
(70,3.968905e-02) (80,2.965756e-02) (90,2.268346e-02) (100,1.731081e-02) (110,1.252673e-02) (120,9.203962e-03)
(130,6.694339e-03) (140,4.947120e-03) (150,3.687833e-03) (160,2.764418e-03) (170,2.064629e-03) (180,1.485464e-03)
(190,1.031776e-03) (200,7.616888e-04)
};
\addplot coordinates {
(10,4.870694e-01) (20,2.582281e-01) (30,9.800719e-02) (40,2.487625e-02) (50,9.209166e-03) (60,6.381139e-03)
(70,4.577490e-03) (80,3.476317e-03) (90,2.574273e-03) (100,1.890909e-03) (110,1.258196e-03) (120,8.127462e-04)
};
\addplot coordinates {
(10,4.088043e-01) (20,2.353055e-01) (30,1.101887e-01) (40,5.470395e-02) (50,2.689806e-02) (60,1.166118e-02)
(70,4.273885e-03) (80,1.114914e-03) (90,3.436526e-04) 
};
\end{semilogyaxis}
\end{tikzpicture}
\caption{$\frac{\|\tilde{\Ab}_{RC}-\Ab\|}{\|\Ab\|}$ at $k/(2\pi)=10$}\label{fig:error_inv_ol_ft_4}
\end{subfigure}

\begin{subfigure}{.45\textwidth}
\begin{tikzpicture}[scale=0.7]
\begin{semilogyaxis}[xmin=10, xmax=280,domain=1:280,ymin=1e-4, ymax=5e+0,domain=1e-4:5e+0,
    xlabel={$N_\lambda$},
    ylabel={$\frac{\|\tilde{\Ab}^{-1}_{RC}-\Ab^{-1}\|}{\|\Ab^{-1}\|}$},
    grid=major,
    legend pos=north east,
    legend entries={$\delta=3$,$\delta=6$,$\delta=9$},
]
\addplot coordinates {
(10,8.989273e-01) (20,8.942845e-01) (30,8.805195e-01) (40,1.644881e+00) (50,8.904516e-01) (60,4.217173e-01)
(70,3.677635e-01) (80,3.851082e-01) (90,3.459411e-01) (100,2.745558e-01) (110,2.109618e-01) (120,1.833317e-01)
(130,1.815177e-01) (140,1.762080e-01) (150,1.612583e-01) (160,1.593467e-01) (170,1.534439e-01) (180,1.233958e-01)
(190,1.028904e-01)(200,5.603056e-02)(210,3.436017e-02)(220,2.370118e-02)(230,2.076777e-02)(240,1.811183e-02)
(250,1.339763e-02)(260,1.214005e-02)(270,7.513163e-03)(280,5.074680e-03)
};
\addplot coordinates {
(10,9.000449e-01) (20,8.942831e-01) (30,8.903315e-01) (40,1.959974e+00) (50,1.235974e+00) (60,4.288238e-01)
(70,1.197226e-01) (80,3.459629e-02) (90,2.272157e-02) (100,2.146363e-02) (110,1.733456e-02) (120,1.268904e-02)
(130,1.093797e-02) (140,1.020040e-02) (150,9.303572e-03) (160,7.655146e-03) (170,7.098124e-03) (180,6.210836e-03)
};
\addplot coordinates {
(10,9.043571e-01) (20,9.049854e-01) (30,9.758641e-01) (40,9.137026e-01) (50,5.433642e-01) (60,2.843520e-01)
(70,2.194978e-01) (80,1.855310e-01) (90,1.042469e-01) (100,4.498632e-02) (110,2.795874e-02) (120,1.213023e-02)
(130,3.029756e-03)
};
\end{semilogyaxis}
\end{tikzpicture}
\caption{$\frac{\|\tilde{\Ab}^{-1}_{RC}-\Ab^{-1}\|}{\|\Ab^{-1}\|}$ at $k/(2\pi)=20$}\label{fig:error_inv_ol_ft_5}
\end{subfigure}
\begin{subfigure}{.45\textwidth}
\begin{tikzpicture}[scale=0.7]
\begin{semilogyaxis}[xmin=10, xmax=280,domain=1:280,ymin=1e-4, ymax=1e+0,domain=1e-4:1e+0,
    xlabel={$N_\lambda$},
    ylabel={$\frac{\|\tilde{\Ab}_{RC}-\Ab\|}{\|\Ab\|}$},
    grid=major,
    legend pos=north east,
    legend entries={$\delta=3$,$\delta=6$,$\delta=9$},
]
\addplot coordinates {
(10,5.931114e-01) (20,4.940502e-01) (30,3.872922e-01) (40,2.729661e-01) (50,1.768531e-01) (60,1.137785e-01)
(70,9.262586e-02) (80,8.292922e-02) (90,7.433133e-02) (100,6.556885e-02) (110,5.764420e-02) (120,5.018116e-02)
(130,4.326327e-02) (140,3.718763e-02) (150,3.138133e-02) (160,2.578856e-02) (170,2.039616e-02) (180,1.554589e-02)
(190,1.150253e-02) (200,8.466861e-03) (210,6.595299e-03) (220,5.038839e-03) (230,3.800461e-03) (240,2.847790e-03)
(250,2.045456e-03) (260,1.492891e-03) (270,1.085145e-03) (280,7.363796e-04)
};
\addplot coordinates {
(10,5.839367e-01) (20,4.838642e-01) (30,3.749859e-01) (40,2.551771e-01) (50,1.495408e-01) (60,6.722001e-02)
(70,2.542829e-02) (80,9.479471e-03) (90,6.796451e-03) (100,5.435378e-03) (110,4.383688e-03) (120,3.662372e-03)
(130,3.074330e-03) (140,2.539485e-03) (150,2.042197e-03) (160,1.539478e-03) (170,1.080006e-03) (180,6.673347e-04)
};
\addplot coordinates {
(10,5.827362e-01) (20,4.668028e-01) (30,3.593322e-01) (40,2.533865e-01) (50,1.654588e-01) (60,9.850985e-02)
(70,5.736401e-02) (80,3.335871e-02) (90,2.063512e-02) (100,8.492702e-03) (110,3.658316e-03) (120,1.582804e-03)
(130,5.675289e-04) 
};
\end{semilogyaxis}
\end{tikzpicture}
\caption{$\frac{\|\tilde{\Ab}_{RC}-\Ab\|}{\|\Ab\|}$ at $k/(2\pi)=20$}\label{fig:error_inv_ol_ft_6}
\end{subfigure}
 \caption{(Experiment I.1 -- Supplemental results) Plots of $\|\tilde{\Ab}^{-1}_{RC}-\Ab^{-1}\|/\|\Ab^{-1}\|$ and $\|\tilde{\Ab}_{RC}-\Ab\|/\|\Ab\|$ for different $N_\lambda$ at wavenumbers: (a) and (b) $k/(2\pi)=5$, (c) and (d) $k/(2\pi)=10$, and (e) and (f) $k/(2\pi)=20$. In each plot, each line represents the error when using overlap parameter $\delta=3$, $6$ and $9$ points. The number of subdomains is $N_s=16$.}\label{fig:error_inv_ol_ft}
\end{figure}
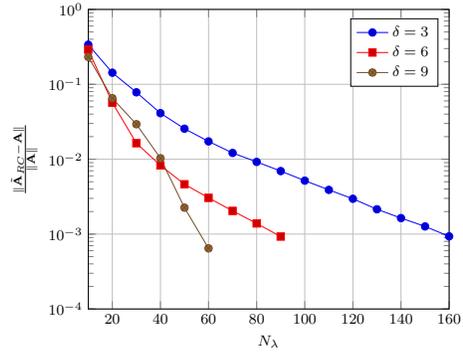
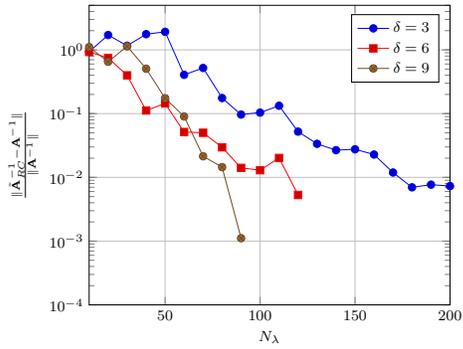
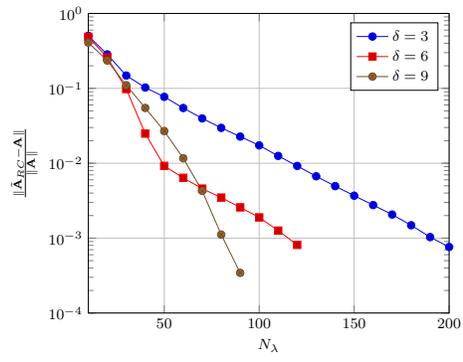
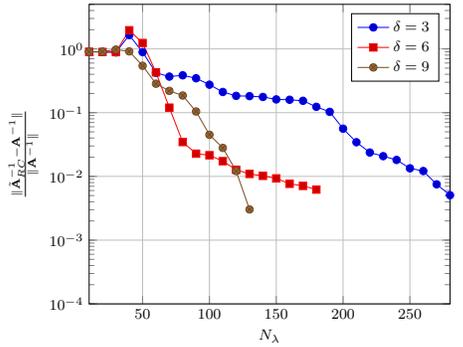
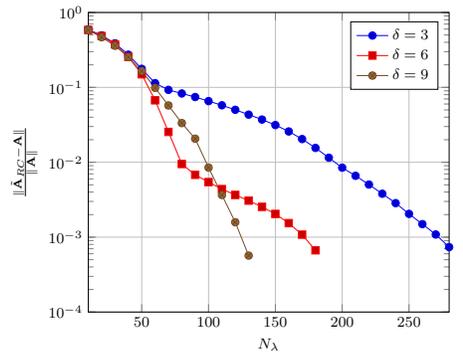


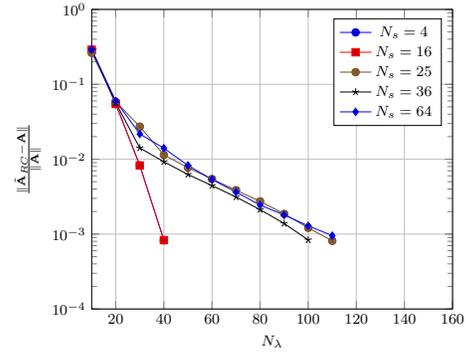
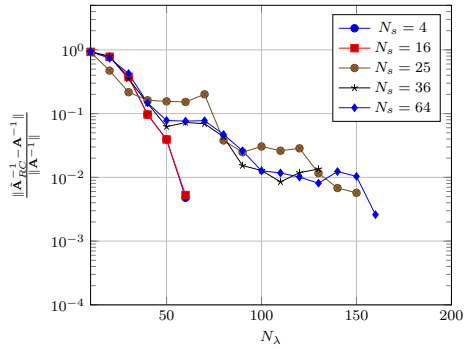
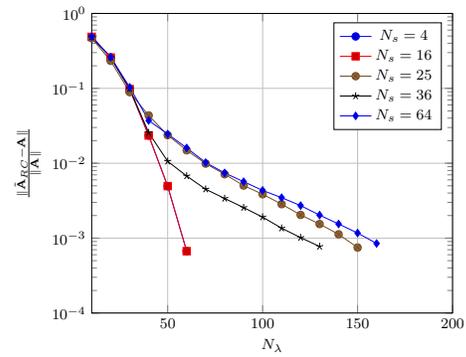
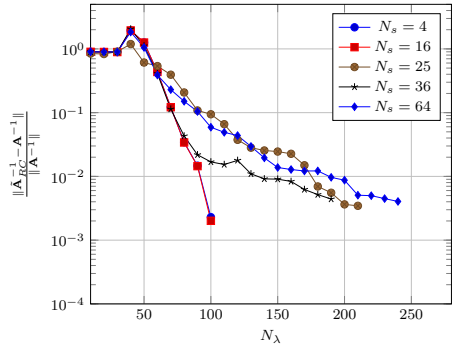
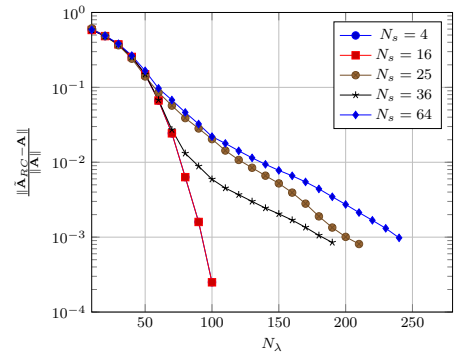
\begin{figure}
\footnotesize
\centering
\begin{subfigure}{.45\textwidth}
\begin{tikzpicture}[scale=0.7]
\begin{semilogyaxis}[xmin=10, xmax=160,domain=1:160,ymin=1e-4, ymax=5e+0,domain=1e-4:5e+0,
    xlabel={$N_\lambda$},
    ylabel={$\frac{\|\tilde{\Ab}^{-1}_{RC}-\Ab^{-1}\|}{\|\Ab^{-1}\|}$},
    grid=major,
    legend pos=north east,
    legend entries={$N_s=4$,$N_s=16$,$N_s=25$,$N_s=36$,$N_s=64$},
]
\addplot coordinates {
(10,4.039273e-01) (20,1.170668e-01) (30,2.233715e-02) (40,3.249538e-03)
};
\addplot coordinates {
(10,4.039060e-01) (20,1.170875e-01) (30,2.235023e-02) (40,3.270202e-03)
};
\addplot coordinates {
(10,3.840428e-01) (20,2.008738e-01) (30,1.004179e-01) (40,4.454047e-02) (50,3.479271e-02) (60,2.853628e-02) 
(70,2.094017e-02) (80,1.435869e-02) (90,1.125832e-02) (100,8.424633e-03) (110,4.960320e-03) 
};
\addplot coordinates {
(10,4.050317e-01) (20,1.236065e-01) (30,4.342313e-02) (40,3.586449e-02) (50,2.477388e-02) (60,2.150342e-02) 
(70,1.808108e-02) (80,1.357595e-02) (90,8.420222e-03) (100,6.257607e-03)
};
\addplot coordinates {
(10,4.060477e-01) (20,1.307975e-01) (30,5.889213e-02) (40,4.420241e-02) (50,3.017183e-02) (60,2.019317e-02)
(70,1.521531e-02) (80,1.030319e-02) (90,7.457172e-03) (100,5.599460e-03) (110,4.262291e-03) 
};
\end{semilogyaxis}
\end{tikzpicture}
\caption{$\frac{\|\tilde{\Ab}^{-1}_{RC}-\Ab^{-1}\|}{\|\Ab^{-1}\|}$ at $k/(2\pi)=5$}\label{fig:error_inv_sd_ft_1}
\end{subfigure}
\begin{subfigure}{.45\textwidth}
\begin{tikzpicture}[scale=0.7]
\begin{semilogyaxis}[xmin=10, xmax=160,domain=1:160,ymin=1e-4, ymax=1e+0,domain=1e-4:1e+0,
    xlabel={$N_\lambda$},
    ylabel={$\frac{\|\tilde{\Ab}_{RC}-\Ab\|}{\|\Ab\|}$},
    grid=major,
    legend pos=north east,
    legend entries={$N_s=4$,$N_s=16$,$N_s=25$,$N_s=36$,$N_s=64$},
]
\addplot coordinates {
(10,2.892341e-01) (20,5.471599e-02) (30,8.248438e-03) (40,8.293137e-04)
};
\addplot coordinates {
(10,2.892403e-01) (20,5.471956e-02) (30,8.249813e-03) (40,8.291689e-04)
};
\addplot coordinates {
(10,2.631831e-01) (20,5.872953e-02) (30,2.725122e-02) (40,1.142284e-02) (50,7.716618e-03) (60,5.440595e-03) 
(70,3.858538e-03) (80,2.734181e-03) (90,1.862698e-03) (100,1.210224e-03) (110,8.121585e-04)
};
\addplot coordinates {
(10,2.896076e-01) (20,5.641939e-02) (30,1.413642e-02) (40,9.235154e-03) (50,6.257580e-03) (60,4.438217e-03) 
(70,3.129285e-03) (80,2.120000e-03) (90,1.386126e-03) (100,8.356812e-04)
};
\addplot coordinates {
(10,2.908993e-01) (20,5.975884e-02) (30,2.171554e-02) (40,1.403383e-02) (50,8.255592e-03) (60,5.369550e-03)
(70,3.643082e-03) (80,2.439101e-03) (90,1.801348e-03) (100,1.288832e-03) (110,9.542202e-04)
};
\end{semilogyaxis}
\end{tikzpicture}
\caption{$\frac{\|\tilde{\Ab}_{RC}-\Ab\|}{\|\Ab\|}$ at $k/(2\pi)=5$}\label{fig:error_inv_sd_ft_2}
\end{subfigure}

\begin{subfigure}{.45\textwidth}
\begin{tikzpicture}[scale=0.7]
\begin{semilogyaxis}[xmin=10, xmax=200,domain=1:200,ymin=1e-4, ymax=5e+0,domain=1e-4:5e+0,
    xlabel={$N_\lambda$},
    ylabel={$\frac{\|\tilde{\Ab}^{-1}_{RC}-\Ab^{-1}\|}{\|\Ab^{-1}\|}$},
    grid=major,
    legend pos=north east,
    legend entries={$N_s=4$,$N_s=16$,$N_s=25$,$N_s=36$,$N_s=64$},
]
\addplot coordinates {
(10,9.306368e-01) (20,7.780412e-01) (30,3.766990e-01) (40,9.656308e-02) (50,3.939553e-02) (60,4.791344e-03)
};
\addplot coordinates {
(10,9.306328e-01) (20,7.785069e-01) (30,3.770790e-01) (40,9.671949e-02) (50,3.928909e-02) (60,5.281532e-03)
};
\addplot coordinates {
(10,9.189302e-01) (20,4.729417e-01) (30,2.168243e-01) (40,1.619401e-01) (50,1.556937e-01) (60,1.526676e-01) 
(70,2.008341e-01) (80,3.792843e-02) (90,2.478205e-02) (100,3.048216e-02) (110,2.616000e-02) (120,2.852909e-02) 
(130,1.163498e-02) (140,6.808715e-03) (150,5.722375e-03) 
};
\addplot coordinates {
(10,9.352570e-01) (20,7.924296e-01) (30,3.501261e-01) (40,1.450787e-01) (50,6.249975e-02) (60,7.247613e-02) 
(70,6.962595e-02) (80,4.316030e-02) (90,1.532382e-02) (100,1.281753e-02) (110,8.514014e-03) (120,1.183270e-02) 
(130,1.352409e-02)
};
\addplot coordinates {
(10,9.381889e-01) (20,7.380760e-01) (30,4.256570e-01) (40,1.473857e-01) (50,7.803782e-02) (60,7.627879e-02)
(70,7.757184e-02) (80,4.706252e-02) (90,2.614307e-02) (100,1.266773e-02) (110,1.172746e-02) (120,1.009579e-02)
(130,8.101820e-03) (140,1.234450e-02) (150,1.035647e-02) (160,2.601058e-03)
};
\end{semilogyaxis}
\end{tikzpicture}
\caption{$\frac{\|\tilde{\Ab}^{-1}_{RC}-\Ab^{-1}\|}{\|\Ab^{-1}\|}$ at $k/(2\pi)=10$}\label{fig:error_inv_sd_ft_3}
\end{subfigure}
\begin{subfigure}{.45\textwidth}
\begin{tikzpicture}[scale=0.7]
\begin{semilogyaxis}[xmin=10, xmax=200,domain=1:200,ymin=1e-4, ymax=1e+0,domain=1e-4:1e+0,
    xlabel={$N_\lambda$},
    ylabel={$\frac{\|\tilde{\Ab}_{RC}-\Ab\|}{\|\Ab\|}$},
    grid=major,
    legend pos=north east,
    legend entries={$N_s=4$,$N_s=16$,$N_s=25$,$N_s=36$,$N_s=64$},
]
\addplot coordinates {
(10,4.869685e-01) (20,2.580756e-01) (30,9.760829e-02) (40,2.344920e-02) (50,4.945121e-03) (60,6.696709e-04)
};
\addplot coordinates {
(10,4.869906e-01) (20,2.580816e-01) (30,9.760798e-02) (40,2.344953e-02) (50,4.945888e-03) (60,6.697620e-04)
};
\addplot coordinates {
(10,4.672614e-01) (20,2.329249e-01) (30,8.908173e-02) (40,4.343666e-02) (50,2.370707e-02) (60,1.488311e-02) 
(70,9.885822e-03) (80,7.128518e-03) (90,5.037522e-03) (100,3.853294e-03) (110,2.838959e-03) (120,2.040565e-03) 
(130,1.541962e-03) (140,1.125048e-03) (150,7.492880e-04) 
};
\addplot coordinates {
(10,4.871432e-01) (20,2.583697e-01) (30,9.825946e-02) (40,2.561158e-02) (50,1.060761e-02) (60,6.754994e-03) 
(70,4.500319e-03) (80,3.385326e-03) (90,2.557125e-03) (100,1.907139e-03) (110,1.357013e-03) (120,1.016193e-03)
(130,7.727455e-04) 
};
\addplot coordinates {
(10,4.886042e-01) (20,2.608583e-01) (30,1.032796e-01) (40,3.724323e-02) (50,2.452482e-02) (60,1.600795e-02)
(70,1.021047e-02) (80,7.423878e-03) (90,5.671452e-03) (100,4.334852e-03) (110,3.467627e-03) (120,2.717651e-03)
(130,2.028074e-03) (140,1.544641e-03) (150,1.169294e-03) (160,8.487230e-04) 
};
\end{semilogyaxis}
\end{tikzpicture}
\caption{$\frac{\|\tilde{\Ab}_{RC}-\Ab\|}{\|\Ab\|}$ at $k/(2\pi)=10$}\label{fig:error_inv_sd_ft_4}
\end{subfigure}

\begin{subfigure}{.45\textwidth}
\begin{tikzpicture}[scale=0.7]
\begin{semilogyaxis}[xmin=10, xmax=280,domain=1:280,ymin=1e-4, ymax=5e+0,domain=1e-4:5e+0,
    xlabel={$N_\lambda$},
    ylabel={$\frac{\|\tilde{\Ab}^{-1}_{RC}-\Ab^{-1}\|}{\|\Ab^{-1}\|}$},
    grid=major,
    legend pos=north east,
    legend entries={$N_s=4$,$N_s=16$,$N_s=25$,$N_s=36$,$N_s=64$},
]
\addplot coordinates {
(10,8.999857e-01) (20,8.947404e-01) (30,8.921807e-01) (40,1.949749e+00) (50,1.258379e+00) (60,4.348366e-01)
(70,1.218250e-01) (80,3.437004e-02) (90,1.470824e-02) (100,2.295859e-03) 
};
\addplot coordinates {
(10,8.999770e-01) (20,8.947446e-01) (30,8.922398e-01) (40,1.947028e+00) (50,1.257346e+00) (60,4.344370e-01)
(70,1.214905e-01) (80,3.405689e-02) (90,1.448083e-02) (100,2.010051e-03)
};
\addplot coordinates {
(10,8.455451e-01) (20,8.309714e-01) (30,9.060578e-01) (40,1.191394e+00) (50,6.112706e-01) (60,5.374652e-01) 
(70,3.937537e-01) (80,2.071560e-01) (90,1.081608e-01) (100,9.432927e-02) (110,6.588901e-02) (120,3.746177e-02) 
(130,2.818890e-02) (140,2.549079e-02) (150,2.451800e-02) (160,2.261516e-02) (170,1.494612e-02) (180,6.951674e-03)
(190,5.537136e-03) (200,3.638279e-03)  (210,3.454920e-03)
};
\addplot coordinates {
(10,8.999483e-01) (20,8.950178e-01) (30,8.935555e-01) (40,2.051003e+00) (50,1.114083e+00) (60,4.091542e-01) 
(70,1.135682e-01) (80,4.303806e-02) (90,2.182953e-02) (100,1.687382e-02) (110,1.536998e-02) (120,1.783047e-02) 
(130,1.097858e-02) (140,9.141661e-03) (150,9.051560e-03) (160,8.333749e-03) (170,6.234958e-03) (180,5.179919e-03) 
(190,4.414493e-03)
};
\addplot coordinates {
(10,9.020168e-01) (20,8.963034e-01) (30,8.956339e-01) (40,1.838463e+00) (50,1.053124e+00) (60,3.888124e-01)
(70,2.294952e-01) (80,1.505668e-01) (90,1.042654e-01) (100,5.902239e-02) (110,4.911457e-02) (120,4.373729e-02)
(130,2.951304e-02) (140,1.953020e-02) (150,1.379273e-02) (160,1.281729e-02) (170,1.216173e-02) (180,1.218698e-02)
(190,9.664138e-03) (200,8.712654e-03) (210,5.069434e-03) (220,4.982937e-03) (230,4.484835e-03) (240,4.043817e-03)
};
\end{semilogyaxis}
\end{tikzpicture}
\caption{$\frac{\|\tilde{\Ab}^{-1}_{RC}-\Ab^{-1}\|}{\|\Ab^{-1}\|}$ at $k/(2\pi)=20$}\label{fig:error_inv_sd_ft_5}
\end{subfigure}
\begin{subfigure}{.45\textwidth}
\begin{tikzpicture}[scale=0.7]
\begin{semilogyaxis}[xmin=10, xmax=280,domain=1:280,ymin=1e-4, ymax=1e+0,domain=1e-4:1e+0,
    xlabel={$N_\lambda$},
    ylabel={$\frac{\|\tilde{\Ab}_{RC}-\Ab\|}{\|\Ab\|}$},
    grid=major,
    legend pos=north east,
    legend entries={$N_s=4$,$N_s=16$,$N_s=25$,$N_s=36$,$N_s=64$},
]
\addplot coordinates {
(10,5.838588e-01) (20,4.837760e-01) (30,3.748901e-01) (40,2.550454e-01) (50,1.493291e-01) (60,6.678851e-02)
(70,2.432874e-02) (80,6.327458e-03) (90,1.596139e-03) (100,2.499537e-04) 
};
\addplot coordinates {
(10,5.838758e-01) (20,4.837910e-01) (30,3.748986e-01) (40,2.550498e-01) (50,1.493312e-01) (60,6.678946e-02)
(70,2.432928e-02) (80,6.327461e-03) (90,1.596161e-03) (100,2.493708e-04)
};
\addplot coordinates {
(10,6.140257e-01) (20,4.849088e-01) (30,3.650730e-01) (40,2.413538e-01) (50,1.399060e-01) (60,8.555059e-02) 
(70,5.705546e-02) (80,3.889042e-02) (90,2.826009e-02) (100,2.036236e-02) (110,1.434206e-02) (120,1.074215e-02)
(130,8.445572e-03) (140,6.621519e-03) (150,5.223647e-03) (160,3.933593e-03) (170,2.801801e-03) (180,1.895330e-03)
(190,1.343492e-03) (200,1.008633e-03) (210,8.096860e-04)
};
\addplot coordinates {
(10,5.841037e-01) (20,4.840533e-01) (30,3.751933e-01) (40,2.554489e-01) (50,1.499203e-01) (60,6.805965e-02) 
(70,2.723095e-02) (80,1.315459e-02) (90,8.841997e-03) (100,5.932500e-03) (110,4.521396e-03) (120,3.687442e-03) 
(130,2.998419e-03) (140,2.444599e-03) (150,2.043978e-03) (160,1.688579e-03) (170,1.352397e-03) (180,1.056196e-03) 
(190,8.500134e-04)
};
\addplot coordinates {
(10,5.893535e-01) (20,4.901633e-01) (30,3.826921e-01) (40,2.659440e-01) (50,1.664967e-01) (60,9.706291e-02)
(70,6.773150e-02) (80,4.630002e-02) (90,3.241731e-02) (100,2.199734e-02) (110,1.785652e-02) (120,1.420909e-02)
(130,1.144187e-02) (140,9.380460e-03) (150,7.792746e-03) (160,6.593513e-03) (170,5.464057e-03) (180,4.409625e-03)
(190,3.461842e-03) (200,2.731477e-03) (210,2.130924e-03) (220,1.678460e-03) (230,1.312523e-03) (240,9.796354e-04)
};
\end{semilogyaxis}
\end{tikzpicture}
\caption{$\frac{\|\tilde{\Ab}_{RC}-\Ab\|}{\|\Ab\|}$ at $k/(2\pi)=20$}\label{fig:error_inv_sd_ft_6}
\end{subfigure}
 \caption{(Experiment I.2 -- Supplemental results) Plots of $\|\tilde{\Ab}^{-1}_{RC}-\Ab^{-1}\|/\|\Ab^{-1}\|$ and $\|\tilde{\Ab}_{RC}-\Ab\|/\|\Ab\|$ for different $N_\lambda$ at wavenumbers: (a) and (b) $k/(2\pi)=5$, (c) and (d) $k/(2\pi)=10$, and (e) and (f) $k/(2\pi)=20$. Each line represents the number of iterations obtained when the preconditioner is obtained using $N_s=4$, $16$, $25$, $36$ and $64$. The overlap parameter is $\delta=8$.}\label{fig:error_inv_sd_ft}
\end{figure}


\begin{table}
\caption{(Experiment I.3 -- Supplemental results) Relative error of the approximation $\|\tilde{\Ab}^{-1}_{RC}-\Ab^{-1}\|/\|\Ab^{-1}\|$ and $\|\tilde{\Ab}_{RC}-\Ab\|/\|\Ab\|$ for domains composed of $N=32^2$, $64^2$ and $128^2$ scatterers at wavenumbers (a) $k/(2\pi)=5$ and (b) $k/(2\pi)=20$. At $k/(2\pi)=5$, the relative error of the matrix approximation is presented for $N_\lambda=40$, $90$, and $120$, while at $k=20(2\pi)$ the relative error of the matrix approximation is presented for $N_\lambda=20$, $40$, $60$, and $80$.}
\label{tab:error_tf}

\begin{subtable}{\textwidth}
\center
\caption{Relative error of the matrix approximation at $k/(2\pi)=5$}
\begin{tabular}{|c||c|c|c|c||c|c|c|c|}\cline{2-9}\cline{2-9}
\multicolumn{1}{c|}{} & \multicolumn{4}{c||}{$\|\tilde{\Ab}_{RC}^{-1}-\Ab^{-1}\|/\|\Ab^{-1}\|$}& \multicolumn{4}{c|}{$\|\tilde{\Ab}_{RC}-\Ab\|/\|\Ab\|$}\\\hline
\backslashbox{$N$}{$N_\lambda$} & 20 & 40 & 60 & 80 & 20 & 40 & 60 & 80\\\hline\hline
$32^2$  					      & 1.8e-1 & 5.8e-2 & 1.5e-2 & 2.3e-3 & 7.6e-2 & 1.8e-2 & 4.3e-3 & 7.7e-4 \\\hline
$64^2$  					      & 1.4e-1 & 4.6e-2 & 1.5e-2 & 8.6e-3 & 5.7e-2 & 8.3e-3 & 3.0e-3 & 1.4e-3 \\\hline
$128^2$					      & 2.9e-1 & 1.3e-1 & 6.2e-2 & 3.5e-2 & 8.5e-2 & 2.6e-2 & 1.1e-2 & 5.8e-3 \\\hline
\end{tabular}
\end{subtable}
\begin{subtable}{\textwidth}
\center
\caption{Relative error of the matrix approximation at $k/(2\pi)=20$}
\begin{tabular}{|c||c|c|c||c|c|c|}\cline{2-7}\cline{2-7}
\multicolumn{1}{c|}{} & \multicolumn{3}{c||}{$\|\tilde{\Ab}^{-1}_RC-\Ab^{-1}\|\|/\Ab^{-1}\|$}& \multicolumn{3}{c|}{$\|\tilde{\Ab}_{RC}-\Ab\|/\|\Ab\|$}\\\hline
\backslashbox{$N$}{$N_\lambda$}   & 40 & 90 & 120 & 40 & 90 & 120\\\hline\hline
$32^2$  						& 5.8e-2 & 6.3e-4 & 1.2e-4 & 1.8e-2 & 2.0e-4 & 6.6e-7 \\\hline
$64^2$					  	& 2.0e0  & 2.3e-2 & 1.3e-2 & 2.6e-1 & 6.8e-3 & 3.7e-3  \\\hline
$128^2$						& 5.6e-1 & 6.3e-2 & 4.6e-2 & 2.5e-1 & 2.0e-2 & 1.1e-2  \\\hline
\end{tabular}
\end{subtable}
\end{table}


\begin{table}
\caption{(Experiment I.3 -- Supplemental results) Relative error of the iterative solution of Equation \eqref{eq:Hessian_prob_1} using the HRC preconditioner with GMRES at (a) $k/(2\pi)=5$ and (b) $k/(2\pi)=20$. This table shows the relative error obtained for the number of iterations in Table \ref{tab:inv_scal_it}. The relative error have order of magnitude $\mathcal{O}(\Phi)$, with $\Phi=10^{-2}$, $10^{-3}$, $10^{-4}$ and $10^{-5}$. The total number of scatterer points in the domain are $N=32^2$, $64^2$ and $128^2$. The number of subdomains is fixed $N_s=16$ and the overlap parameter is $\delta=3$, $6$ and $12$ respectively for the domain with $N=32^2$, $64^2$ and $128^2$. The number of singular values used is $N_\lambda=20$, $40$, $60$ and $80$ at $k/(2\pi)=5$ and $N_\lambda=40$, $90$ and $120$ at $k/(2\pi)=20$.} \label{tab:inv_scal_error}
\begin{subtable}{\textwidth}
\center
\caption{Relative error for $k/(2\pi)=5$}\label{tab:inv_scal_error_20}
\begin{tabular}{|c||c||*{6}{c|}}\cline{4-7}\cline{4-7}
\multicolumn{3}{c|}{} & \multicolumn{4}{c|}{$N_\lambda$}\\\hline
$N$ & $\mathcal{O}(Error)$ & GMRES & 20 & 40 & 60 & 80 \\\hline\hline
\multirow{4}{*}{$32^2$}   & $10^{-2}$ & 1.8e-2 & 5.8e-2 & 5.0e-2 & 1.3e-2 & 6.4e-2 \\
                                        & $10^{-3}$ & 2.4e-3 & 4.7e-3 & 4.2e-3 &  2.6e-3  & 3.7e-4 \\
                                        & $10^{-4}$ & 2.2e-4 & 5.9e-4 & 6.2e-4 &  3.3e-4  & 3.7e-4 \\
                                        & $10^{-5}$ & 1.9e-5 & 7.4e-5 & 6.6e-5 &  3.5e-5  & 6.5e-5 \\\hline
\multirow{4}{*}{$64^2$}   & $10^{-2}$ & 1.9e-2 & 3.0e-2 & 3.6e-2 & 1.3e-2  & 1.3e-2 \\
                                        & $10^{-3}$ & 2.8e-3 & 4.9e-3 & 5.4e-3 & 3.4e-3  & 3.3e-4 \\
                                        & $10^{-4}$ & 2.4e-4 & 1.8e-4 & 4.4e-4 & 4.9e-4 & 3.3e-4\\
                                        & $10^{-5}$ & 3.0e-5 & 2.1e-5 & 8.8e-5 & 2.9e-5 & 1.5e-5\\\hline
\multirow{4}{*}{$128^2$} & $10^{-2}$ & 1.8e-2 & 6.9e-2 & 1.8e-2  & 3.7e-2 & 2.4e-2\\
                                        & $10^{-3}$ & 3.8e-3 & 3.8e-3 & 1.5e-3  & 1.8e-3 & 2.0e-3\\
                                        & $10^{-4}$ & 2.6e-4 & 4.4e-4 & 6.5e-4  & 1.1e-4 & 2.9e-4\\
                                        & $10^{-5}$ & 3.2e-5 & 5.1e-5 & 8.2e-5  & 1.1e-5 & 2.3e-5\\\hline
\end{tabular}
\end{subtable}
\begin{subtable}{\textwidth}
\center
\caption{Relative error for $k/(2\pi)=20$}\label{tab:inv_scal_error_20}
\begin{tabular}{|c||c||*{5}{c|}}\cline{4-6}\cline{4-6}
\multicolumn{3}{c|}{} & \multicolumn{3}{c|}{$N_\lambda$}\\\hline
$N$ & $\mathcal{O}(Error)$ & GMRES & 40 & 90 & 120  \\\hline\hline
\multirow{4}{*}{$32^2$}   & $10^{-2}$ & 4.0e-2 & 3.2e-2 &  1.2e-2  & 2.2e-2 \\
                                        & $10^{-3}$ & 7.6e-3 & 3.1e-3 &  1.9e-3  & 1.7e-3 \\
                                        & $10^{-4}$ & 8.9e-4 & 3.0e-4 &  1.4e-4  & 1.5e-5 \\
                                        & $10^{-5}$ & 8.9e-5 & 1.8e-5 &  2.2e-5  & 1.5e-5 \\\hline
\multirow{4}{*}{$64^2$}   & $10^{-2}$ & 4.1e-2 & 1.6e-2 &  3.8e-2  & 3.7e-2 \\
                                        & $10^{-3}$ & 5.8e-3 & 1.2e-3 &  3.0e-3  & 3.4e-3 \\
                                        & $10^{-4}$ & 8.1e-4 & 1.3e-4 &  2.3e-4  & 3.6e-4 \\
                                        & $10^{-5}$ & 8.5e-5 & 4.6e-5 &  4.2e-5  & 9.3e-5 \\\hline
\multirow{4}{*}{$128^2$} & $10^{-2}$ & 2.7e-2 & 3.8e-2 &  1.8e-2  & 1.6e-2 \\
                                        & $10^{-3}$ & 5.8e-3 & 2.2e-3 &  1.4e-3  & 1.1e-3  \\
                                        & $10^{-4}$ & 9.0e-4 & 2.3e-4 &  1.2e-4  & 4.4e-5  \\
                                        & $10^{-5}$ & 8.2e-5 & 2.1e-5 &  1.1e-5  & 4.4e-5  \\\hline

\end{tabular}
\end{subtable}
\end{table}


\begin{figure}
  \centering
  \begin{subfigure}[t]{0.48\textwidth}
  \input{H_k5_example4.tex}
\caption{Singular Values for $\Hb$ at $k/(2\pi)=5$}
\end{subfigure}
\begin{subfigure}[t]{0.48\textwidth}
  \input{LR_k5_lambda140_example4.tex}
\caption{LR at $k=5/(2\pi)$ and $N_\lambda=140$.}\label{fig:Heig_example_I4:a}
\end{subfigure}

\begin{subfigure}[t]{0.48\textwidth}
\input{HRC_k5_Ns16_delta4_lambda140_example4.tex}
\caption{HRC-4 with $k=5/(2\pi)$, $\delta=4$ and $N_\lambda=140$.}\label{fig:Heig_example_I4:c}
\end{subfigure}
\begin{subfigure}[t]{0.48\textwidth}
\input{HRC_k5_Ns16_delta8_lambda40_example4.tex}
\caption{HRC-8 with $k=5/(2\pi)$, $\delta=8$ and $N_\lambda=40$.}\label{fig:Heig_example_I4:e}
\end{subfigure}
\caption{(Experiment I.4 -- Supplemental results) Plot of the singular values of (a) $\Hb$, (b) $\tilde{\Hb}_{LR}$ with $N_\lambda=140$, (c) $\tilde{\Hb}_{HRC}$ with $N_\lambda=140$, $N_s=16$ and $\delta=4$, and (d) $\tilde{\Hb}_{HRC}$ with $N_\lambda=40$, $N_s=16$ and $\delta=8$ when the incoming incident waves have wavenumber $k/(2\pi)=5$.}
\label{fig:Heig_example_I4_k5}
\end{figure}

\begin{figure}
  \centering
\begin{subfigure}[t]{0.48\textwidth}
  \input{H_k20_example4.tex}
\caption{Singular values for $\Hb$ at $k/(2\pi)=20$}
\end{subfigure}
\begin{subfigure}[t]{0.48\textwidth}
  \input{LR_k20_lambda240_example4.tex}
\caption{LR at $k/(2\pi)=20$ and $N_\lambda=240$.}\label{fig:Heig_example_I4:b}
\end{subfigure}

\begin{subfigure}[t]{0.48\textwidth}
\input{HRC_k20_Ns16_delta4_lambda240_example4.tex}
\caption{HRC-4 with $k/(2\pi)=20$, $\delta=4$ and $N_\lambda=240$.}\label{fig:Heig_example_I4:d}
\end{subfigure}
\begin{subfigure}[t]{0.48\textwidth}
\input{HRC_k20_Ns16_delta8_lambda100_example4.tex}
\caption{HRC-8 with $k/(2\pi)=20$, $\delta=8$ and $N_\lambda=100$.}\label{fig:Heig_example_I4:f}
\end{subfigure}
\caption{(Experiment I.4 -- Supplemental results) Plot of the singular values of (a) $\Hb$, (b) $\tilde{\Hb}_{LR}$ with $N_\lambda=240$, (c) $\tilde{\Hb}_{HRC}$ with $N_\lambda=240$, $N_s=16$ and $\delta=4$, and (d) $\tilde{\Hb}_{HRC}$ with $N_\lambda=100$, $N_s=16$ and $\delta=8$ when the incoming incident waves have wavenumber $k/(2\pi)=20$.}
\label{fig:Heig_example_I4_k20}
\end{figure}